\documentclass[leqno,11pt]{article}
\usepackage{latexsym}
\usepackage{amssymb}
\usepackage{amsfonts}
\usepackage{amsmath}
\usepackage{color}
\usepackage{epsfig, graphics, hyperref}

\makeatletter
\long\def\unmarkedfootnote#1{{\long\def\@makefntext##1{##1}\footnotetext{#1}}}
\makeatother

\newtheorem{definition}{Definition}[section]

\newtheorem{lemma}[definition]{Lemma}

\newtheorem{theorem}[definition]{Theorem}

\newtheorem{proposition}[definition]{Proposition}

\newtheorem{remark}[definition]{Remark}

\newtheorem{example}[definition]{Example}

\def\o{\Omega}

\def\m2{|\o| /2}
\def\M2{\frac{\hn (\Omega )}{2}}
\def\u+{u_+^*}

\def\-p{\overline{p}}

\def\w0{{W_0^{m,A}(\Omega)}}
\def\i{{(0, |\Omega|)}}

\def\R{\mathbb R}
\def\N{\mathbb N}
\def\Z{\mathbb Z}
\def\S{\mathbb S}
\def\rn{{{\R}^n}}
\def\sn{{{\S}^{n-1}}}

\newcommand{\hh}{{\mathcal H}^{n-1}}
\newcommand{\hi}{{\mathcal H}^{n-1}_\infty}
\newcommand{\hn}{{\mathcal H}^{n}}
\newcommand{\medint}{-\kern  -,395cm\int}
\newcommand{\medintinrigo}{-\kern  -,310cm\int}
\newcommand{\medelle}{-\kern  -,235cm L}
\newcommand{\medellenrigo}{-\kern  -,180cm L}
\newcommand{\qed}{\thinspace\null\nobreak\hfill
\hbox{\vbox{\kern-.2pt\hrule height.2pt
depth.2pt\kern-.2pt\kern-.2pt \hbox to1.8mm {\kern-.2pt\vrule
width.4pt \kern-.2pt\raise1.8mm\vbox to.2pt{} \lower0pt\vtop
to.2pt{}\hfil\kern-.2pt \vrule
width.4pt\kern-.2pt}\kern-.2pt\kern-.2pt \hrule height.2pt
depth.2pt \kern-.2pt}}\par\medbreak}

\title{
Sobolev inequalities in arbitrary domains}
\numberwithin{equation}{section}
\author{
  Andrea Cianchi\\
 {\it Dipartimento di Matematica e Informatica \lq\lq U.Dini", Universit\`a di Firenze}\\ {\it Piazza Ghiberti
27, 50122 Firenze, Italy
}
\bigskip
 \and
 Vladimir  Maz'ya \\
{\it   Department of Mathematics, Link\"oping University, SE-581 83
Link\"oping, Sweden}
\\ and \\
{\it  Department of Mathematical Sciences, M\&O Building}\\ {\it
University of Liverpool, Liverpool L69 3BX,
 UK}}
\begin{document}
\maketitle
%

\begin{abstract} A theory of Sobolev inequalities in arbitrary open
sets in $\rn$ is established. Boundary regularity of domains is
replaced with information on boundary traces of trial functions and
of their derivatives up to some explicit minimal order. The relevant
Sobolev inequalities involve constants independent of the geometry
of the domain, and exhibit the same critical exponents as in the
classical inequalities on regular domains. Our approach relies upon
 new representation formulas for Sobolev functions, and on ensuing pointwise
estimates which hold in any open set.
\end{abstract}

\unmarkedfootnote {
\par {\it
 This research was partly supported by MIUR
(Italian Ministry of Education, University and Research) via  the
research project ``Geometric aspects of partial differential
equations and related topics" 2008, and by  GNAMPA of the Italian
INdAM (National Institute of High Mathematics).}}

\section{Introduction }\label{sec1}

$\mathcal a$ $\mathcal b$

The aim of this paper is to develop a theory of Sobolev embeddings,
of any order $m \in \N$, in \emph{arbitrary}
 open sets $\Omega$ in $\rn$. As usual, by an $m$-th order
Sobolev embedding we mean an inequality between a norm of the $h$-th
order weak derivatives ($0 \leq h \leq m-1$) of  any $m$ times
weakly differentiable  function in $\Omega$ in terms of  norms of
some of its derivatives up to the order $m$.
\par
The classical theory of Sobolev embeddings involves ground domains
$\Omega$ satisfying suitable regularity assumptions. For instance, a
formulation of the original theorem by Sobolev reads as follows.
Assume that $\Omega$ is a bounded domain satisfying the cone
property,  $m \in \N$,   $1 < p <\tfrac nm$, and $ \mathcal F (\cdot
)$  is any continuous seminorm in $W^{m,p}(\Omega)$ which does not
vanish on any polynomial of degree not exceeding $m-1$.  Then there
exists a constant $C=C(\Omega)$ such that
\begin{equation}\label{sobolevm}
\|u\|_{L^{\frac {np}{n-mp}}(\Omega)} \leq C \big(\|\nabla ^m
u\|_{L^p(\Omega )} + \mathcal F (u)\big)
\end{equation}
for every  $u\in W^{m,p}(\Omega)$. Here,  $W^{m,p}(\Omega)$ denotes
the usual  Sobolev space of those functions in $\Omega$ whose weak
derivatives up to the order $m$ belong to $L^p(\Omega)$, and $\nabla
^m u$ stands for the vector of all (weak) derivatives of $u$ of
order $m$.
\par It is well known that standard Sobolev
embeddings are spoiled in presence of domains with \lq\lq bad"
boundaries. In particular, inequalities of the form \eqref{sobolevm}
do not hold, at least with the same critical exponent $\frac
{np}{n-mp}$, in irregular domains. This is the case, for instance,
of  domains with outward cusps. A theory of Sobolev embeddings,
including possibly irregular domains, was initiated in the papers
\cite{Ma1960} and \cite{Ma1961}, and is systematically exposed in
the monograph \cite{Mabook}, where   classes of Sobolev inequalities
are characterized in terms of geometric properties of the domain.
Specifically, they are shown to be equivalent to either
isoperimetric or isocapacitary inequalities relative to the domain.
The interplay between the geometry of the domain and Sobolev
inequalities, even in  frameworks more general than the Euclidean
one, has over the years been the subject of extensive
investigations, along diverse directions, by a number of authors.
Their results are the object of a rich literature,
 which includes the papers
  \cite{AFT, Aubin, BCR1, BL,
BWW, BH, BLbis, BK, BK1, Cheeger, Ci_ind, Ci1,  CFMP1, CP_gauss,
EKP, EFKNT, Gr, HK, HS1, KP, KM, Kl, Ko, LPT, LYZ, M, Moser, Ta, Zh}
and the monographs \cite{BZ, CDPT, chavel, Heb, Mabook, Saloff}. An
updated bibliography on the area of Sobolev type inequalities can be
found in \cite{Mabook}.

\par
In order to remove any a priori regularity assumption on $\Omega$,
we consider Sobolev inequalities from an unconventional perspective.
The underling idea of our results is that suitable information on
boundary traces of trial functions
can replace boundary regularity of the domain
 in
Sobolev inequalities.

\par The inequalities that will be established have the form
\begin{equation}\label{cianchimazyah}
\|\nabla ^h u\|_{Y(\Omega , \mu)} \leq C \big(\|\nabla ^m
u\|_{X(\Omega )} +  \mathcal N _{\partial \Omega} (u)\big),
%
%
\end{equation}
where
 $m \in \N$, $h \in \N _0$,
$\|\cdot\|_{X(\Omega )}$ is a Banach function norm on $\Omega$ with
respect to  Lebesgue measure $\mathcal L^n$, $\|\cdot\|_{Y(\Omega ,
\mu)}$ is a Banach function norm with respect to a possibly more
general measure $\mu$, and $\mathcal N _{\partial \Omega} (\cdot)$
is a (non-standard) seminorm on $\partial \Omega$, depending on the
trace of $u$, and of its derivatives up to
 the order $\big[\tfrac {m-1}2\big]$. Here, $\N _0 = \N \cup \{0\}$,
 and
$[\cdot]$ stands for integer part. Moreover, $\nabla ^0 u$ stands
just for $u$, and we shall denote $\nabla ^1u$ also by $\nabla u$.

\smallskip
\par Some distinctive features of the inequalities
 to be presented  can be itemized as follows:

\smallskip
 \par\noindent
$\bullet$ No regularity on $\Omega$ is a priori assumed. In
particular, the constants in \eqref{cianchimazyah} are independent
of the geometry of $\Omega$.

\smallskip
\par\noindent
$\bullet$ The  critical Sobolev exponents,
or, more generally, the optimal target norms, are the same as in the
case of regular domains.

\smallskip
\par\noindent
$\bullet$ The order $\big[\tfrac {m-1}2\big]$ of the derivatives, on
which the seminorm $\mathcal N _{\partial \Omega} (\cdot)$ depends,
is minimal for an inequality of the form \eqref{cianchimazyah} to
hold without any additional assumption on $\Omega$.

\smallskip
\par
A first-order Sobolev inequality on arbitrary domains $\Omega$ in
$\rn$ of the form \eqref{cianchimazyah}, where $X(\Omega ) =
L^p(\Omega)$, $Y(\Omega , \mu)= L^q(\Omega)$, and  $\mathcal N
_{\partial \Omega} (\cdot)= \|\cdot\|_{L^r(\partial \Omega)}$, with
$1 \leq p < n$, $r\geq 1$ and $q = \min\{\tfrac{rn}{n-1},
\tfrac{np}{n-p}\}$ was established in
 \cite{Ma1960} via isoperimetric inequalities. Sobolev inequalities of this kind,
  but still  involving only first-order derivatives and  Lebesgue measure, have received a renewed
attention in recent years. In particular, the paper
\cite{MaggiVillani1}  makes use of mass transportation techniques to
address the problem of the optimal constants  for $p \in (1,n)$, the
problem when $p=1$ having already been solved in \cite{Ma1960}.
Sharp constants in  inequalities  in the borderline case when $p=n$
are exhibited in \cite{MaggiVillani2}.

\par
In the present paper, we develop a completely different approach,
which not only enables us to establish arbitrary-order inequalities,
which cannot just be derived via iteration of first-order ones, but
also augments the first-order theory, in that more general measures
and norms are allowed.
\par Our point of departure  is a
new pointwise estimate for  functions, and their derivatives, on
arbitrary -- possibly unbounded and with infinite measure -- domains
$\Omega$.  Such estimate involves  a novel class of
 double-integral operators, where integration is extended over $\Omega \times \mathbb
S^{n-1}$. The relevant operators act on a kind of higher-order
difference quotients of the traces of functions and of their
derivatives on $\partial \Omega$.
 \\ In view of applications to norm inequalities, the next
step   calls for an analysis of boundedness properties of these
operators in function spaces. To this purpose, we prove their
boundedness  between optimal endpoint spaces. In combination with
interpolation arguments based on the use of Peetre $K$-functional,
these endpoint results lead to pointwise  bounds, for Sobolev
functions, in rearrangement form. As a consequence, Sobolev
inequalities on an arbitrary $n$-dimensional domain are reduced to
considerably simpler one-dimensional inequalities for Hardy type
operators.
\par With this apparatus at disposal, we are able to
establish inequalities involving Lebesgue norms, with respect to
quite general measures, as well as Yudovich-Pohozaev-Trudinger type
inequalities in exponential Orlicz
 spaces for limiting situations. The compactness of  corresponding
 Reillich-Kondrashov type embeddings, with subcritical exponents, is also shown.
 Inequalities for other rearrangement-invariant
norms, such as Lorentz and Orlicz norms, could be derived. However,
in order to avoid unnecessary additional technical complications,
this issue is not   addressed here.
\par
The paper is organized as follows. In the next section we offer a
brief overview of some  Sobolev type inequalities, in basic cases,
which follow from our results, and discuss their novelty and
optimality. Section \ref{notation} contains some preliminary
definitions and results.
 The statement
of our main results starts with Section \ref{proofs},  which is
devoted to our key pointwise inequalities for Sobolev functions on
arbitrary open sets. Estimates in rearrangement form are derived in
the subsequent Section \ref{sobolev}. In Section \ref{ineq}, Sobolev
type inequalities in arbitrary open sets are shown to follow via
such estimates. Examples which demonstrate the sharpness of our
results are exhibited in Section \ref{sharp}. In particular,
Example \ref{ex4} shows that inequalities of the form
\eqref{cianchimazyah} may possibly
 fail if $\mathcal N_{\partial \Omega}(u)$ only depends on derivatives of $u$ on $\partial \Omega$ up  to an order smaller than $[\frac {m
-1}2]$.
 Finally, in the
Appendix,  some new notions,
which are introduced in the definitions of the seminorms $\mathcal N
_{\partial \Omega}(\cdot)$,
are linked to classical properties of Sobolev functions.

\section{A taste of results}\label{over}

In order to give an overall idea
of the content of this paper, we enucleate hereafter a few basic
instances
of the inequalities that can be derived via our approach.
\par We begin with two examples which demonstrate
that our conclusions  lead to new results also in the case of
first-order inequalities, namely in the case when $m=1$ in
\eqref{cianchimazyah}.
\\
Let $\Omega$ be any open set in $\rn$, and let $\mu$ be a Borel
measure on $\Omega$ such that $\mu (B_r \cap \Omega) \leq Cr^\alpha$
for some  $C>0$, and $\alpha \in (n-1, n]$,  and for every ball
$B_r$ radius $r$. Clearly, if $\mu = \mathcal L^n$, then this
condition holds with $\alpha =n$.
\\ Assume  that $1<p<n$ and
$r>1$, and let $s=\min\{\tfrac{r\alpha}{n-1}, \tfrac{\alpha
p}{n-p}\}$. Then
\begin{align}\label{firstmeas}
\|u\|_{L^s(\Omega, \mu)}\leq C \big(\|\nabla  u\|_{L^p(\Omega )} +
\|u\|_{L^r(\partial \Omega )}\big)
\end{align}
for some constant $C$ and every function $u$ with bounded support,
provided that $\mathcal L^n (\Omega)< \infty$, $\mu (\Omega)<
\infty$ and $\hh (\partial \Omega) < \infty$. Here, $\hh$ denotes
the $(n-1)$-dimensional Hausdorff measure. In particular, if $r=
\tfrac{p(n-1)}{n-p}$, and hence $s= \tfrac {\alpha p}{n-p}$, then
\eqref{firstmeas}  holds even if the assumption on the finiteness of
these measures is dropped; in this case, the constant $C$ depends
only on $n$. Inequality \eqref{firstmeas} follows via a general
principle contained in Theorem \ref{reduction}, Section \ref{ineq}.
It extends a version  of the  Sobolev inequality for measures, on
regular domains \cite[Theorem 1.4.5]{Mabook}. It also augments, at
least for $p>1$, the results for general domains of \cite{Ma1960}
and \cite{MaggiVillani1}, whose approach is confined to norms
evaluated with respect to the Lebesgue measure. Let us point out
that, by contrast, our method, being based on representation
formulas, need not lead to optimal inequalities for $p=1$.
\\
Consider now the borderline case corresponding to  $p=n$. As a
consequence of Theorem \ref{reduction} again, one can show that
\begin{align}\label{trudfirst}
\| u\|_{\exp L^{\frac {n}{n-1}}(\Omega, \mu)}  \leq C\Big( \|\nabla
 u\|_{L^{n}(\Omega )} +
 \|u\|_{\exp L^{\frac
{n}{n-1}}(\partial \Omega)}\Big),
\end{align}
for some constant $C$ and every function $u$ with bounded support,
provided that $\mathcal L^n (\Omega)< \infty$, $\mu (\Omega)<
\infty$ and $\hh (\partial \Omega) < \infty$. Here, $\|
\cdot\|_{\exp L^{\frac {n}{n-1}}(\Omega, \mu)}$ and $
\|\cdot\|_{\exp L^{\frac {n}{n-1}}(\partial \Omega )}$ denote norms
in Orlicz spaces of exponential type on $\Omega$ and $\partial
\Omega$, respectively.
 Inequality
\eqref{trudfirst} on the one hand extends  the
 Yudovich-Pohozaev-Trudinger  inequality to possibly irregular
 domains;  on the other hand, it improves a result of
 \cite{MaggiVillani2}, where estimates for the weaker norm in  $\exp L(\Omega)$
 are established, and just for the Lebesgue measure.

\smallskip
\par
Let us now turn to higher-order inequalities. Focusing, for the time
being, on second-order inequalities
may help to grasp the quality and sharpness
of our conclusions in this framework. In the remaining part of this
section, we thus assume that $m =2$ in \eqref{cianchimazyah}; we
also assume, for simplicity, that $\mu = \mathcal L^n$.
\par
 First, assume that $h=0$. Then
 we can prove (among other possible choices of the exponents) that, if $1<p< \tfrac n2$,  then
\begin{align}\label{mainsecond0}
\| u\|_{L^{\frac{pn}{n-2p}}(\Omega)}  \leq C \big(\|\nabla ^{2}
u\|_{L^p(\Omega )} +   \|u\|_{\mathcal
V^{1,0}L^{\frac{p(n-1)}{n-p}}(\partial \Omega )}  +
\|u\|_{L^{\frac{p(n-1)}{n-2p}}(\partial \Omega )}\big),
\end{align}
 for some constant $C=C(p,n)$, in particular independent of
 $\Omega$, and every function $u$ with bounded support. Note that $\tfrac{pn}{n-2p}$ is the same critical Sobolev
 exponent as in the case of regular domains.
 Here, $\|\cdot\|_{\mathcal V^{1,0}L^{r}(\partial \Omega
 )}$ denotes, for $r \in [1,\infty]$, the seminorm given by
\begin{equation}\label{hajlasznorm}
\|u\|_{\mathcal V^{1,0}L^{r}(\partial \Omega )} = \inf
_{g}\|g\|_{L^r (\partial \Omega)},
\end{equation}
where the infimum is taken among all Borel functions $g$ on
$\partial \Omega$ such that
\begin{equation}\label{hajlasz}
|u(x) - u(y)| \leq |x -y| (g (x) + g   (y)) \quad \hbox{for
$\hh$-a.e. $x, y \in
\partial \Omega$,}
\end{equation}
and $L^{r}(\partial \Omega )$ denotes a Lebesgue space on $\partial
\Omega$ with respect to the  measure $\hh$.
%
The function $g$ appearing in \eqref{hajlasz}  is an upper gradient,
in the sense of \cite{Hajlasz}, for the restriction of $u$ to
$\partial \Omega$, endowed with the metric inherited from the
Euclidean metric in $\rn$, and with the measure $\hh$. In
\cite{Hajlasz}, a definition of this kind, and an associated
seminorm given as in \eqref{hajlasznorm}, were introduced to define
first-order Sobolev type spaces on arbitrary metric measure spaces.
 In the last two decades, various notions of upper gradients and of Sobolev spaces of
functions defined on metric measure spaces, have been the object of
 investigations and applications. They constitute the topic of a number of papers and
 monographs, including  \cite{AT, BB, FHK, HK, Hein, HeKo, Koskela}.
\\ Let us emphasize that, although the new term $ \|u\|_{\mathcal V^{1,0}L^{\frac{p(n-1)}{n-p}}(\partial \Omega )}$
on the right-hand side of \eqref{mainsecond0} can be dropped  when
$\Omega$ is a regular, say Lipschitz, domain,
it is indispensable in an arbitrary domain. This can be shown by a
domain as in Figure 1 (see Example \ref{ex3}, Section \ref{sharp}).

\begin{figure}
\begin{center}
\includegraphics[height=10cm]{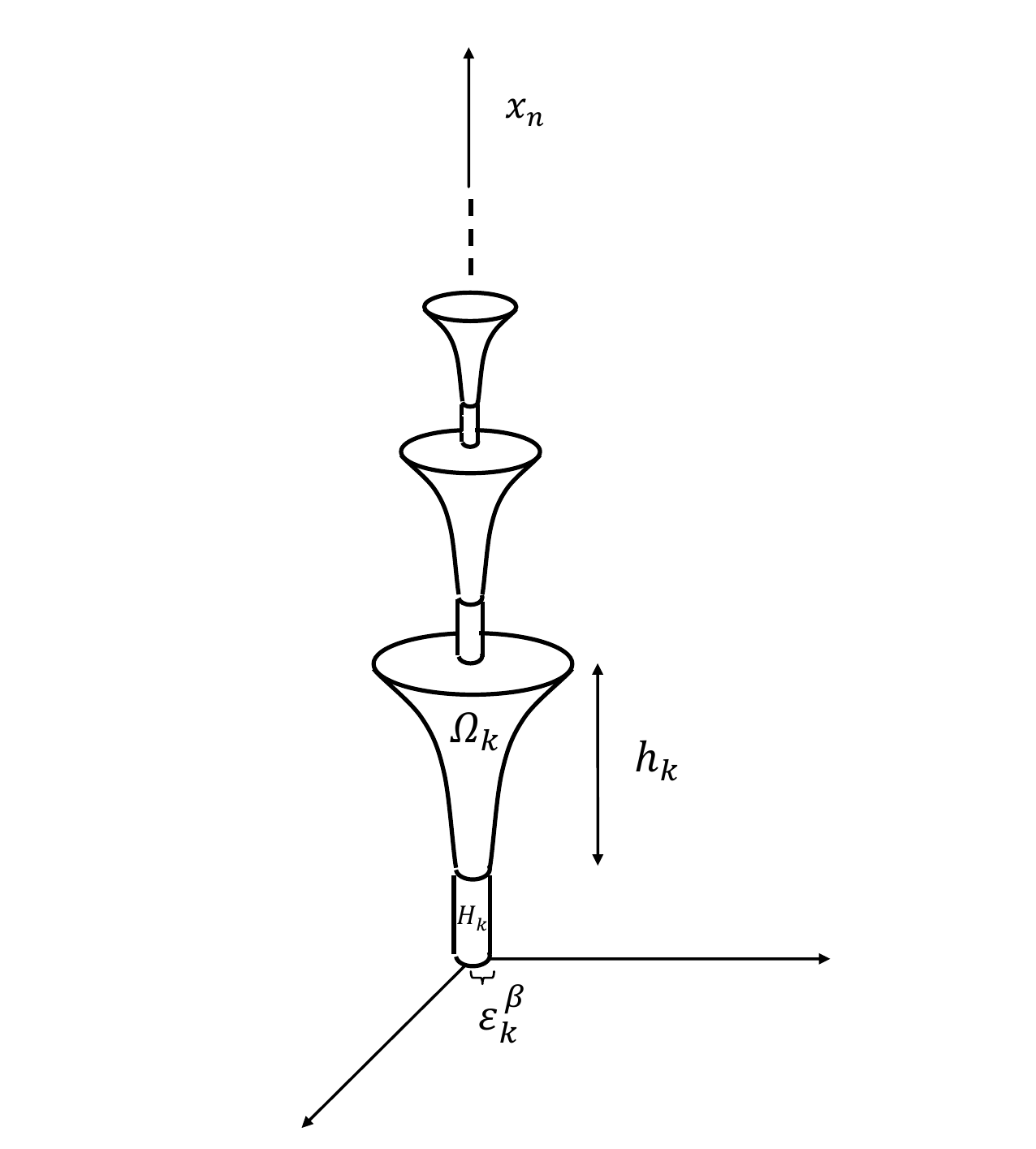}
\end{center}
        \label{Fig3}
        \caption{Example \ref{ex3}, Section \ref{sharp}}
\end{figure}


\par As in the case of regular domains,  if $p >\tfrac n2$,
 then the Lebesgue norm on the left-hand side of \eqref{mainsecond0} can
be replaced by the norm in $L^\infty$. Indeed, if  $r>n-1$, then
\begin{align}\label{secondinf}
\| u\|_{L^{\infty}(\Omega)}  \leq C \big(\|\nabla ^{2}
u\|_{L^p(\Omega )} +   \|u\|_{\mathcal V^{1,0}L^{r}(\partial \Omega
)} + \|u\|_{L^{\infty}(\partial \Omega )}\big)
\end{align}
 for any open set $\Omega$ such that $\mathcal
L^n (\Omega) < \infty$ and $\hh (\partial \Omega) < \infty$,  for
some constant $C$, and for any function $u$ with bounded support. In
particular, the constant $C$ depends on $\Omega$ only through
$\mathcal L^n (\Omega)$ and $\hh (\partial \Omega)$.
 \par
 In the limiting situation when $n \geq 3$, $p =\tfrac n2$ and $r>n-1$,
 a
 Yudovich-Pohozaev-Trudinger type
 inequality of the form
\begin{align}\label{secondexp}
\| u\|_{\exp L^{\frac {n}{n-2}}(\Omega)} \leq C \big(\|\nabla ^{2}
u\|_{L^{\frac n2}(\Omega )} +   \|u\|_{\mathcal
V^{1,0}L^{r}(\partial \Omega )} + \|u\|_{\exp L^{\frac
{n}{n-2}}(\partial \Omega )}\big)
\end{align}
holds for some constant $C$ independent of the regularity of
$\Omega$, and every function $u$ with bounded support, provided that
$\mathcal L^n (\Omega)< \infty$ and $\hh (\partial \Omega) <
\infty$. The norms $\| \cdot\|_{\exp L^{\frac {n}{n-2}}(\Omega)}$
and $ \|\cdot\|_{\exp L^{\frac {n}{n-2}}(\partial \Omega )}$
 are the same exponential
norms appearing in the Yudovich-Pohozaev-Trudinger inequality on
regular domains, and in its boundary trace counterpart.

%
%
%

\medskip
\par
Consider next the case when still $m=2$  in \eqref{cianchimazyah},
but $h=1$. Then one can infer from our estimates that, if $1<p < n$
and $r \geq 1$, and $\Omega$ is any open set with $\mathcal L ^n
(\Omega )< \infty$ and $\hh (\partial \Omega)< \infty$, then
\begin{align}\label{mainsecond1}
\|\nabla u\|_{L^{q}(\Omega)}  \leq C \big(\|\nabla ^{2}
u\|_{L^p(\Omega )}   +   \| u\|_{\mathcal V ^{1,0}L^{r}(\partial
\Omega )} \big)
\end{align}
for some constant $C$ independent of the geometry of $\Omega$, and
every function $u$ with bounded support, where
\begin{equation}\label{q}
q = \min\big\{\tfrac{rn}{n-1}, \tfrac{np}{n-p}\big\}.
\end{equation}
In particular, if $r=\tfrac{p(n-1)}{n-p}$, and hence
$q=\frac{np}{n-p}$, then the constant $C$ in \eqref{mainsecond1}
depends only on $n$ and $p$.
\\ Inequality \eqref{mainsecond1}
is optimal under various respects. For instance,  if
 $\Omega$ is  regular,  then, as a consequence of \eqref{sobolevm},
 the seminorm $\| u\|_{\mathcal V ^{1,0}L^{r}(\partial
\Omega )}$ can  be replaced just with $\| u\|_{L^{r}(\partial \Omega
)}$ on the right-hand
 side. By contrast, a domain $\Omega$ as in Figure 2
 shows that  this is impossible for every $q \in [1, \tfrac
{np}{n-p}]$, whatever $r$ is -- see Example
 \ref{ex1}, Section \ref{sharp}.

\begin{figure}
\begin{center}
\includegraphics[height=10cm]{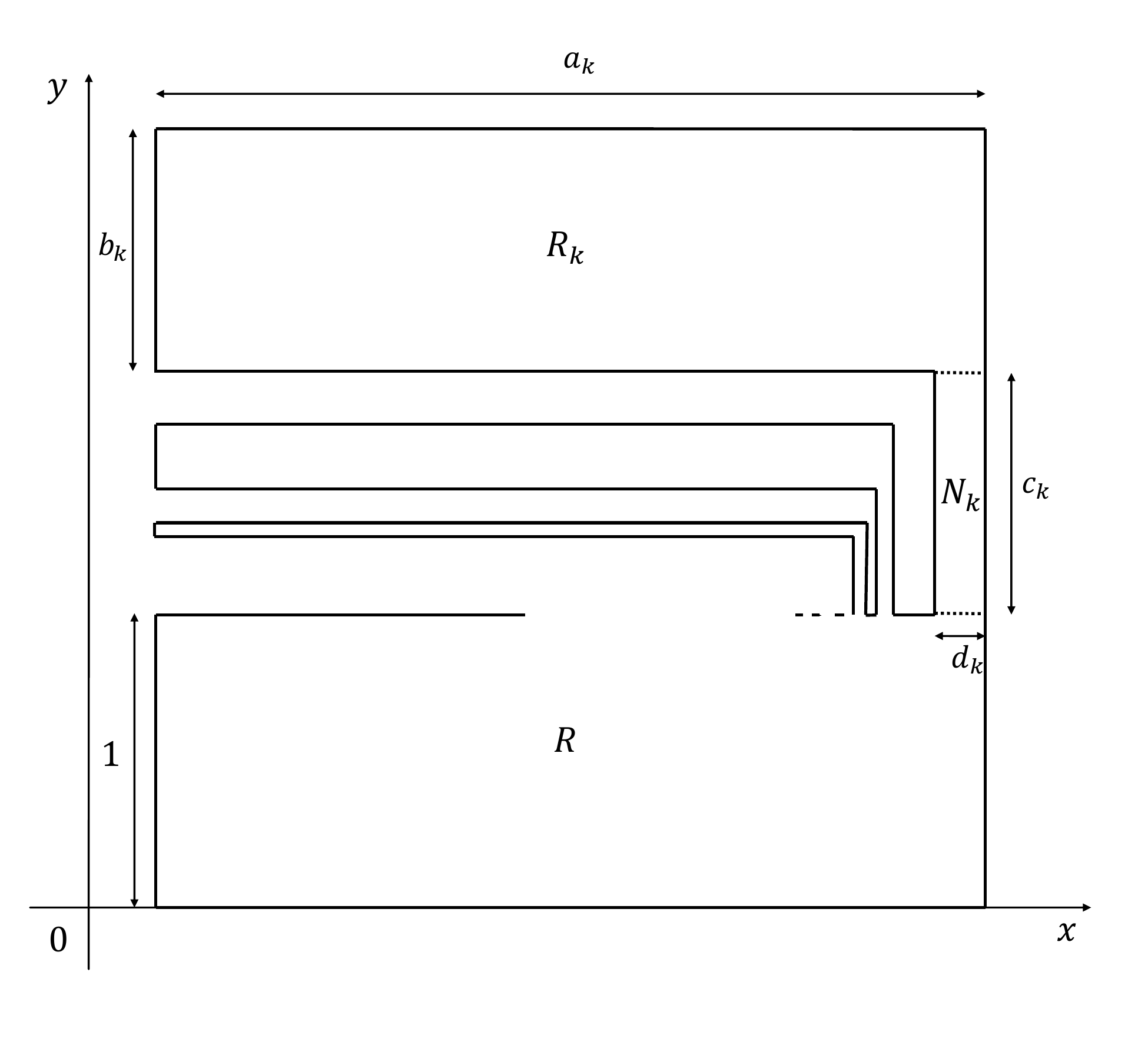}
\end{center}
        \label{Fig1}
        \caption{ Example \ref{ex1}, Section \ref{sharp}}
\end{figure}

\par\noindent
The question of the optimality of the exponent $q$ given by
\eqref{q} can also be raised. The answer is affirmative. Actually,
domains like that of Figure 3 show that such exponent $q$ is the
largest possible in \eqref{mainsecond1} if no regularity is imposed
on $\Omega$ (Example \ref{ex2}, Section \ref{sharp}).

\begin{figure}[ht]
\begin{center}
\includegraphics[height=10cm]{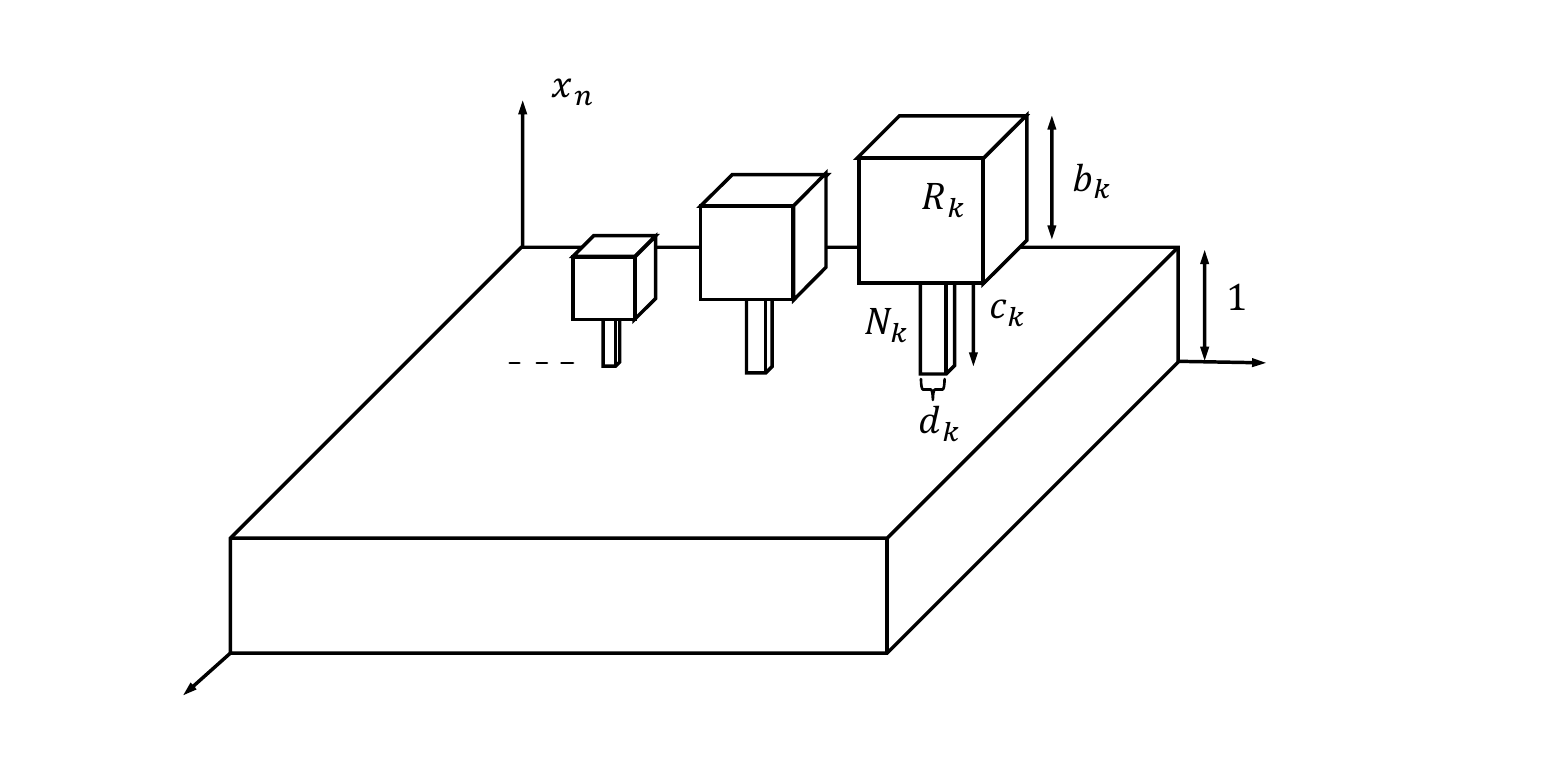}
\end{center}
        \label{Fig2}
        \caption{Example \ref{ex2}, Section \ref{sharp}}
\end{figure}

\par
When  $p > n$, inequality \eqref{mainsecond1} can be replaced with
\begin{align}\label{mainsecond1inf}
\|\nabla u\|_{L^{\infty}(\Omega)}  \leq C \big(\|\nabla ^{2}
u\|_{L^p(\Omega )}   +   \| u\|_{\mathcal V
^{1,0}L^{\infty}(\partial \Omega )} \big),
\end{align}
for some constant $C$ independent of the regularity of $\Omega$, and
every function $u$ with bounded support, provided that $\mathcal L
^n (\Omega ) < \infty$ and $\hh (\partial \Omega)< \infty$.
\par
Finally, in the borderline case corresponding to $p=n$, an
exponential norm is involved again. Under the assumption that
$\mathcal L ^n (\Omega ) < \infty$ and $\hh (\partial \Omega)<
\infty$, one has that
\begin{align}\label{mainsecond1esp}
\|\nabla  u\|_{\exp L^{\frac {n}{n-1}}(\Omega)}  \leq C
\big(\|\nabla ^{2} u\|_{L^{n}(\Omega )} +  \|u\|_{\mathcal V^{1,
0}\exp L^{\frac {n}{n-1}}(\partial \Omega)}\big)
\end{align}
for some constant $C$, depending on $\Omega$ only through $\mathcal
L ^n (\Omega )$ and $\hh (\partial \Omega)$,  and for every function
$u$ with bounded support. Here, the seminorm $\|\cdot\|_{\mathcal
V^{1, 0}\exp L^{\frac {n}{n-1}}(\partial \Omega)}$ is defined as in
\eqref{hajlasznorm}, with the norm $\|\cdot \|_{L^r(\partial
\Omega)}$ replaced with the norm  $\|\cdot \|_{\exp L^{\frac
{n}{n-1}}(\partial \Omega)}$. Again, the exponential norms in
\eqref{mainsecond1esp} are the same optimal Orlicz target norms for
Sobolev and trace inequalities, respectively, on regular domains.

\section{Preliminaries}\label{notation}

 Let $\Omega$ be any open  set in $\rn$, $n \geq 2$. Given $x \in
 \Omega$,
 define
\begin{equation}\label{omegax}
\Omega _x = \{y \in \Omega: (1-t)x + ty \subset \Omega\,\,\,
\hbox{for every $t \in (0,1)$}\},
\end{equation}
and
\begin{equation}\label{partomegax}
(\partial \Omega )_x = \{y \in \partial \Omega: (1-t)x + ty \subset
\Omega\,\,\, \hbox{for every $t \in (0,1)$}\}.
\end{equation}
They are the largest subset of $\Omega$ and $\partial \Omega $,
respectively, which can be \lq\lq seen" from $x$.
 It is easily verified that $\Omega _x$ is
an open set. The following proposition tells us that $(\partial
\Omega )_x$ is a Borel
%
set.

\begin{proposition}\label{borel}  Assume that $\Omega$ is an open set in $\rn$, $n \geq 2$.  Let $x \in \Omega$. Then the set
$(\partial \Omega) _x$, defined by \eqref{partomegax}, is  Borel measurable.
\end{proposition}
{\bf Proof}. Given any $r \in \mathbb Q \cap (0,1)$, define
$$(\partial \Omega )_x(r) = \{y \in \partial \Omega: (1-t)x + ty
\subset \Omega\,\,\, \hbox{for every $t \in (0,r)$}\}.$$ If $y \in
(\partial \Omega )_x(r)$, then there exists $\delta >0$ such that
$B_\delta (y) \cap \partial \Omega \subset (\partial \Omega )_x(r)$.
Thus, for each $r \in \mathbb Q \cap (0,1)$, the set $(\partial
\Omega )_x(r)$ is  open in $\partial \Omega$, in the topology
induced by $\rn$. The conclusion then follows from the fact that
$(\partial \Omega) _x = \cap _{r \in \mathbb Q \cap (0,1)}(\partial
\Omega )_x(r)$. \qed

\medskip
\par
Next, we define the sets
\begin{equation}\label{finite}
(\Omega \times \mathbb S ^{n-1})_0 = \{(x, \vartheta ) \in \Omega
\times \mathbb S ^{n-1}: x+t\vartheta \in \partial \Omega\,\,
\hbox{for some $t>0$}\},
\end{equation}
and
\begin{equation}\label{infinite}
(\Omega \times \mathbb S ^{n-1})_\infty = (\Omega \times \mathbb S
^{n-1}) \setminus (\Omega \times \mathbb S ^{n-1})_0.
\end{equation}
Clearly,
\begin{equation}\label{omegabound}
(\Omega \times \mathbb S ^{n-1})_0= \Omega \times \mathbb S ^{n-1}
\quad \hbox{if $\Omega $ is bounded.}
\end{equation}
 Let
 \begin{equation}\label{z}
 \zeta  : (\Omega \times \mathbb S ^{n-1})_0
\to \rn
\end{equation}
 be the function defined as
$$\zeta (x, \vartheta) = x + t \vartheta, \quad \hbox{where $t$ is
such that $x + t \vartheta \in (\partial \Omega) _x$}.$$
  In other
words, $\zeta (x, \vartheta)$ is the first point of intersection of
the half-line $\{ x + t \vartheta: t >0\}$ with $\partial \Omega$.
%
%
%
\par
 Given a function $g : \partial \Omega \to \R$, with compact
 support, we adopt the convention that  $g(\zeta (x, \vartheta))$ is defined for every
 $(x, \vartheta) \in \Omega \times \mathbb S ^{n-1}$, on extending it by $0$ on  $(\Omega \times \mathbb S ^{n-1})_\infty$; namely, we set
\begin{equation}\label{conv}
g(\zeta (x, \vartheta))= 0 \quad  \hbox{if $(x, \vartheta) \in
(\Omega \times \mathbb S ^{n-1})_\infty$.}
\end{equation}
%
Let us next introduce the functions
\begin{equation}\label{abdef}
\mathfrak a: \Omega \times \mathbb S ^{n-1} \to [-\infty , 0) \quad
\hbox{and}\quad \mathfrak b: \Omega \times \mathbb S ^{n-1} \to (0,
\infty]
\end{equation}
given by
\begin{equation}\label{rbx}
\mathfrak  b(x, \vartheta ) = \begin{cases} |\zeta(x, \vartheta )-x|
& \hbox{if $(x, \vartheta ) \in  (\Omega \times \mathbb S
^{n-1})_0$,}
\\ \infty & \hbox{otherwise},
\end{cases}
\end{equation}
and
\begin{equation}\label{rax}
\mathfrak a(x, \vartheta ) =   -\mathfrak b(x, -\vartheta )\quad
\hbox{if $(x, \vartheta ) \in \Omega \times \mathbb S ^{n-1}$.}
\end{equation}
%

\begin{proposition}\label{zx}
The function $\zeta$ is Borel measurable. Hence, the functions
$\mathfrak a$ and $\mathfrak b$ are Borel measurable as well.
%
%
\end{proposition}
{\bf Proof}. Assume first that $\Omega$ is bounded, so that $(\Omega
\times \mathbb S^{n-1})_0=\Omega \times \mathbb S^{n-1}$. Consider a
sequence of nested polyhedra $\{Q_k\}$ invading $\Omega$, and the
corresponding sequence of functions $\{\zeta_k\}$,  defined as
$\zeta$, with $\Omega$ replaced with $Q_k$. Such functions are Borel
measurable, by elementary considerations, and hence $\zeta$ is also
Borel measurable, since $\zeta_k$ converges to $\zeta$ pointwise.
\\ Next, assume that $\Omega$ is unbounded. For each $h \in \N$, consider the set
$\Omega _h = \Omega \cap B_h(0)$, where $B_h(0)$ is the ball,
centered at $0$, with radius $h$. Let $\zeta_h$ and $\mathfrak b_h$
be the functions, defined as $\zeta$ and $\mathfrak b$, with
$\Omega$ replaced with $\Omega _h$. Since $\Omega _h$ is bounded,
then we already know that $\mathfrak b_h$ is Borel measurable.
Moreover, $\mathfrak b_h$ converges to $\mathfrak b$ pointwise.
Hence, $\mathfrak b$ is Borel measurable as well, and in particular
the set $(\Omega \times \mathbb S^{n-1})_0$, which agrees with
$\{\mathfrak b<\infty\}$, is Borel measurable. Finally, the function
$\zeta_h$ is Borel  measurable, inasmuch as $\Omega _h$ is a bounded
set. Moreover,  $\zeta_h$ converges to $\zeta$ pointwise to $\zeta$
on the Borel set $(\Omega \times \mathbb S^{n-1})_0$. Thus, $\zeta$
is Borel measurable. \qed

Given $m \in \N$ and $ p \in [1, \infty]$, we denote by $V^{m, p}
(\Omega )$ the Sobolev type space defined as
\begin{equation}\label{sobolevV}
V^{m, p} (\Omega ) =   \big\{u: \hbox{$u$ is $m$-times weakly
differentiable in $\Omega$, and $|\nabla ^m u| \in
L^p(\Omega)$}\big\}.
\end{equation}
Let us notice that, in the definition of $V^{m, p} (\Omega )$, it is
only required that the derivatives of the highest order $m$ of $u$
belong to $L^p(\Omega )$. Replacing $L^p(\Omega)$ in
\eqref{sobolevV} with a more general Banach function space
$X(\Omega)$ leads to the notion of $m$-th order Sobolev type space
$V^mX(\Omega)$ built upon $X(\Omega)$.

For $k \in \N _0$, we denote as usual by $C^k(\overline \Omega)$ the
space of  real-valued functions whose $k$-th order derivatives in
$\Omega$  are continuous up to the boundary. We also set
\begin{equation}\label{cdot}
C^k_{\rm b}(\overline \Omega ) = \{u\in C^k(\overline \Omega) :
u\,\, \hbox{has bounded support} \}.
\end{equation}
Clearly, $$C^k_{\rm b}(\overline \Omega ) = C^k(\overline \Omega)
\quad \hbox{if $\Omega$ is bounded}.$$
\par
Let $\alpha = (\alpha _1, \dots , \alpha _n)$ be a multi-index with
$\alpha _i \in \N _0$ for $i=1, \dots , n$. We adopt the notations
$|\alpha|= \alpha _1 + \cdots + \alpha _n$, $\alpha ! = \alpha _1!
\cdots \alpha _n!$, and $\vartheta^\alpha = \vartheta^{\alpha _1}_1
\cdots \vartheta^{\alpha _n}_n$ for $\vartheta \in \rn$. Moreover,
we set $D^\alpha u = \frac{\partial ^{|\alpha |}u}{\partial
x_1^{\alpha _1} \dots \partial x_n^{\alpha _n}}$ for $u : \Omega \to
\R$.
\\
We need to extend the notion of upper gradient $g$ for the
restriction of $u$ to $\partial \Omega$ appearing in \eqref{hajlasz}
to the case of higher-order derivatives. To this purpose,  let us
denote by  $g^{k,j}$, where $k \in \N_0$ and $j=0,1$, $(k,j)\neq
(0,0)$,  any Borel function on $\partial \Omega$ fulfilling the
following property:
\par\noindent (i) If $k \in \N $, $j=0$, and $u
\in C^{k-1}_{\rm b}(\overline \Omega )$,
%
%
 %
%
%
\begin{equation}\label{M}
 \bigg|\sum _{|\alpha|\leq k-1} \frac{(2k-2-|\alpha|)! }{
(k-1-|\alpha|)!\alpha !} \frac{(y - x)^\alpha}{|y -
x|^{2k-1}}\Big[(-1)^{|\alpha |}D^\alpha u(y) - D ^\alpha
u(x)\Big]\bigg|\leq g^{k,0}  (x ) + g^{k,0}  (y )
\end{equation}
for $\hh$-a.e. $x, y \in \partial \Omega$.
\par\noindent (ii)
If $k \in \N $, $j=1$, and $u \in C^{k}_{\rm b}(\overline \Omega )$,
%
 %
%
%
\begin{equation}\label{M1}
 \sum _{i=1}^n\bigg|\sum _{|\alpha|\leq k-1} \frac{(2k-2-|\alpha|)! }{
(k-1-|\alpha|)!\alpha !} \frac{(y - x)^\alpha}{|y -
x|^{2k-1}}\Big[(-1)^{|\alpha |}D^\alpha \tfrac{\partial u}{\partial
x_i}(y) - D ^\alpha\tfrac{\partial u}{\partial
x_i}(x)\Big]\bigg|\leq g^{k,1}  (x ) + g^{k,1} (y)
\end{equation}
for $\hh$-a.e. $x, y \in \partial \Omega$.
\par\noindent (iii) If $k=0$, $j=1$, and  $u \in C^0_{\rm b}(\overline \Omega)$,
\begin{equation}\label{M0}
 |u(x)| \leq g^{0,1}(x)  
\end{equation}
for $\hh$-a.e. $x \in \partial \Omega$. Note that inequality
\eqref{M}, with $k=1$, agrees with \eqref{hajlasz}, and hence
$g^{1,0}$ has the same role as $g$ in \eqref{hajlasz}. Let us also
point out that,  as  \eqref{hajlasz} extends a classical property of
the gradient of weakly differentiable functions in $\rn$, likewise
its higher-order versions \eqref{M} and \eqref{M1} extend a parallel
property of functions in $\rn$ endowed with higher-order weak
derivatives. This is shown in Proposition \ref{maximaln} of the
Appendix.
\par\noindent
In analogy with \eqref{hajlasznorm}, we introduce the  seminorm
given, for $r \in [1, \infty]$, by
\begin{equation}\label{besov}
\|u\|_{\mathcal V^{k,j}L^r(\partial \Omega )}= \inf _{g^{k,j}
}\|g^{k,j} \|_{L^r (\partial \Omega)}
\end{equation}
where $k$, $j$ and $u$ are as above, and the infimum is extended
over all functions $g^{k,j}$ fulfilling the appropriate  definition
among \eqref{M}, \eqref{M1} and \eqref{M0}. More generally, given a
Banach function space $Z(\partial \Omega)$ on $\partial \Omega$ with
respect to the Hausdorff measure $\hh$, we define
\begin{equation}\label{besovZ}
\|u\|_{\mathcal V^{k,j}Z(\partial \Omega )}= \inf _{g^{k,j}
}\|g^{k,j} \|_{Z (\partial \Omega)}.
\end{equation}
Observe that, in particular,
 $$\|u\|_{\mathcal V^{0,1}Z(\partial
\Omega )}= \|u \|_{Z (\partial \Omega)}.$$

\section{Pointwise estimates}\label{proofs}

In the present section we establish our first main result: a
pointwise estimate for Sobolev functions, and their derivatives, in
arbitrary open sets. In what follows we define, for $k \in \N$,
\begin{equation}\label{natural} \natural (k) =
\begin{cases}   0 & \hbox{if $k$ is
odd,}
\\ 1 & \hbox{if $k$ is even.}
\end{cases}
\end{equation}

\begin{theorem}\label{intermest}
{\bf [Pointwise estimate]} Let $\Omega$ be any  open set in $\rn$,
$n \geq 2$. Assume that $m \in \N$ and $h  \in \N _0$ are such that
$0<m-h < n$. Then there exists a constant $C=C(n, m)$ such that
\begin{align}\label{hfund1ell}
|\nabla ^h u(x)| & \leq C \bigg(\int _\Omega \frac{|\nabla
^{m}u(y)|}{|x-y|^{n-m + h}}\, dy + \sum _{k=1}^{m -h- 1}\int
_{\Omega}  \int _{\mathbb S^{n-1}}\frac{g^{[\frac{k+h+1}2], \natural
( k+h )}(\zeta (y,\vartheta ))}{|x-y|^{n-k}} \, d\hh (\vartheta ) \,
dy
\\ \nonumber & \qquad \qquad +   \int _{\mathbb S^{n-1}}g^{[\frac{h+1}2], \natural ( h )} (\zeta (x, \vartheta ))\,
d\hh (\vartheta )\bigg) \qquad \hbox{for a.e. $x \in \Omega$,}
\end{align}
for every  $u \in V^{m , 1}(\Omega ) \cap C^{[\frac {m -1}2]}_{\rm
b}(\overline \Omega )$. Here, $g^{[\frac{k+h+1}2], \natural ( k+h
)}$ is any function as in \eqref{M}--\eqref{M0}, and convention
\eqref{conv} is adopted.
\end{theorem}

\begin{remark}\label{log} {\rm In the
case when $m-h =n$, and $\Omega$ is bounded, an estimate analogous
to \eqref{hfund1ell} can be proved, with the kernel
$\frac{1}{|x-y|^{n-m+h}}$ in the first integral on the right-hand
side replaced with $\log \frac C{|x-y|}$. The constant $C$ depends
on $n$ and the diameter of $\Omega$. If $m-h
>n$, and $\Omega$ is bounded, then the kernel is bounded by a constant depending on $n$, $m$ and the diameter of $\Omega$. }
\end{remark}

\begin{remark}\label{0}
{\rm Under the assumption that
\begin{equation}\label{zero}
 u=\nabla u= \dots \nabla ^{[\frac {m-1}2]}u=0 \quad \hbox{on $\partial
\Omega$, }
\end{equation}
one can choose $g^{[\frac{h+1}2], \natural ( h )}=0$ for $k=0, \dots
, m-h-1$ in \eqref{hfund1ell}. Hence,
\begin{align}\label{hfund1ell0}
|\nabla ^hu(x)| & \leq C \int _\Omega \frac{|\nabla
^{m}u(y)|}{|x-y|^{n-m+h}}\, dy   \qquad \hbox{for a.e. $x \in
\Omega$.}
\end{align}
A special case of \eqref{hfund1ell0}, corresponding to $h=m-1$, is
the object of \cite[Theorem 1.6.2]{Mabook}.}
\end{remark}

\begin{remark}\label{m-1mezzi}
{\rm As already mentioned in Section \ref{sec1}, the order
$\big[\tfrac {m-1}2\big]$ of the derivatives prescribed on $\partial
\Omega$, which appears on the right-hand side of \eqref{hfund1ell},
is minimal for Sobolev type inequalities to hold in arbitrary
domains. This issue is discussed in Example \ref{ex4}, Section
\ref{sharp} below.}
\end{remark}

\medskip
\par

A key step in the proof of Theorem \ref{intermest} is Lemma
\ref{fundest} below, which deals with the case when $h=m-1$ in
Theorem \ref{intermest}.

provides us with
 estimates for the $h$-th order
derivatives of  a function in terms of its $(h+1)$-th order
derivatives.

%

\begin{lemma}\label{fundest}
Let $\Omega$ be any  open set in $\rn$, $n \geq 2$.
\par\noindent
(i) If $u \in V^{2\ell-1,1}(\Omega ) \cap C^{\ell-1}_{\rm
b}(\overline \Omega )$ for some  $\ell \in \N$,  then
\begin{align}\label{fund1disp}
|\nabla ^{2\ell-2} u(x)| \leq C \bigg(\int _\Omega \frac{|\nabla
^{2\ell-1} u (y)|}{|x-y|^{n-1}}\, dy +  \int _{\mathbb S^{n-1}}
g^{\ell-1, 1}(\zeta  (x, \vartheta ))\, d\hh (\vartheta )\bigg) \,
\,\,\hbox{for a.e. $x \in \Omega$,}
\end{align}
 for some constant $C=C(n,\ell)$.
\par\noindent
(ii) If $u \in V^{2\ell,1}(\Omega ) \cap C^{\ell-1}_{\rm
b}(\overline \Omega )$ for some  $\ell \in \N$, then
\begin{align}\label{fund1}
|\nabla ^{2\ell-1} u(x)| \leq C\bigg( \int _\Omega \frac{|\nabla
^{2\ell} u (y)|}{|x-y|^{n-1}}\, dy +  \int _{\mathbb S^{n-1}}
g^{\ell, 0}(\zeta (x, \vartheta ))\, d\hh (\vartheta )\bigg)
\end{align}
 for some constant $C=C(n,\ell)$.
Here, $g^{\ell-1, 1}$ and  $g^{\ell, 0}$ are functions as in
 \eqref{M} -- \eqref{M0},and convention
\eqref{conv} is adopted.
\end{lemma}

Our proof of Lemma \ref{fundest} in turn requires the following
representation formula for the $(2\ell -1)$-th order derivative of a
one-dimensional function in an interval, in terms of its $2\ell$-th
derivative in the relevant interval, and of its  derivatives up to
the order $\ell -1$ evaluated at the endpoints.

\begin{lemma}\label{1dim}
Let $-\infty < a < b < \infty$. Assume that   $\psi \in W^{2\ell,
1}(a,b)$ for some  $\ell \in \N$. Then
\begin{align}\label{1dim1}
\psi ^{(2\ell-1)}(t)& = \int _a^t Q _{2\ell-1}\Big(\frac{2\tau
-a-b}{b-a}\Big) \psi ^{(2\ell)}(\tau )\, d\tau - \int _t^b Q
_{2\ell-1}\Big(\frac{a+b-2\tau }{b-a}\Big) \psi ^{(2\ell)}(\tau )\,
d\tau
\\ \nonumber & \quad + (2\ell-1)! (-1)^\ell\sum _{k=0}^{\ell-1} \frac{(2\ell-k-2)!}{k!
(\ell-k-1)!}\frac 1{(b-a)^{2\ell-k-1}}\big[(-1)^{k+1}\psi ^{(k)}(b)+
\psi ^{(k)}(a)\big]
\end{align}
for  $t \in (a,b)$. Here, $Q _{2\ell-1}$ is the polynomial of degree
$2\ell-1$, obeying
 \begin{equation}\label{Pi1}Q _{2\ell-1}(t) +
Q _{2\ell-1}(-t) = 1 \quad \hbox{for $t \in \R$,}
\end{equation}
 and
\begin{equation}\label{Pi2}
Q _{2\ell-1} (-1)= Q _{2\ell-1} ^{(1)}(-1) = \cdots = Q _{2\ell-1}
^{(\ell-1)}(-1) =0\,.
\end{equation}
\end{lemma}
\par\noindent
{\bf Proof}. Let us represent $\psi$ as
 \begin{equation}\label{1dim2}
 \psi (t) =
\varpi (t) + \varsigma (t) \quad
\hbox{for $t \in (a,
b)$,}
\end{equation}
 where $\varpi$ and $\varsigma$ are the solutions to the problems
 \begin{equation}\label{1dim2'}\begin{cases} \varpi ^{(2\ell)}(t) =  \psi ^{(2\ell)}(t) \quad \hbox{ in $(a, b)$,} \\
\varpi ^{(k)}(a) = \varpi ^{(k)}(b)=0 \quad \hbox{for $k=0, 1, \dots
, \ell-1$,}
\end{cases}
\end{equation}
and
 \begin{equation}\label{1dim2''}
 \begin{cases} \varsigma ^{(2\ell)}(t) =  0 \quad \hbox{in $(a, b)$,} \\
\varsigma ^{(k)}(a) =  \psi ^{(k)}(a)\,, \quad \varsigma ^{(k)}(b) =
\psi ^{(k)}(b) \quad \hbox{for $k=0, 1, \dots , \ell-1$,}
\end{cases}
\end{equation}
respectively. Let us first focus on problem \eqref{1dim2'}. We claim
that
\begin{equation}\label{maz1}
\varpi ^{(2\ell-1)}(t) = \int _a^t Q _{2\ell-1}\Big(\frac{2\tau
-a-b}{b-a}\Big) \psi ^{(2\ell)}(\tau )\, d\tau - \int _t^b Q
_{2\ell-1}\Big(\frac{a+b-2\tau }{b-a}\Big) \psi ^{(2\ell)}(\tau )\,
d\tau \quad \hbox{for $t \in (a, b)$,}
\end{equation}
where $Q _{2\ell-1}$ is  as in the statement. In order to verify
\eqref{maz1}, let us consider the auxiliary problem
\begin{equation}\label{maz2}
 \begin{cases} \omega  ^{2\ell}(s) = \phi (s) \quad \hbox{in $(-1, 1)$,} \\
 \omega ^{(k)} (\pm 1) =0, \quad \hbox{$k=0, 1, \dots , \ell -1$,}
 \end{cases}
\end{equation}
where $\phi \in L^1(a,b)$ is any given function. Let $\kappa :
[-1,1]^2 \to \R$ be the Green function associated with  problem
\eqref{maz2}, so that
\begin{equation}\label{maz3}
\omega (s) = \int _{-1}^1 \kappa (s,r)\phi(r)\, dr \quad \hbox{for
$s \in [-1,1]$.}
\end{equation}
The function $\kappa$ takes an explicit form (\cite{boggio}; see
also \cite[Section 2.6]{GGS}), given by
\begin{equation}\label{kappa}
\kappa (s,r) =  C |s-r|^{2\ell-1} \int _1^{\frac{1-sr}{|s-r|}}
(t^2-1)^{\ell-1}\, dt  \qquad \hbox{for $s \neq r$,}
\end{equation}
where  $C=C(\ell)$ is a suitable constant.
One can easily see from formula \eqref{kappa} that  $\kappa (s,r)$
is a polynomial of degree $2\ell-1$ in $s$ for fixed $r$, and a
polynomial of degree $2\ell-1$ in $r$ for fixed $s$, both in
$\{(s,r) \in [-1,1]^2: s>r\}$, and in $\{(s,r) \in [-1,1]^2: s<r\}$.
Moreover, $\kappa (s,r) = \kappa (-s,-r)$. In particular,  if $s>r$,
one has that
\begin{align}\label{kappa1}
 \kappa (s,r) = C \bigg[\sum _{j=0}^{\ell-1}
 \binom{\ell-1}{j}
 \frac{(-1)^{\ell-1-j}}{2j+1} (1-sr)^{2j+1} (s-r)^{2\ell-2k-2} -
 (s-r)^{2\ell-1} \sum _{j=0}^{\ell-1}
  \binom{\ell-1}{j}
 \frac{(-1)^{\ell-1-j}}{2j+1}\bigg].
 \end{align}
Thus, if $s>r$,
\begin{equation}\label{kappa2}
\frac {\partial ^{2\ell-1} \kappa }{\partial s^{2\ell-1}}(s,r) =
C(2\ell-1)! \bigg[\sum _{j=0}^{\ell-1}  \binom{\ell-1}{j}
 \frac{(-1)^{\ell+j}}{2j+1} r^{2j+1} -
 \sum _{j=0}^{\ell-1}  \binom{\ell-1}{j}
 \frac{(-1)^{\ell-1-j}}{2j+1}\bigg],
 \end{equation}
 a polynomial of degree $2\ell-1$ in $r$, depending only on odd powers of
 $r$.  Let us denote this polynomial by $Q _{2\ell-1}(r)$. It follows
 from \eqref{kappa2} that
 $Q _{2\ell-1}(-1) =0$. Moreover,
\begin{equation}\label{kappa3}
\frac {\partial Q _{2\ell-1}}{\partial r}(r) =
 \frac
{\partial}{\partial r} \bigg(\frac {\partial ^{2\ell-1} \kappa
}{\partial s^{2\ell-1}}(s,r)\bigg) = C (2\ell-1)! \sum
_{j=0}^{\ell-1} \binom{\ell-1}{j}
 (-1)^{\ell+j} r^{2j}  = - C(2\ell-1)!(r^2 -1)^{\ell-1}.
 \end{equation}
Thus, $Q _{2\ell-1}$ vanishes, together with all its derivatives up
to the order $\ell-1$, at $-1$, namely  $Q _{2\ell-1}$ fulfills
\eqref{Pi2}. Equation \eqref{kappa2} also tells us that
$Q_{2\ell-1}(s) - Q _{2\ell-1}(0)$ is an odd function, and hence
\begin{equation}\label{kappa3bis}
Q _{2\ell-1}(s) + Q_{2\ell-1}(-s)=2Q_{2\ell-1}(0)\quad \hbox{for $s
\in \R$.}
\end{equation}
%
%
%
%
\\ Since $\kappa$ is an even function,
$$\tfrac {\partial ^{2\ell-1} \kappa
}{\partial s^{2\ell-1}}(s,r) = \tfrac {\partial ^{2\ell-1} \kappa
}{\partial s^{2\ell-1}}(-s,-r)=  - Q _{2\ell-1}(-r) \quad \hbox{if
$-1\leq s<r \leq 1$.}$$ Thus,
 $(2\ell-1)$-times differentiation of equation \eqref{maz3} yields
\begin{equation}\label{maz4}
\omega ^{(2\ell-1)}(s) = \int _{-1}^s Q _{2\ell-1}(r) \phi(r) \, dr
- \int _s^1 Q _{2\ell-1}(-r) \phi(r) \, dr \quad \hbox{for $s \in
[-1,1]$.}
\end{equation}
Since $\omega ^{(2\ell)}=\phi$, an integration by parts in
\eqref{maz4}, equation \eqref{kappa3bis}, and the the fact that $Q
_{2\ell-1}(-1) =0$, tell us that
\begin{align}\label{maz4bis}
\omega ^{(2\ell-1)}(s) =  2Q _{2\ell-1}(0) \omega ^{(2\ell-1)}(s) -
\int _{-1}^s Q _{2\ell-1}'(r) \omega ^{(2\ell-1)}(r) \, dr - \int
_s^1 Q _{2\ell-1}'(-r) \omega ^{(2\ell-1)}(r) \, dr \quad \hbox{for
$s \in [-1,1]$.}
\end{align}
Owing to the arbitrariness of $\omega$, equation \eqref{maz4bis}
ensures that $2Q _{2\ell-1}(0)=1$. Equation \eqref{Pi1} thus follows
from \eqref{kappa3bis}.
\\
The function $\varpi$ defined as
$$\varpi (t) = \omega \Big(\tfrac{2t-a-b}{b-a}\Big)\quad \hbox{for
$t \in [a,b]$,}$$  is thus the solution to problem \eqref{1dim2'},
and the representation formula \eqref{maz1} follows via a change of
variables in \eqref{maz4}.
\par Consider next problem
\eqref{1dim2''}. The function $\varsigma$ is a polynomial of degree
$2\ell -1$, and $\varsigma ^{(2\ell-1)}$ is a constant which, owing
to the two-point Taylor interpolation formula (see e.g.
\cite[Chapter 2, Section 2.5, Ex. 3]{davis}), is given by
\begin{align}\label{1dim3}
\varsigma ^{(2\ell-1)} & = (2\ell-1)!
\bigg[\frac{d^{\ell-1}}{dt^{\ell-1}}\bigg(\frac{\varsigma
(t)}{(t-b)^\ell}\bigg)_{|t=a} +
\frac{d^{\ell-1}}{dt^{\ell-1}}\bigg(\frac{\varsigma
(t)}{(t-a)^\ell}\bigg)_{|t=b}\bigg] \\ \nonumber & = (2\ell-1)!
\bigg[\frac{d^{\ell-1}}{dt^{\ell-1}}\bigg(\frac{\psi
(t)}{(t-b)^\ell}\bigg)_{|t=a} +
\frac{d^{\ell-1}}{dt^{\ell-1}}\bigg(\frac{\psi
(t)}{(t-a)^\ell}\bigg)_{|t=b}\bigg].
\end{align}
Leibnitz' differentiation rule for products yields
\begin{multline}\label{1dim4}
 \bigg[\frac{d^{\ell-1}}{dt^{\ell-1}}\bigg(\frac{\psi
(t)}{(t-b)^\ell}\bigg)_{|t=a} +
\frac{d^{\ell-1}}{dt^{\ell-1}}\bigg(\frac{\psi
(t)}{(t-a)^\ell}\bigg)_{|t=b}\bigg] \\ = (-1)^\ell\sum
_{k=0}^{\ell-1} \frac{(2\ell-k-2)!}{k! (\ell-k-1)!}\frac
1{(b-a)^{2\ell-k-1}}\big[(-1)^{k+1}\psi ^{(k)}(b)+ \psi
^{(k)}(a)\big].
\end{multline}
 Equation \eqref{1dim1} follows from    \eqref{maz1}, \eqref{1dim3} and \eqref{1dim4}. \qed

\medskip
\par\noindent
{\bf Proof of Lemma \ref{fundest}}. Given $x \in \Omega$ and
$\vartheta \in \mathbb S^{n-1}$, let $\mathfrak a(x, \vartheta )$
and $\mathfrak b(x, \vartheta )$ be defined as in \eqref{rax} and
\eqref{rbx}, respectively.
\par\noindent
We begin with the proof of \eqref{fund1disp} for $\ell =1$. If $u
\in V^{1,1}(\Omega ) \cap C^{0}_{\rm b}(\overline \Omega )$, then,
by a standard property of Sobolev functions, for a.e. $x \in \Omega$
  the function
$$[0, \mathfrak b(x,\vartheta )] \ni t \mapsto u(x+t\vartheta
)$$
  belongs to $V^{1,1}(0 , \mathfrak b(x, \vartheta
  ))$  for $\hh$-a.e. $\vartheta \in \mathbb S^{n-1}$,
and
$$\frac{d}{dt} u(x+t\vartheta) = \nabla u(x + t\vartheta ) \cdot \vartheta \quad \hbox{for a.e.
$t \in [0, \mathfrak b(x,\vartheta )]$.}$$
 Hence, for any such $x$ and
$\vartheta$,
\begin{equation}\label{repr0}
 u(\zeta(x, \vartheta )) - u(x) = \int _0^{\mathfrak b(x,\vartheta )} \nabla u(x + t\vartheta ) \cdot \vartheta \, dt,
\end{equation}
where   convention \eqref{conv} is adopted.
 Integrating both sides of equation
\eqref{repr0} over $\mathbb S^{n-1}$ yields
\begin{equation}\label{repr2}
n\omega _{n} u(x) = \int _{\mathbb S^{n-1}}  u(\zeta(x, \vartheta
))\, d\hh (\vartheta ) - \int _{\mathbb S^{n-1}} \int _0^{\mathfrak
b(x,\vartheta )} \nabla u(x + t\vartheta ) \cdot \vartheta \, dt\,
d\hh (\vartheta ),
\end{equation}
where $\omega _n = \Pi ^{\frac n2}/\Gamma (1+\tfrac n2)$, the
Lebesgue measure of the unit ball in $\rn$. One has that
\begin{align}\label{repr4}
\int _{\mathbb S^{n-1}} \int _0^{\mathfrak b(x,\vartheta )} \nabla
u(x + t\vartheta ) \cdot \vartheta \, dt\, d\hh (\vartheta )
 &
= \int _{\mathbb S^{n-1}} \int _0^{\mathfrak b(x,\vartheta ) } \frac
1{t ^{n-1}}\nabla u(x +\varrho \vartheta ) \cdot \vartheta \,
t^{n-1} d\varrho \, d\hh (\vartheta )
\\ \nonumber & =
\int _{\Omega _x} \frac{\nabla u (y)\cdot (y-x)}{|x-y|^{n}}\, dy.
\end{align}
Inequality \eqref{fund1disp}, with $\ell =1$,  follows from
\eqref{repr2} and \eqref{repr4}.
\par\noindent
Let us next  prove \eqref{fund1}.
%
 If $u \in V^{2\ell,1}(\Omega ) \cap C^{\ell-1}_{\rm b}(\overline
\Omega )$, then for a.e. $x \in \Omega$, the function
$$[\mathfrak a(x, \vartheta ) , \mathfrak b(x,\vartheta )] \ni t \mapsto u(x+t\vartheta
)$$
  belongs to $V^{2\ell,1}(\mathfrak a(x, \vartheta ) , \mathfrak b(x, \vartheta
  ))$  for $\hh$-a.e. $\vartheta \in \mathbb S^{n-1}$. Consider any such $x$ and $\vartheta$. If  $(x,
  \vartheta)   \in (\Omega \times \mathbb
S^{n-1})_0$, then,  by
  Lemma \ref{1dim},
 \begin{align}\label{fund2}
 &  \frac{d^{2\ell-1}}{dt^{2\ell-1}}u(x+t\vartheta )  = \int _{\mathfrak a(x,\vartheta )}^t Q _{2\ell-1}\Big(\frac{2\tau
-\mathfrak a(x,\vartheta )-\mathfrak b(x,\vartheta )}{\mathfrak
b(x,\vartheta )-\mathfrak a(x,\vartheta
)}\Big)  \frac{d^{2\ell}}{d\tau^{2\ell}}u(x+\tau\vartheta )\, d\tau \\
\nonumber & - \int _t^{\mathfrak b(x,\vartheta )} Q
_{2\ell-1}\Big(\frac{\mathfrak a(x,\vartheta )+\mathfrak
b(x,\vartheta )-2\tau }{\mathfrak b(x,\vartheta )-\mathfrak
a(x,\vartheta )}\Big)
 \frac{d^{2\ell}}{d\tau ^{2\ell}}u(x+\tau\vartheta )\, d\tau
\\ \nonumber & \quad + (2\ell-1)! (-1)^\ell\sum _{k=0}^{\ell-1} \frac{(2\ell-k-2)!}{k!
(\ell-k-1)!}\frac{\big[(-1)^{k+1}\big(\frac{d^{k}u(x+t\vartheta
)}{dt^k}\big)_{|t=\mathfrak b(x,\vartheta )}+
\big(\frac{d^{k}u(x+t\vartheta )}{dt^k}\big)_{|t=\mathfrak
a(x,\vartheta )}\big]}{(\mathfrak b(x,\vartheta )-\mathfrak
a(x,\vartheta ))^{2\ell-k-1}}
\end{align}
for  $t \in (\mathfrak a(x,\vartheta ),\mathfrak b(x,\vartheta ))$.
If, instead, $(x, \vartheta)   \in (\Omega \times \mathbb
S^{n-1})_\infty$, then,
\begin{align}\label{fund2inf}
\frac{d^{2\ell-1}}{dt^{2\ell-1}}u(x+t\vartheta ) = \begin{cases} -
\int _t^\infty \frac{d^{2\ell}}{d\tau^{2\ell}}u(x+\tau\vartheta )\,
d\tau, &
\hbox{if $\mathfrak b(x,\vartheta ) = \infty$,} \\
 \int _{-\infty} ^t \frac{d^{2\ell}}{d\tau^{2\ell}}u(x+\tau\vartheta )\,
d\tau, & \hbox{if $\mathfrak a(x,\vartheta ) = -\infty$,}
\end{cases}
\end{align}
for  $t \in (\mathfrak a(x,\vartheta ),\mathfrak b(x,\vartheta ))$
(if both $\mathfrak b(x,\vartheta ) = \infty$ and $\mathfrak
a(x,\vartheta ) = -\infty$, then either expression on the right-hand
side of \eqref{fund2inf} can be exploited).
\\
We have that
\begin{equation}\label{fund3}
\frac{d^{k}}{dt^k}u(x+t\vartheta ) = \sum _{|\alpha |=k}
\frac{k!}{\alpha !}\vartheta ^\alpha D^\alpha u(x+ t\vartheta )
\quad \hbox{for a.e. $t \in (\mathfrak a(x, \vartheta ) , \mathfrak
b(x,\vartheta ))$,}
\end{equation}
for $k =1, \dots , 2\ell$. From \eqref{fund2}--\eqref{fund3} we
infer that
\begin{align}\label{fund4}
 & \sum _{|\alpha |=2\ell-1}
\frac{(2\ell-1)!}{\alpha !}\vartheta ^\alpha D^\alpha u(x)  =
-\chi_{(\Omega \times \mathbb S^{n-1})_\infty}(x, \vartheta ) \sum
_{|\alpha |=2\ell-1} \frac{(2\ell-1)!}{\alpha !}\vartheta ^\alpha
\int _0^{\pm\infty} \sum _{|\gamma |=1}\vartheta ^\gamma D^{\alpha +
\gamma} u(x+ \tau\vartheta )\, d\tau
\\
\nonumber  & \qquad \qquad \quad + \chi_{(\Omega \times \mathbb
S^{n-1})_0}(x, \vartheta )\bigg[\int _{\mathfrak a(x,\vartheta )}^0
Q _{2\ell-1}\Big(\frac{2\tau -\mathfrak a(x,\vartheta )-\mathfrak
b(x,\vartheta )}{\mathfrak b(x,\vartheta )-\mathfrak a(x,\vartheta
)}\Big) \frac{d^{2\ell}}{d\tau^{2\ell}}u(x+\tau\vartheta )
d\tau \\
\nonumber & \qquad \qquad \quad - \int _0^{\mathfrak b(x,\vartheta
)} Q _{2\ell-1}\Big(\frac{\mathfrak a(x,\vartheta )+\mathfrak
b(x,\vartheta )-2\tau }{\mathfrak b(x,\vartheta )-\mathfrak
a(x,\vartheta )}\Big)
\frac{d^{2\ell}}{d\tau^{2\ell}}u(x+\tau\vartheta )
\, d\tau
\\ \nonumber & \quad + (2\ell-1)! (-1)^\ell\sum _{k=0}^{\ell-1} \frac{(2\ell-k-2)!}{k!
(\ell-k-1)!} \sum _{|\alpha |=k} \frac{k!}{\alpha !}\vartheta
^\alpha \frac{\big[(-1)^{k+1}D^\alpha u(x+ \mathfrak b(x,\vartheta )
\vartheta ) + D^\alpha u(x+ \mathfrak a(x,\vartheta ) \vartheta )
 \big]}{(\mathfrak b(x,\vartheta )-\mathfrak a(x,\vartheta ))^{2\ell-k-1}}\bigg]\,,
\end{align}
where the signs $+$ or $-$ in the first integral on the right-hand
side depend on whether $\mathfrak b(x, \vartheta)=\infty$ or
$\mathfrak a(x, \vartheta)=-\infty$, respectively.
\par
 Denote by
$\{P_\beta \}$ the system of all homogeneous polynomials of degree
$2\ell-1$ in the variables $\vartheta_1, \dots , \vartheta _n$ such
that
$$\int _{\mathbb S^{n-1}}P_\beta (\vartheta ) \vartheta ^\alpha d\hh
(\vartheta) = \delta _{\alpha \beta},$$ where $\delta _{\alpha
\beta}$ stands for the Kronecker delta. On multiplying equation
\eqref{fund4} by $P_\beta (\vartheta )$, dividing through by $(2
\ell -1)!$, and integrating over $\sn$ one obtains that
\begin{align}\label{fund5}
  & \frac{D^\beta u(x)}{\beta !}   =
- \int _{\sn} \chi_{(\Omega \times \mathbb S^{n-1})_\infty}(x,
\vartheta ) P_\beta (\vartheta )  \sum _{\substack{ |\alpha
|=2\ell-1 \\ |\gamma |=1}} \frac{\vartheta ^{\alpha +
\gamma}}{\alpha !} \int _0^{\pm\infty}  D^{\alpha + \gamma} u(x+
\tau\vartheta )\, d\tau \, d\hh (\vartheta )
\\
\nonumber & \qquad \qquad \quad + \int _{\sn} \chi_{(\Omega \times
\mathbb S^{n-1})_0}(x, \vartheta ) \frac{P_\beta (\vartheta
)}{(2\ell-1)!}\bigg[ \int _{\mathfrak a(x,\vartheta )}^0 Q
_{2\ell-1}\Big(\frac{2\tau -\mathfrak a(x,\vartheta )-\mathfrak
b(x,\vartheta )}{\mathfrak b(x,\vartheta )-\mathfrak a(x,\vartheta
)}\Big) \frac{d^{2\ell}}{d\tau^{2\ell}}u(x+\tau\vartheta )
d\tau \\
\nonumber & \qquad \qquad \quad -  \int _0^{\mathfrak b(x,\vartheta
)} Q _{2\ell-1}\Big(\frac{\mathfrak a(x,\vartheta )+\mathfrak
b(x,\vartheta )-2\tau }{\mathfrak b(x,\vartheta )-\mathfrak
a(x,\vartheta )}\Big)
\frac{d^{2\ell}}{d\tau^{2\ell}}u(x+\tau\vartheta )
\, d\tau\bigg]d\hh
 (\vartheta) \,
\\ \nonumber &  +  (-1)^\ell \int _{\sn}P_\beta (\vartheta )\sum _{|\alpha | \leq \ell-1} \frac{(2\ell-|\alpha|-2)!\,\vartheta ^\alpha}{
(\ell-|\alpha|-1)!\alpha !} \frac{\big[(-1)^{|\alpha|+1}D^\alpha
u(x+ \mathfrak b(x,\vartheta ) \vartheta ) + D^\alpha u(x+ \mathfrak
a(x,\vartheta ) \vartheta )\big]}{(\mathfrak b(x,\vartheta
)-\mathfrak a(x,\vartheta ))^{2\ell-|\alpha|-1}} d\hh
 (\vartheta).
\end{align}
There exist a constants $C=C(n,\ell)$ and $C'=C'(n,\ell)$ such that
\begin{align}\label{fund5inf}
\bigg|& \int _{\sn}\chi_{(\Omega \times \mathbb S^{n-1})_\infty}(x,
\vartheta ) P_\beta (\vartheta )  \sum _{\substack{ |\alpha
|=2\ell-1 \\ |\gamma |=1}} \frac{\vartheta ^{\alpha +
\gamma}}{\alpha !} \int _0^{\pm\infty }D^{\alpha + \gamma} u(x+
\tau\vartheta )\, d\tau \, d\hh (\vartheta )\bigg|
\\ \nonumber &
\leq C \int _{\sn} \int _0^\infty |\nabla ^{2\ell} u(x+
\tau\vartheta )|\, d\tau \, d\hh (\vartheta ) = C \int _{\sn} \int
_0^\infty \frac{|\nabla ^{2\ell} u(x+ \tau\vartheta )|}{\tau ^{n-1}}
\tau ^{n-1}\, d\tau \, d\hh (\vartheta )
\\ \nonumber & \leq C' \int _\Omega \frac{|\nabla ^{2\ell}
u(y)|}{|x-y|^{n-1}}\, dy.
\end{align}
Next, we claim that
there exists a constant $C=C(\ell, n)$ such that
\begin{align}\label{fund6}
\bigg|&\int _{\sn}  \chi_{(\Omega \times \mathbb S^{n-1})_0}(x,
\vartheta ) \bigg[P_\beta (\vartheta ) \int _{\mathfrak
a(x,\vartheta )}^0 Q _{2\ell-1}\Big(\frac{2\tau -\mathfrak
a(x,\vartheta )-\mathfrak b(x,\vartheta )}{\mathfrak b(x,\vartheta
)-\mathfrak a(x,\vartheta )}\Big)
\frac{d^{2\ell}}{d\tau^{2\ell}}u(x+\tau\vartheta )
d\tau \\
\nonumber & -P_\beta (\vartheta ) \int _0^{\mathfrak b(x,\vartheta
)} Q _{2\ell-1}\Big(\frac{\mathfrak a(x,\vartheta )+\mathfrak
b(x,\vartheta )-2\tau }{\mathfrak b(x,\vartheta )-\mathfrak
a(x,\vartheta )}\Big)
\frac{d^{2\ell}}{d\tau^{2\ell}}u(x+\tau\vartheta ) d\tau \bigg]\,
d\hh
 (\vartheta)\bigg|
 \\ \nonumber & \leq C  \int _\Omega \frac{|\nabla ^{2\ell} u
(y)|}{|x-y|^{n-1}}\, dy.
\end{align}
In order to prove \eqref{fund6}, observe that
\begin{align*}
&\int _{\sn}  \chi_{(\Omega \times \mathbb S^{n-1})_0}(x, \vartheta
) P_\beta (\vartheta ) \int _{\mathfrak a(x,\vartheta )}^0 Q
_{2\ell-1}\Big(\frac{2\tau -\mathfrak a(x,\vartheta )-\mathfrak
b(x,\vartheta )}{\mathfrak b(x,\vartheta )-\mathfrak a(x,\vartheta
)}\Big) \frac{d^{2\ell}}{dr^{2\ell}}u(x+\tau\vartheta  )
d\tau\, d\hh
 (\vartheta)
 \\ \nonumber & = \int _{\sn}  \chi_{(\Omega \times \mathbb S^{n-1})_0}(x, \vartheta
) P_\beta (\vartheta ) \int _0^{-\mathfrak a(x,\vartheta )} Q
_{2\ell-1}\Big(-\frac{2r +\mathfrak a(x,\vartheta )+\mathfrak
b(x,\vartheta )}{\mathfrak b(x,\vartheta )-\mathfrak a(x,\vartheta
)}\Big) \frac{d^{2\ell}}{d\tau^{2\ell}}u(x-r\vartheta )
dr\, d\hh
 (\vartheta)
 \\ \nonumber & = \int _{\sn}  \chi_{(\Omega \times \mathbb S^{n-1})_0}(x,
 -\theta
) P_\beta (-\theta ) \int _0^{-\mathfrak a(x,-\theta )} Q
_{2\ell-1}\Big(-\frac{2r +\mathfrak a(x,-\theta )+\mathfrak
b(x,-\theta )}{\mathfrak b(x, -\theta )-\mathfrak a(x, -\theta
)}\Big) \frac{d^{2\ell}}{dr^{2\ell}}u(x+r\theta )
dr\, d\hh
 (\theta)
 \\ \nonumber & = -\int _{\sn}  \chi_{(\Omega \times \mathbb S^{n-1})_0}(x,
 \theta
)P_\beta (\theta ) \int _0^{\mathfrak b(x,\theta )} Q
_{2\ell-1}\Big(\frac{\mathfrak a(x,\theta )+\mathfrak b(x,\theta
)-2r}{\mathfrak b(x, \theta )-\mathfrak a(x, \theta )}\Big)
\frac{d^{2\ell}}{dr^{2\ell}}u(x+r\theta )
dr\, d\hh
 (\theta),
 \end{align*}
 where we have made use of the fact that $P_\beta (-\theta ) = -
 P_\beta (\theta)$ if $|\beta |= 2\ell-1$, and of \eqref{rax}. Thus,
\begin{align*}
\bigg|\int _{\sn}&  \chi_{(\Omega \times \mathbb S^{n-1})_0}(x,
\vartheta ) \bigg[P_\beta (\vartheta ) \int _{\mathfrak
a(x,\vartheta )}^0 Q _{2\ell-1}\Big(\frac{2\tau -\mathfrak
a(x,\vartheta )-\mathfrak b(x,\vartheta )}{\mathfrak b(x,\vartheta
)-\mathfrak a(x,\vartheta )}\Big)
\frac{d^{2\ell}}{d\tau^{2\ell}}u(x+\tau\vartheta )
d\tau \\
\nonumber & -P_\beta (\vartheta ) \int _0^{\mathfrak b(x,\vartheta
)} Q _{2\ell-1}\Big(\frac{\mathfrak a(x,\vartheta )+\mathfrak
b(x,\vartheta )-2\tau }{\mathfrak b(x,\vartheta )-\mathfrak
a(x,\vartheta )}\Big)
\frac{d^{2\ell}}{d\tau^{2\ell}}u(x+\tau\vartheta ) d\tau \bigg]\,
d\hh
 (\vartheta)\bigg|
\\
\nonumber & = 2 \bigg|\int _{\sn}P_\beta (\vartheta ) \int
_0^{\mathfrak b(x,\vartheta )} Q _{2\ell-1}\Big(\frac{\mathfrak
a(x,\vartheta )+\mathfrak b(x,\vartheta )-2r }{\mathfrak
b(x,\vartheta )-\mathfrak a(x,\vartheta )}\Big)
\frac{d^{2\ell}}{d\tau^{2\ell}}u(x+r\vartheta ) \,dr d\hh
 (\vartheta)\bigg|
\\
\nonumber & = 2 \bigg|\int _{\sn}P_\beta (\vartheta ) \int
_0^{\mathfrak b(x,\vartheta )} Q _{2\ell-1}\Big(\frac{\mathfrak
a(x,\vartheta )+\mathfrak b(x,\vartheta )-2r }{\mathfrak
b(x,\vartheta )-\mathfrak a(x,\vartheta )}\Big) \sum _{|\alpha
|=2\ell}\frac{(2\ell)!}{\alpha !} \vartheta ^\alpha D^\alpha
u(x+r\vartheta)
 \, drd\hh
 (\vartheta)\bigg|
\\
\nonumber & \leq C\int _{\sn}\int _0^{\mathfrak b(x,\vartheta )}  |
\nabla ^{2\ell} u(x+r\vartheta)|
 \, d\hh
 (\vartheta)
 \leq   C' \int_{\Omega _x}\frac{|\nabla ^{2\ell}
u(y)|}{|x-y|^{n-1}}\, dy
\end{align*}
for some constants $C=C(n,\ell)$ and $C'=C'(n, \ell)$. Hence,
inequality \eqref{fund6}  follows.
%
%
\\ Finally, by   definition \eqref{M}, there exists a
constant $C=C(n,\ell)$ such that
\begin{align}\label{fund7}
\bigg|& \int _{\sn} \chi_{(\omega \times \mathbb S^{n-1})_0}(x,
\vartheta )  (-1)^\ell
 P_\beta (\vartheta )
 \\ \nonumber & \times\sum _{|\alpha
| \leq \ell-1} \frac{(2\ell-|\alpha|-2)!}{ (\ell-|\alpha|-1)!\alpha
!}  \vartheta ^\alpha \frac{\big[(-1)^{|\alpha|+1}D^\alpha u(x+
\mathfrak b(x,\vartheta ) \vartheta ) + D^\alpha u(x+ \mathfrak
a(x,\vartheta ) \vartheta )\big]}{(\mathfrak b(x,\vartheta
)-\mathfrak a(x,\vartheta ))^{2\ell-|\alpha|-1}}\, d\hh
 (\vartheta)\bigg|
 \\ \nonumber & \leq C \int _{\sn} \chi_{(\omega \times \mathbb S^{n-1})_0}(x,
\vartheta )  \big[g^{\ell, 0}   (x+ \vartheta \mathfrak a(x,
\vartheta) )+ g^{\ell, 0}   (x+ \vartheta \mathfrak b(x, \vartheta)
)\big]\, d\hh (\vartheta )
\\ \nonumber & =
2C \int _{\sn} g ^{\ell, 0} (\zeta(x, \vartheta ) )\, d\hh
(\vartheta ).
\end{align}
%
%
%
%
%
%
Combining \eqref{fund5}--\eqref{fund7} yields
 \eqref{fund1}.
\par
Inequality \eqref{fund1disp}, with $\ell \geq 2$, follows on
applying \eqref{fund1} with $u$ replaced with its first-order
derivatives.
 \qed

\medskip
\par\noindent
{\bf Proof of Theorem \ref{intermest}}.  For simplicity of notation,
we consider the case when $h=0$, the proof in the general case being
analogous. Let $u \in V^{m , 1}(\Omega ) \cap C^{[\frac {m
-1}2]}_{\rm b}(\overline \Omega )$. By inequality \eqref{fund1disp}
with $\ell =1$,
\begin{align}\label{hfund1}
|u(x)| \leq C \bigg(\int _\Omega \frac{|\nabla  u
(y)|}{|x-y|^{n-1}}\, dy +   \int _{\mathbb S^{n-1}}g^{0, 1}(\zeta(x,
\vartheta ))\, d\hh (\vartheta )\bigg)\, \,\quad \hbox{for a.e. $x
\in \Omega$.}
\end{align}
From \eqref{hfund1} and an application of inequality \eqref{fund1}
with $\ell=1$ one obtains that
\begin{align}\label{hfund2}
|u(x)| & \leq C \bigg(\int _\Omega \int _\Omega \frac{|\nabla ^2 u
(z)|}{|y-z|^{n-1}}\,  \, \frac {dz\,dy}{|x-y|^{n-1}} \\ \nonumber &
\quad + \int _\Omega \int _{\mathbb S^{n-1}}g^{1,0}(\zeta ( y,
\vartheta ))\, \frac {d\hh (\vartheta )\, dy}{|x-y|^{n-1}} + \int
_{\mathbb S^{n-1}}g^{0, 1}(\zeta(x, \vartheta ))\, d\hh (\vartheta
)\bigg)
\\ \nonumber &
\leq C' \bigg(\int _\Omega  \frac{|\nabla ^2 u (z)|}{|x-z|^{n-2}}\,
dz  + \int _\Omega \int _{\mathbb S^{n-1}}g^{1,0}(\zeta(y, \vartheta
))\, \frac
{d\hh (\vartheta )\, dy}{|x-y|^{n-1}}  \\
\nonumber & \quad + \int _{\mathbb S^{n-1}}g^{0, 1}(\zeta(x,
\vartheta ))\, d\hh (\vartheta )\bigg) \,\qquad \hbox{for a.e. $x
\in \Omega$,}
\end{align}
for some constants $C=C(n)$ and $C'=C'(n)$. Note that in the last
inequality we have made use of a special case of the well known
identity
\begin{equation}\label{hf2}
\int _\rn \frac{1}{|x-y|^{n-\sigma }}\int _\rn
\frac{f(z)}{|y-z|^{n-\gamma }}\, dz\, dy = C \int _{\rn}
\frac{f(z)}{|x-z|^{n- \sigma-\gamma}}\, dz \quad \hbox{for a.e. $x
\in \rn$,}
\end{equation}
which holds for some constant $C=C(n,\sigma, \gamma)$ and for every
compactly supported integrable function $f$, provided that $\sigma
>0$, $\gamma >0$ and $\sigma + \gamma < n$.
\\ Inequality \eqref{hfund2} in turn yields, via an application of inequality \eqref{fund1disp} with
$m=2$,
\begin{align}\label{hfund3}
|u(x)|  \leq C \bigg(& \int _\Omega  \frac{|\nabla ^3 u
(z)|}{|x-z|^{n-3}}\, dz  + \int _\Omega \int _{\mathbb
S^{n-1}}g^{1,1}(\zeta(y, \vartheta ))\, \frac {d\hh (\vartheta
)\,dy}{|x-y|^{n-2}}
\\
\nonumber & + \int _\Omega \int _{\mathbb S^{n-1}}g^{1,0}(\zeta (y,
\vartheta ))\, \frac {d\hh (\vartheta )\,dy}{|x-y|^{n-1}}
 + \int _{\mathbb S^{n-1}}g^{0, 1}(\zeta (x, \vartheta ))\, d\hh (\vartheta )\bigg) \,\,\hbox{for a.e.
$x \in \Omega$,}
\end{align}
for some constant $C=C(n)$. A finite induction argument, relying
upon an alternate iterated use of inequalities \eqref{fund1} and
\eqref{fund1disp} as above, eventually leads to \eqref{hfund1ell}. \qed

\section{Estimates in rearrangement form}\label{sobolev}

The pointwise bounds established in the previous section enable us
to derive rearrangement estimates for functions, and their
derivatives, with respect to any Borel measure $\mu$ on $\Omega$
such that
\begin{equation}\label{measure}
\mu (B_r(x) \cap \Omega ) \leq C_\mu r^\alpha\quad \hbox{for $x \in
\Omega$ and  $r>0$,}
\end{equation}
for some $\alpha \in (n-1, n]$ and some constant $C_\mu>0$. Here,
$B_r(x)$ denotes the ball, centered at $x$, with radius $r$.
\par
Recall that, given a  measure space $\mathcal R $, endowed with a
positive measure $\nu$,  the decreasing rearrangement $\phi_\nu\sp*
: [0, \infty ) \to [0, \infty]$ of a $\nu$-measurable function $\phi
: \mathcal R \to \R$
%
%
%
 is defined as
$$
\phi_\nu\sp*(s)=\sup\{{t\geq 0}:\,\nu (|\phi|>t\})>s\}\quad
\textup{for}\ s\in [0,\infty).
$$
The operation of decreasing rearrangement is not linear. However,
one has that
\begin{equation}\label{f+g}
(\phi + \psi)_\nu\sp*(s) \leq \phi _\nu\sp*(s/2) + \psi
_\nu\sp*(s/2)\quad \hbox{for $s \geq 0$,}
\end{equation}
for every measurable functions $\phi$ and $\psi$ on $\mathcal R$.
\\ Any function $\phi$ shares its integrability properties with its decreasing rearrangement $\phi _\nu^*$, since
$$\nu (\{|\phi|>t\}) = \mathcal L^1(\{\phi _\nu^* >t\}) \quad
\hbox{for every $t\geq 0$.}$$
 As a consequence, any norm inequality, involving
 rearrangement-invariant norms, between the rearrangements of the derivatives of Sobolev
 functions and the rearrangements of its lower-order derivatives,
 immediately yields a corresponding  inequality for the original Sobolev functions. Thus, the
 rearrangement inequalities to be established hereafter reduce the
 problem of $n$-dimensional Sobolev type inequalities in arbitrary
 open sets to considerably simpler one-dimensional Hardy type
 inequalities -- see Theorem \ref{reduction}, Section \ref{ineq}
 below.

\begin{theorem}\label{rearrintermest} {\bf [Rearrangement estimates]}
Let $\Omega$ be any open bounded open set in $\rn$, $n \geq 2$. Let
$m \in \N$ and $h \in \N _0$ be such that  $0<m-h < n$. Assume that
$\mu$ is a Borel measure in $\Omega$ fulfilling \eqref{measure} for
some $\alpha \in (n-1, n]$ and for some $C_\mu
>0$. Then there exists  constants $c=c(n,m)$ and $C=C(n,m, \alpha ,
C_\mu)$ such that
\begin{align}\label{rearresthell}
|\nabla ^h u|_\mu^*(cs)  & \leq C\bigg[s^{-\frac
{n-m+h}{\alpha}}\int _0^{s^{\frac n\alpha}}|\nabla ^{m} u|_{\mathcal
L^n}^*(r) dr + \int _{s^{\frac n\alpha}}^\infty r^{-\frac{n-m+h}n}
|\nabla ^{m} u|_{\mathcal L^n}^*(r) dr
\\ \nonumber & \quad \quad +
\sum _{k=1}^{m-h-1}  \bigg(s^{-\frac {n-1-k}{\alpha}}\int
_0^{s^{\frac {n-1}\alpha}}  \big[g^{[\frac{k+h+1}2], \natural (
k+h)}  \big]_{\hh}^*(r)dr
\\ \nonumber & \quad \quad +  \int_{s^{\frac
{n-1}\alpha}}^\infty r^{-\frac {n-1-k}{n-1}}
\big[g^{[\frac{k+h+1}2], \natural (k+h )}  \big]_{\hh}^*(r)dr\bigg)
\\ \nonumber & \quad \quad + s^{-\frac {n-1}{\alpha}}\int _0^{s^{\frac {n-1}\alpha}}
\big[g^{[\frac{h+1}2], \natural ( h )} \big]_{\hh}^*(r)dr\bigg]\quad
\quad \hbox{for $s>0$,}
\end{align}
for every $u \in V^{m , 1}(\Omega ) \cap C^{[\frac {m -1}2]}_{\rm
b}(\overline \Omega )$. Here, $\natural (\cdot)$ is defined as in
\eqref{natural}, and $g^{[\frac{k+h+1}2], \natural ( k+h)}$ denotes
any Borel function on $\partial \Omega$ fulfilling the appropriate
condition from \eqref{M}--\eqref{M0}.
\end{theorem}

\begin{remark}\label{traces} {\rm In inequality \eqref{rearresthell}, and
in what follows, when considering rearrangements and norms with
respect to a measure $\mu$, Sobolev functions and their derivatives
have to be interpreted as their traces with respect to $\mu$. Such
traces are well defined, thanks to standard (local) Sobolev
inequalities with measures, owing to the assumption that $\alpha \in
(n-1, n]$ in \eqref{measure}. An analogous convention applies to the
integral operators to be considered below.}
\end{remark}

\par
In preparation for the proof of Theorem  \ref{rearrintermest}, we
introduce a few integral operators, and  pointwise estimate their
rearrangements.
\par\noindent
Let $\Omega$ be any open set in $\rn$. We define the operator $T$ as
\begin{equation}\label{T} Tg(x) = \int _{\mathbb S^{n-1}}|g(\zeta(x,
\vartheta ))|\,d\hh ( \vartheta ) \quad \hbox{ for $x \in \Omega$,}
\end{equation}
at any function Borel function $g:
\partial \Omega \to \R$. Here, and in what follows, we adopt convention
\eqref{conv}. Note that, owing to Fubini's theorem, $Tg$ is a
measurable function with respect to any Borel measure in $\Omega$.
\par\noindent
For $\gamma \in (0, n)$, we denote by $I_\gamma$ the classical Riesz
potential operator given by
\begin{equation}\label{I}
I_\gamma f (x) = \int _\Omega \frac{f(y)}{|y-x|^{n-\gamma}}\, dy
\quad \hbox{for   $x \in \Omega$,}
\end{equation}
at any $f \in L^1(\Omega)$, and we call $N_\gamma$ the operator
defined as
\begin{equation}\label{N}
N_\gamma g (x) =  \int _{\partial \Omega }
\frac{g(y)}{|x-y|^{n-\gamma}}\, d\hh (y) \, \,\,\hbox{for   $x \in
\Omega$,}
\end{equation}
at any function $g\in L^1(\partial \Omega)$.
\par\noindent
Finally,   we define the operator $Q_\gamma$ as the composition
\begin{equation}\label{comp}
Q_\gamma  = I_\gamma \circ T\,.
\end{equation}
Namely, \begin{equation}\label{Q} Q_\gamma g (x) = \int _{\Omega }
\int _{\mathbb S^{n-1}}|g(\zeta(y, \vartheta ))|\, \frac{d\hh
(\vartheta )\, dy}{|x-y|^{n-\gamma }}\qquad \hbox{for   $x \in
\Omega$,}
\end{equation}
for any Borel function $g:
\partial \Omega \to \R$.

Our analysis of these operators requires a few notations and
properties from interpolation theory.
\par
Assume that $\mathcal R $ is a  measure space, endowed with a
positive measure $\nu$.
%
%
%
Given a pair $X_1(\mathcal R)$ and $X_2(\mathcal R)$ of normed
function spaces, a function $\phi \in X_1(\mathcal R) + X_2(\mathcal
R)$ and $s \in \R$, we denote by $K(\phi, s ;X_1(\mathcal R ),
X_2(\mathcal R )) $ the associated Peetre's $K$-functional, defined
as
$$
  K(s,\phi; X_1(\mathcal R ), X_2(\mathcal R ))=\inf_{\phi=\phi_1+\phi_2}\left(\|\phi_1\|_{X_1(\mathcal R )}+s\|\phi_2\|_{X_2(\mathcal R )}\right)
\qquad \textup{for}\,\, s>0.
$$
We need an expression for the $K$-functional (up to equivalence) in
the case when $X_1(\mathcal R )$ and $X_2(\mathcal R )$ are certain
Lebesgue or Lorentz spaces, and $\mathcal R$ is one of the measure
spaces mentioned above. Recall that, given $\sigma >1$, the Lorentz
space $L^{\sigma , 1}(\mathcal R )$ is the Banach function space of
those measurable functions $\phi$ on $\mathcal R$ for which the norm
$$\|\phi \|_{L^{\sigma , 1}(\mathcal R )} =
\int_0^{\infty}\phi_\nu^* (s) s^{-1+\frac 1\sigma}ds$$ is finite.
The Lorentz space $L^{\sigma , \infty}(\mathcal R )$, also called
Marcinkiewicz space or weak-$L^\sigma$ space, is the Banach function
space of those measurable functions $\phi$ on $\mathcal R$ for which
the quantity
$$\|\phi \|_{L^{\sigma , \infty}(\mathcal R )} =
\sup _{s>0} s^{\frac 1\sigma}\phi_\nu^* (s)$$ is finite. Note that,
in spite of the notation, this is not a norm. However, it is
equivalent to a norm, up to multiplicative constants depending on
$\sigma$, obtained on replacing $\phi_\nu^* (s)$ with $\tfrac 1s
\int_0^s \phi_\nu^* (r)dr$.
\par\noindent It is well known that
\begin{equation}\label{k1}
K(\phi, s; L^1(\mathcal R ), L^\infty (\mathcal R )) = \int _0^s
\phi_\nu^*(r)\, dr \quad \hbox{for $s>0$,}
\end{equation}
for every $\phi \in L^1(\mathcal R) + L^\infty (\mathcal R)$
\cite[Chapter 5, Theorem 1.6]{BS}. If $\sigma >1$, then
\begin{equation}\label{k2}
K(\phi, s; L^{\sigma , \infty}(\mathcal R) , L^\infty (\mathcal R ))
\approx \|r^{\frac 1\sigma} \phi_\nu^*(r)\|_{L^\infty (0, s^\sigma
)} \quad \hbox{for $s>0$,}
\end{equation}
for every $\phi \in L^{\sigma , \infty}(\mathcal R) + L^\infty
(\mathcal R)$ \cite[Equation (4.8)]{H}. Moreover,
\begin{equation}\label{k3}
K(\phi, s; L^{1}(\mathcal R) , L^{\sigma , 1}(\mathcal R)) \approx
\int _0^{s^{\sigma '}} \phi _\nu ^*(r)\, dr + s \int _{s^{\sigma
'}}^\infty r^{-\frac 1{\sigma '}}\phi _\nu^*(r)\, dr \quad \hbox{for
$s>0$.}
\end{equation}
 for every
  $\phi \in L^1(\mathcal R ) + L^{\sigma ,
1}(\mathcal R)$ \cite[Theorem 4.2]{H}.
 In \eqref{k2} and \eqref{k3}, the notation $\lq\lq
\approx "$ means that the two sides are bounded by each other up to
multiplicative constants depending on $\sigma$.
\par\noindent
Let $\mathcal R$ and $\mathcal S$ be positive measure spaces. An
operator $L$ defined on a linear space of measurable functions on
$\mathcal R $, and taking values into the space of measurable
functions in $\mathcal S$, is called sub-linear if,  for every $\phi
_1$ and $\phi_2$ in the domain of $L$ and every $\lambda_1 , \lambda
_2 \in \R$,
$$
|L(\lambda_1\phi_1 + \lambda_2\phi_2)| \leq
|\lambda_1||L\phi_1|+|\lambda_2||L\phi_2|.
$$
A basic result  in the theory of real interpolation tells us what
follows. Assume that
%
$L$   is a sub-linear operator as above, and $X_i(\mathcal R )$ and
$Y_i(\mathcal S )$, $i=1,2$,
 are normed function function spaces on $\mathcal R$ such that
\begin{equation}\label{101}
L : X_i(\mathcal R ) \to Y_i(\mathcal S)
\end{equation}
 with norms not exceeding $N_i$,   $i=1,2$. Here, the arrow $"\to"$ denotes a bounded
 operator.
 Then,
\begin{equation}\label{k4}
K(L\phi, s; Y_1 (\mathcal S), Y_2(\mathcal S)) \leq   \max \{N_1,
N_2\} K(\phi, s; X_1 (\mathcal R), X_2(\mathcal R)) \qquad \hbox{for
$s>0$,}
\end{equation}
for every $\phi\in X_1 (\mathcal R)+ X_2 (\mathcal R)$.

\begin{lemma}\label{lemma1N} Let $\Omega$ be an open set in $\rn$, $n \geq 2$, and let $\gamma \in (0, n)$.
Assume that $\mu$ is any Borel  measure in $\Omega$ fulfilling
\eqref{measure} for some $\alpha \in (n-\gamma , n]$ and for some
$C_\mu
>0$.
 Then there
exists a constant $C=C(n,\alpha , \gamma , C_\mu)$ such that
\begin{equation}\label{lemma1N.1}
\|N_\gamma g\|_{L^{\frac \alpha {n-\gamma }, \infty}(\Omega , \mu)}
\leq C \|g\|_{L^1(\partial \Omega )}
\end{equation}
for every $g \in L^1(\partial \Omega)$.
\end{lemma}
{\bf Proof}.
We make use of an argument related to \cite{Adams1, Adams2}. Given
$g \in L^1(\partial \Omega)$ and  $t>0$, define $E _t = \{x\in
\Omega : |N_\gamma g (x)|
>t\}$, and denote by $\mu _t$ the restriction of the measure $\mu$
to  $E _t$. By Fubini's Theorem, one has that
\begin{align}\label{lemma1.3}
t \mu (\mathcal L _t)  = t \int _\Omega d \mu_t (x) \leq \int
_\Omega |N_\gamma g(x)| d \mu_t (x)
 & \leq
\int _\Omega \int _{\partial \Omega }
\frac{|g(y)|}{|x-y|^{n-\gamma}}\, d\hh (y)\, d \mu_t (x)
\\ \nonumber &\leq  \int _{\partial \Omega} |g(y)|
\int _\Omega  \frac{d \mu _t(x) d\hh (y)}{|x-y|^{n-\gamma}}.
\end{align}
Next,
\begin{align}\label{lemma1.4}
\int _\Omega  \frac{d \mu _t(x)}{|x-y|^{n-\gamma}}\, & =
(n-\gamma)\int _0^\infty \varrho^{-n+\gamma -1} \int _{\{x\in
\Omega: |x-y|^{\gamma -n}
> \varrho^{\gamma -n}\}} \,d \mu _t(x)\,d\varrho \\ \nonumber & \leq
(n-\gamma)\int _0^\infty \varrho^{-n+\gamma -1} \mu _t(B_\varrho
(y))\,d\varrho \quad \hbox{for $t
>0$.}
\end{align}
From \eqref{lemma1.3} and \eqref{lemma1.4} we deduce that, for each
fixed $r >0$,
\begin{align}\label{lemma1.5}
t \mu (\mathcal L _t) & \leq (n-\gamma) \int _{\partial \Omega}
|g(y)|\int _0^\infty \varrho^{-n+\gamma -1} \mu _t(B_\varrho (y))\,d\varrho \, d\hh (y)\\
\nonumber & = (n-\gamma) \int _0^r \varrho^{-n+\gamma -1}  \int
_{\partial
\Omega} |g(y)|\mu _t(B_\varrho (y))\, d\hh (y)\,d\varrho \\
\nonumber & \quad + (n-\gamma) \int _r^\infty \varrho^{-n+\gamma -1}
\int _{\partial \Omega} |g(y)| \mu _t(B_\varrho (y))\, d\hh
(y)\,d\varrho \quad \hbox{for $t
>0$.}
\end{align}
We have that
\begin{multline}\label{lemma1.6}
\int _0^r \varrho^{-n+\gamma -1}  \int _{\partial \Omega} |g(y)|\mu
_t(B_\varrho (y))\, d\hh (y) \,d\varrho \\
\leq C_\mu \|g\|_{L^1(\partial \Omega)} \int _0^r \varrho^{-n+\gamma
-1 +\alpha}\,d\varrho = C_\mu \tfrac {r^{\alpha - n +\gamma}}{\alpha
- n +\gamma} \|g\|_{L^1(\partial \Omega)} .
\end{multline}
On the other hand,
\begin{align}\label{lemma1.7}
\int _r^\infty \varrho^{-n+\gamma -1} \int _{\partial \Omega} |g(y)|
\mu _t(B_\varrho (y))\, d\hh (y) \,d\varrho & \leq \mu (\mathcal L
_t) \int _r^\infty \varrho^{-n+\gamma -1} \int _{\partial \Omega}
|g(y)| \, d\hh (y) \,d\varrho \\ \nonumber & \leq \mu (\mathcal L
_t) \tfrac {r^{- n +\gamma}}{n -\gamma}\|g\|_{L^1(\partial \Omega)}
.
\end{align}
Combining \eqref{lemma1.5}--\eqref{lemma1.7}, and choosing $r=
\Big(\frac{\mu (\mathcal L _t)}{C_\mu}\Big)^{\frac 1\alpha}$, yield
\begin{align}\label{lemma1.8}
t \mu (\mathcal L _t) &\leq (n-\gamma) \|g\|_{L^1(\partial
\Omega)}\Big(C_\mu \frac {r^{\alpha - n +\gamma}}{\alpha - n
+\gamma} +
\mu (\mathcal L _t) \frac {r^{- n +\gamma}}{n -\gamma}\Big) \\
\nonumber & = \tfrac{\alpha }{\alpha -n +\gamma} \|g\|_{L^1(\partial
\Omega)}C_\mu^{\frac{n-\gamma}\alpha} \mu (\mathcal L _t)^{1 -
\frac{n-\gamma}\alpha}.
\end{align}
Thus,
\begin{align}\label{lemma1.9}
t \mu (\mathcal L _t)^{\frac{n-\gamma}\alpha} &\leq \tfrac{\alpha
}{\alpha -n +\gamma}C_\mu^{\frac{n-\gamma}\alpha}\|g\|_{L^1(\partial
\Omega)} \quad \hbox{for $t>0$.}
\end{align}
Hence, inequality \eqref{lemma1N.1} follows. \qed

\begin{proposition}\label{boundary} Let $\Omega$ be an  open
set in $\rn$.
Then
\begin{align}\label{bound1}
\int _{\sn} |g(\zeta(x, \vartheta ))|\, d \hh (\vartheta )  \leq
2^n \int _{\partial  \Omega} \frac{|g(y)|}{|x-y|^{n-1}}d \hh (y)
\quad \hbox{for a.e. $x \in \Omega$,}
\end{align}
for every Borel function $g : \partial \Omega \to \R$.
Here,  convention \eqref{conv} is adopted.
\end{proposition}
{\bf Proof}. We spilt the proof in steps.
\par\noindent
\emph{Step 1}. Denote by $\Pi : \rn \setminus \{x\} \to \sn$ the
projection function into $\mathbb S^{n-1}$ given by
$$\Pi (y) = \frac{y-x}{|y-x|} \quad \hbox{for $y \in \rn \setminus \{x\}$.}$$
Then
\begin{equation}\label{bound2}
\hh (\Pi (E)) \leq \frac 1{{\rm dist}(x, E)^{n-1}} \hh (E)
\end{equation}
for every $E \subset \partial \Omega$.
\par
The function $\Pi$ is differentiable, and $|\nabla \Pi(y)| \leq
|y-x|^{-1}$ for $y \in \rn \setminus \{x\}$. Thus, the restriction
of $\Pi$ to $E$ is Lipschitz continuous, and
\begin{equation}\label{bound2bis} |\nabla \Pi(y)| \leq \frac 1{{\rm
dist}(x, E)} \quad \hbox{for $y \in E$.}
\end{equation}
 Inequality \eqref{bound2bis} implies \eqref{bound2}, by a standard
 property of Hausdorff measure -- see e.g. \cite[Theorem
 7.5]{Mattila}.
\par\noindent
\emph{Step 2}. We have that
\begin{equation}\label{bound5}
\hh (\Pi (E)) \leq 2^n \int _E \frac {d\hh (y)}{|x-y|^{n-1}}
 \quad \hbox{for every Borel set $E \subset \partial
\Omega$.}
\end{equation}
\par
The following chain holds:
\begin{align}\label{bound6}
\int _E \frac {d\hh (y)}{|x-y|^{n-1}}  & = \int _0^\infty \hh
(\{y\in E: |x-y|^{-n+1}>t\})\, dt \\ \nonumber & = \int _0^\infty
\hh (\{y\in E: |x-y|<\tau\})\, d(-\tau ^{1-n})
\\ \nonumber & =
\sum _{k \in \Z} \int _{2^k}^{2^{k+1}} \hh (\{y\in E:
|x-y|<\tau\})\, d(-\tau ^{1-n})
\\ \nonumber & \geq
\sum _{k \in \Z}  \hh (\{y\in E: |x-y|<2^k\})\int _{2^k}^{2^{k+1}}\,
d(-\tau ^{1-n})
\\ \nonumber & =
(1- 2^{-n+1})\sum _{k \in \Z} 2^{-k(n-1)}  \hh (\{y\in E:
|x-y|<2^k\})
\\ \nonumber & \geq
(1- 2^{-n+1}) \sum _{k \in \Z} 2^{-k(n-1)} \hh (\{y\in E:
2^{k-1}\leq  |x-y|<2^k\}).
\end{align}
Since ${\rm dist }(x, \{y\in E: 2^{k-1}\leq  |x-y|<2^k\}) \geq
2^{k-1}$, by \eqref{bound2}
\begin{equation}\label{bound7}
\hh (\{y\in E: 2^{k-1}\leq  |x-y|<2^k\}) \geq C 2^{(k-1)(n-1)} \hh
(\Pi (\{y\in E: 2^{k-1}\leq  |x-y|<2^k\}))
\end{equation}
for $k \in \Z$. From \eqref{bound6} and \eqref{bound7} we deduce
that
\begin{align}\label{bound8}
\int _E \frac {d\hh (y)}{|x-y|^{n-1}} & \geq \tfrac 12 \sum _{k \in
\Z} 2^{-k(n-1)}
2^{(k-1)(n-1)} \hh (\Pi (\{y\in E: 2^{k-1}\leq  |x-y|<2^k\})) \\
\nonumber & = 2^{-n}  \sum _{k \in \Z}  \hh (\Pi (\{y\in E:
2^{k-1}\leq  |x-y|<2^k\}))\\
 \nonumber & \geq  2^{-n}   \hh \big(\cup _{k \in \Z} \Pi (\{y\in E: 2^{k-1}\leq  |x-y|<2^k\})\big)
 \\ \nonumber & =  2^{-n}   \hh \big( \Pi (\cup _{k \in \Z} \{y\in E: 2^{k-1}\leq  |x-y|<2^k\})\big)
\\ \nonumber & =  2^{-n}   \hh ( \Pi (E)).
\end{align}
Inequality \eqref{bound5} is thus established.
\par\noindent
{\emph Step 3}. Conclusion.
\par\noindent
Fix $x \in \Omega$. We have that
$$\Pi (\{y \in (\partial \Omega)_x: |g(y)|>t\}) = \{\vartheta \in \sn : |g(\zeta(x, \vartheta
))|>t\} \quad \hbox{for  $t> 0$.}$$
 Thus, by
\eqref{bound5},
\begin{align*}\int _{\{y\in (\partial \Omega )_x:
|g(y)|>t\}} \frac {d\hh (y)}{|x-y|^{n-1}} & \geq
  2^{-n} \hh (\{\vartheta \in \sn : |g(\zeta(x, \vartheta ))|>t\})
\quad \hbox{for $t> 0$.}
\end{align*}
 Hence,
\begin{align}\label{bound9}
\int _{(\partial \Omega)_x} \frac {|g(y)|}{|x-y|^{n-1}}d\hh (y)& =
\int _0^\infty  \int _{\{y\in (\partial \Omega )_x: |g(y)|>t\}}
\frac {d\hh (y)}{|x-y|^{n-1}}\, dt \\ \nonumber & \geq 2^{-n} \int
_0^\infty \hh (\{\vartheta \in \sn : |g(\zeta(x, \vartheta
))|>t\})\,dt
\\ \nonumber & = 2^{-n} \int _\sn |g(\zeta(x, \vartheta ))|\, d\hh
(\vartheta ).
\end{align}
Inequality \eqref{bound1} is thus established. \qed

\begin{lemma}\label{lemma1} Let $\Omega$ be an open set in $\rn$, $n \geq 2$.
Assume that $\mu$ is any Borel  measure in $\Omega$ fulfilling
\eqref{measure} for some $\alpha \in (n-1 , n]$ and for some $C_\mu
>0$. Then there exists a
constant $C=C(n,\alpha , C_\mu)$ such that
\begin{align}\label{lemma1ineq}
(Tg)^*_\mu(s) \leq C s^{-\frac{n-1}{\alpha}}  \int
_0^{s^{\frac{n-1}{\alpha}}} g_{\hh}^*(r)\, dr \quad \hbox{for
$s>0$,}
\end{align}
%
for every Borel function $g : \partial \Omega \to \R$.
\end{lemma}
\par\noindent
{\bf Proof}. By Proposition \ref{boundary}, there exists a constant
$C=C(n)$ such that
\begin{equation}\label{1.2'}
Tg(x) \leq C N_1|g|(x) \quad \hbox{for a.e. $x \in \Omega$,}
\end{equation}
for every Borel function $g : \partial \Omega \to \R$. Hence, by
Lemma \ref{lemma1N}, with $\gamma =1$, there exists a constant
$C=C(n,\alpha , C_\mu)$ such that
\begin{equation}\label{lemma1T.1}
\|T g\|_{L^{\frac \alpha {n-1 }, \infty}(\Omega , \mu)} \leq C
\|g\|_{L^1(\partial \Omega )}
\end{equation}
for every for every Borel function $g : \partial \Omega \to \R$.
\par\noindent
On the other hand,
\begin{align}\label{lemma1T.2bis}
0 \leq T g (x)  \leq \|g\|_{L^\infty (\partial \Omega )} \int _{\sn}
d\hh (\vartheta) = n\omega _{n} \|g\|_{L^\infty (\partial \Omega )}
 \quad
\hbox{for every $x \in \Omega$,}
\end{align}
 and hence
\begin{align}\label{lemma1T.2}
\|T g\|_{L^{ \infty}(\Omega , \mu)} \leq n\omega _{n}
\|g\|_{L^\infty (\partial \Omega )}
\end{align}
for every Borel function $g : \partial \Omega \to \R$.
 We thus
deduce from \eqref{k1}, \eqref{k2}, \eqref{k4}, \eqref{lemma1T.1}
and \eqref{lemma1T.2} that
\begin{align}\label{k5}
s (Tg)^*_\mu(s^{\frac{\alpha}{n-1}}) & \leq \|r^{\frac {n-1}\alpha}
(Tg)^*_\mu(r)\|_{L^\infty (0, s^{\frac{\alpha}{n-1}} )}
 \approx K(Tg, s;
L^{\frac{\alpha}{n-1} , \infty}(\Omega , \mu) , L^\infty (\Omega ,
\mu ))
 \\
\nonumber & \leq C K(g, s; L^1(\partial \Omega ), L^\infty (\partial
\Omega )) = C \int _0^s g^*_{\hh}(r)\, dr \quad \hbox{for $s>0$,}
\end{align}
for some constant $C=C(n, \alpha, C_\mu )$, and for every  Borel
function $g : \partial \Omega \to \R$. Hence, inequality
\eqref{lemma1ineq} follows. \qed

\bigskip
\par\noindent

\begin{lemma}\label{lemma2} Let $\Omega$ be an open set in $\rn$, $n \geq 2$, and let $\gamma \in (0, n)$.
Assume that $\mu$ is any Borel  measure in $\Omega$ fulfilling
\eqref{measure} for some $\alpha \in (n-\gamma , n]$ and for some
$C_\mu
>0$. Then, there exists a constant
$C=C(n,\gamma, \alpha, C_\mu) $ such that
\begin{align}\label{k7}
(I_\gamma f)_\mu^*(s) \leq C \bigg(s^{-\frac{n-\gamma}{\alpha}} \int
_0^{s^{\frac{n}{\alpha}}} f^*_{\mathcal L^n}(r)\, dr + \int
_{s^{\frac{n}{\alpha}}}^\infty r^{-\frac {n-\gamma}{n}}f
^*_{\mathcal L^n}(r)\, dr\bigg) \quad \hbox{for $s>0$,}
\end{align}
for every $f \in L^1(\Omega)$.
\end{lemma}

\medskip
\par\noindent
{\bf Proof }. A standard weak-type inequality for Riesz potentials
tells us that there exists a constant $C_1=C_1(n, \gamma, \alpha ,
C_\mu)$ such that
\begin{equation}\label{lemma1.1}
\|I _\gamma f\|_{L^{\frac \alpha {n-\gamma}, \infty}(\Omega , \mu)}
\leq C_1 \|f\|_{L^1( \Omega )}
\end{equation}
for every $f \in L^1(\Omega)$ (a proof of inequality
\eqref{lemma1.1} follows, in fact, along the same lines as that of
\eqref{lemma1N.1}).
 Furthermore,
there exists a constant $C_2=C_2(n,\gamma, \alpha, C_\mu)$ such that
\begin{equation}\label{lemma1.2}
\|I_\gamma  f\|_{L^{\infty}(\Omega , \mu )} \leq C_2 \|f\|_{L^{\frac
n\gamma,1}( \Omega )}
\end{equation}
for every $f \in L^{\frac n\gamma,1}( \Omega )$. Inequality
\eqref{lemma1.2}  can be derived from \eqref{lemma1.1}, applied with
$\mu = \mathcal L^n$  and  $\alpha =n$, via a duality argument.
Indeed,
\begin{align}\label{lemma1.2proof}
\|I_\gamma f\|_{L^{\infty}(\Omega )} &\leq \|I_\gamma
|f|\|_{L^{\infty}(\Omega )}
 = \sup _{\|h\|_{L^1(\Omega )}\leq 1} \int
_{\Omega}|h(x)| \int _\Omega \frac{|f(y)|}{|y-x|^{n-\gamma}}\, dy\,
dx
\\ \nonumber & =
\sup _{\|h\|_{L^1(\Omega )}\leq 1} \int _{\Omega}|f(y)| \int _\Omega
\frac{|h(x)|}{|y-x|^{n-\gamma}}\, dy\, dx
 \leq \sup _{\|h\|_{L^1(\Omega )}\leq 1} C \|f\|_{L^{\frac n\gamma , 1}(\Omega
)}\|I_\gamma h\|_{L^{\frac {n}{n-\gamma}, \infty}(\Omega )} \\
\nonumber & \leq \sup _{\|h\|_{L^1(\Omega )}\leq 1} C'
\|f\|_{L^{\frac n\gamma , 1}(\Omega )}\| h\|_{L^{1}(\Omega )} \leq
C' \|f\|_{L^{\frac n\gamma , 1}(\Omega )}
\end{align}
for some constants $C=C(n,\gamma)$ and $C'=C'(n, \gamma, \alpha ,
C_\mu)$, and for every $f \in L^{\frac n\gamma , 1}(\Omega )$. Note
that the first inequality holds owing to a H\"older type inequality
in Lorentz spaces. As shown by a standard convolution argument, the
space of continuous functions is dense in $L^{\frac n\gamma ,
1}(\Omega )$. Inequality \eqref{lemma1.2proof} then  implies that
$I_\gamma f$ is continuous for $f \in L^{\frac n\gamma , 1}(\Omega
)$. Thus, $\|I_\gamma f\|_{L^{\infty}(\Omega , \mu )} \leq
\|I_\gamma f\|_{L^{\infty}(\Omega )}$, and \eqref{lemma1.2} follows
from \eqref{lemma1.2proof}.
\par By \eqref{lemma1.1} and \eqref{lemma1.2}, via
\eqref{k2}, \eqref{k3} and \eqref{k4}, we deduce that there exists a
constant $C=C(n, \gamma, \alpha , C_\mu)$ such that
\begin{align}\label{k6bis}
s &(I_\gamma f)_\mu^*(s^{\frac{\alpha}{n-\gamma}})  \leq \|r^{\frac
{n-\gamma}\alpha} (I_\gamma f)^*_\mu(r)\|_{L^\infty (0,
s^{\frac{\alpha}{n-\gamma}} )}
 \approx K(I_\gamma f, s;
L^{\frac{\alpha}{n-\gamma} , \infty}(\Omega , \mu) , L^\infty
(\Omega , \mu ))\\ \nonumber &  \leq C K(f, s; L^{1}(\Omega) ,
L^{\frac n\gamma , 1}(\Omega)) \approx \int _0^{s^{\frac
n{n-\gamma}}} f ^*_{\mathcal L^n}(r)\, dr + s \int _{s^{\frac
n{n-\gamma}}}^\infty r^{-\frac {n-\gamma}n}f ^*_{\mathcal L^n}(r)\,
dr \quad \hbox{for $s>0$,}
\end{align}
where the equivalence is up to multiplicative constants depending on
$n, \gamma, \alpha , C_\mu$. Hence, \eqref{k7} follows. \qed

\begin{lemma}\label{lemma3Q}  Let $\Omega$ be an open set in $\rn$, $n \geq 2$, and let $\gamma \in (0, n-1)$.
Assume that $\mu$ is any Borel  measure in $\Omega$ fulfilling
\eqref{measure} for some $\alpha \in (n-\gamma , n]$ and for some
$C_\mu
>0$.  Then there exists a
constant $C=C(n,\gamma, \alpha , C_\mu)$ such that
\begin{align}\label{lemma1ineq'}
(Q_\gamma g)^*_\mu(s) \leq C \bigg(s^{-\frac{n-1-\gamma}{\alpha}}
\int _0^{s^{\frac{n-1}{\alpha}}} g^*_{\hh}(r)\, dr + \int
_{s^{\frac{n-1}{\alpha}}}^\infty r^{-\frac {n-1-\gamma}{n-1}}g
^*_{\hh}(r)\, dr\bigg) \quad \hbox{for $s>0$,}
\end{align}
%
for every  Borel function $g : \partial \Omega \to \R$.
\end{lemma}
\par\noindent
{\bf Proof}. By inequality \eqref{lemma1T.1}, with $\mu = \mathcal
L^n$, there exists a constant $C=C(n)$ such that
\begin{equation}\label{lemma3Q.1}
\|T g\|_{L^{\frac n {n-1 }, \infty}(\Omega )} \leq C
\|g\|_{L^1(\partial \Omega )}
\end{equation}
for every  Borel function $g : \partial \Omega \to \R$. Moreover,
there exists a constant $C=C(n,\gamma, \alpha , C_\mu)$ such that
\begin{equation}\label{lemma3Q.2}
\|I_\gamma f\|_{L^{\frac \alpha{n-1-\gamma}, \infty}(\Omega , \mu)}
\leq C \|f\|_{L^{\frac n {n-1}, \infty}(\Omega )}
\end{equation}
for every $f \in L^{\frac n {n-1 }, \infty}(\Omega )$. Indeed, by
\eqref{k7}, for any such $f$,
\begin{align}\label{lemma3Q.2proof}
\|I_\gamma f\|_{L^{\frac \alpha{n-1-\gamma}, \infty}(\Omega , \mu)}
& = \sup _{s >0} s^{\frac{n-\gamma -1}\alpha} (I_\gamma f)^*_\mu (s)
\\ \nonumber & \leq C \sup _{s >0}\bigg( s^{-\frac 1\alpha}
\int _0^{s^{\frac{n}{\alpha}}} f^*_{\mathcal L^n}(r)\, dr +
s^{\frac{n-\gamma -1}\alpha} \int _{s^{\frac{n}{\alpha}}}^\infty
r^{-\frac {n-\gamma}{n}}f ^*_{\mathcal L^n}(r)\, dr\bigg)
\\ \nonumber & \leq C \|f\|_{L^{\frac n {n-1}, \infty}(\Omega )} \sup _{s >0}\bigg( s^{-\frac 1\alpha}
\int _0^{s^{\frac{n}{\alpha}}} r^{-\frac{n-1}n}\, dr +
s^{\frac{n-\gamma -1}\alpha} \int _{s^{\frac{n}{\alpha}}}^\infty
r^{-\frac{n-1}n-\frac {n-\gamma}{n}}\, dr\bigg)
\\ \nonumber & = C' \|f\|_{L^{\frac n {n-1}, \infty}(\Omega )},
\end{align}
where $C$ is the constant appearing in \eqref{k7}, and
$C'=C'(n,\gamma, \alpha , C_\mu)$ . If follows from
\eqref{lemma3Q.1} and \eqref{lemma3Q.2} that
\begin{equation}\label{lemma3Q.3}
\|Q_\gamma g\|_{L^{\frac \alpha{n-1-\gamma}, \infty}(\Omega , \mu)}
\leq C \|g\|_{L^1(\partial \Omega )},
\end{equation}
for some constant $C=C(n,\gamma, \alpha , C_\mu)$, and for every
Borel function $g : \partial \Omega \to \R$.
\par\noindent On the other hand, by \eqref{lemma1ineq}, applied with
$\mu = \mathcal L^n$, there exists a constant $C=C(n, \gamma)$ such
that
\begin{align}\label{lemma3Q.4}
\|Tg\|_{L^{\frac n\gamma , 1}(\Omega )} & = \int _0^\infty
(Tg)^*_{\mathcal L^n} (s)s^{-1+\frac \gamma n}\, ds \leq C \int
_0^\infty s^{-1+\frac \gamma n-\frac{n-1}{n}} \int
_0^{s^{\frac{n-1}{n}}} g^*_{\hh}(r)\, dr\, ds
\\ \nonumber & = \int _0^\infty g^*_{\hh}(r) \int _{r^{\frac
n{n-1}}}^\infty s^{-1+\frac \gamma n-\frac{n-1}{n}} \, ds\, dr
 = \tfrac n{n-\gamma -1}\int _0^\infty g^*_{\hh}(r) r^{-1+\frac \gamma {n-1}}\, dr
 \\ \nonumber & = \tfrac n{n-\gamma -1}\|g\|_{L^{\frac {n-1}\gamma ,
1}(\partial \Omega)}
\end{align}
for every  Borel function $g : \partial \Omega \to \R$. Coupling
inequalities  \eqref{lemma3Q.4} and \eqref{lemma1.2} tells us that
there exists a constant $C=C(n,\gamma, \alpha , C_\mu)$ such that
\begin{equation}\label{lemma3Q.5}
\|Q_\gamma g\|_{L^{\infty}(\Omega , \mu)} \leq C \|g\|_{L^{\frac
{n-1}\gamma , 1}(\partial \Omega)}
\end{equation}
for every  Borel function $g : \partial \Omega \to \R$.
\par\noindent
Now,  by \eqref{k2}, \eqref{k3}, \eqref{k4}, \eqref{lemma3Q.3} and
\eqref{lemma3Q.5}, there exists a constant $C=C(n,\gamma, \alpha ,
C_\mu)$ such that
\begin{multline*}
s (Q_\gamma g)^*_\mu(s^{\frac{\alpha}{n- 1- \gamma}})  \leq
\|r^{\frac {n- 1- \gamma}\alpha} (Q_\gamma g)^*_\mu (r)\|_{L^\infty
(0, s^{\frac{\alpha}{n-\gamma - 1}} )}
 \approx K(Q_\gamma g, s;
L^{\frac{\alpha}{n- 1- \gamma} , \infty}(\Omega , \mu) , L^\infty
(\Omega , \mu ))\\   \leq C K(g, s; L^{1}(\partial \Omega) ,
L^{\frac {n-1}\gamma , 1}(\partial \Omega)) \approx \int
_0^{s^{\frac {n-1}{n- 1- \gamma}}} g ^*_{\hh}(r)\, dr + s \int
_{s^{\frac {n-1}{n- 1- \gamma}}}^\infty r^{-\frac {n- 1-
\gamma}{n-1}}g ^*_{\hh}(r)\, dr\quad \hbox{for $s>0$,}
\end{multline*}
 for every  Borel function $g : \partial \Omega \to \R$,
where equivalence holds up to multiplicative constants depending on
$n, \gamma, \alpha, C_\mu$. Inequality \eqref{lemma1ineq'} follows.
\qed

\medskip
\par\noindent
{\bf Proof of Theorem \ref{intermest}}. Inequality \eqref{hfund1ell}
can be written as
$$|\nabla ^h u(x)| \leq C\bigg(I_{m-h}(|\nabla ^m u|)(x) + \sum
_{k=1}^{m - h -1} Q_k\big(g^{[\frac{k+h+1}{2}], \natural ( k+h
)}\big)(x) + T \big(g^{[\frac{h+1}{2}], \natural ( h
)}\big)(x)\bigg) \quad \hbox{for a.e. $x \in \Omega$.}$$
 Hence, \eqref{rearresthell}  follows via Lemmas \ref{lemma1} --
 \ref{lemma3Q},
 owing to
property \eqref{f+g} of rearrangements. \qed

\section{Sobolev inequalities}\label{ineq}

We present here a sample of Sobolev type inequalities that can be
established via  the universal pointwise and rearrangement estimates
of Sections \ref{proofs} and \ref{sobolev}, respectively. We limit
ourselves to  inequalities for standard norms, such as Lebesgue
norms and Orlicz norms of exponential or logarithmic type, which
naturally come into play in borderline situations. Measures $\mu$
satisfying \eqref{measure} will be included in our results. Let us
emphasize, however, that inequalities for more general norms can be
derived from the relevant pointwise bounds. Virtually, any Sobolev
type inequality for rearrangement-invariant norms, which holds in
regular domains, has a counterpart in arbitrary domains, provided
that appropriate boundary seminorms are employed.
 \par A key tool in our approach is the reduction principle to
 one-dimensional inequalities stated in Theorem \ref{reduction}
 below for Sobolev inequalities involving arbitrary
 rearrangement-invariant norms.
 Recall that a
rearrangement-invariant space $X(\mathcal R)$ on a measure space
$\mathcal R$, endowed with a positive measure $\nu$, is a Banach
function space (in the sense of Luxemburg) endowed with a norm
$X(\mathcal R)$ such that
\begin{equation}\label{ri1}
\|\phi\|_{X(\mathcal R)} = \|\psi\|_{X(\mathcal R)}\quad
\hbox{whenever} \quad \phi^*_\nu = \psi^*_\nu.
\end{equation}
Every rearrangement-invariant space $X(\mathcal R)$ admits a
representation space $\overline X(0, \infty)$, namely another
rearrangement-invariant space  on $(0, \infty)$ such that
\begin{equation}\label{ri2}
\|\phi\|_{X(\mathcal R)} = \|\phi^*_\nu\|_{\overline X(0,
\infty)}\quad \hbox{for every $\phi \in X(\mathcal R)$.}
\end{equation}
In customary situations, an expression for the norm
$\|\cdot\|_{\overline X(0, \infty)}$ immediately follows from that
of $\|\cdot\|_{X(\mathcal R)}$. The Lebesgue spaces and the Lorentz
spaces, whose definition has been recalled above, are standard
instances of rearrangement-invariant spaces. The exponential spaces,
which have already been mentioned in Section \ref{over}, can be
regarded as special examples of Orlicz spaces. The Orlicz space
$L^A(\mathcal R)$ built upon a Young function $A: [0, \infty) \to
[0, \infty]$, namely  a left-continuous convex function which is
neither identically equal to $0$ nor to $\infty$, is a
rearrangement-invariant space equipped the Luxemburg norm given by
\begin{equation}\label{Orlicz}
\| \phi \|_{_{L^{A}(\mathcal R)}} =  \ \inf \Bigg\{ \lambda > 0 \ :
\ \int_{\mathcal R} \bigg(\frac{|\phi (x)|}{\lambda} \bigg)\, dx \
\leq \ 1 \Bigg\}.
\end{equation}
The class of Orlicz spaces includes that of Lebesgue spaces, since
$L^{A}(\mathcal R)= L^{p}(\mathcal R)$ if $A(t)=t^p$ for $p \in [1,
\infty[$, and $L^{A}(\mathcal R)= L^{\infty}(\mathcal R)$ if
$A(t)=\infty \chi_{_{(1, \infty)}}(t)$. Given $\sigma
>0$, we denote by $\exp L^\sigma (\mathcal R)$ the Orlicz space
built upon  the Young function $A(t)=e^{t^\sigma } -1$, and by
$L^p(\log L)^\sigma (\mathcal R)$ the Orlicz space built upon the
Young function $A(t) = t^p \log ^\sigma (c+t)$, where $c$ is a
sufficiently large positive number.
\par\noindent We refer to \cite{BS} for a comprehensive account of
rearrangement-invariant spaces.

\begin{theorem}\label{reduction} {\bf [Reduction principle for Sobolev inequalities]}
Let $\Omega$ be any  open set in $\rn$, $n \geq 2$. Assume that
$\mu$ is a measure in $\Omega$ fulfilling \eqref{measure} for some
$\alpha \in (n-1, n]$, and for some constant $C_\mu$. Let $m \in
\N$, and $h \in \N _0$ be such that
 $0<m-h < n$. Assume that $X(\Omega)$, $Y(\Omega , \mu)$ and
 $X_k(\partial \Omega)$, $k=0, \cdots, m-h-1$, are rearrangement-invariant spaces such that
\begin{equation}\label{red1}
\left\| s^{-\frac {n-m +h)}{\alpha}}\int _0^{s^{\frac
n\alpha}}\varphi(r) dr \right\|_{\overline Y(0, \infty )} \leq C
\|\varphi\|_{\overline X(0, \infty )},
\end{equation}
\begin{equation}\label{red2}
\left\| \int _{s^{\frac n\alpha}}^\infty r^{-\frac{n-m +h)}n}
\varphi(r) dr \right\|_{\overline Y(0, \infty )} \leq C
\|\varphi\|_{\overline X(0, \infty )},
\end{equation}
\begin{equation}\label{red3}
\left\|s^{-\frac {n-k-1}{\alpha}}\int _0^{s^{\frac {n-1}\alpha}}
\varphi(r)dr \right\|_{\overline Y(0, \infty )} \leq C
\|\varphi\|_{\overline X_k(0, \infty )}, \quad k=1, \dots , m - h
-1,
\end{equation}
\begin{equation}\label{red4}
\left\| \int_{s^{\frac {n-1}\alpha}}^\infty  r^{-\frac {n-k-1}{n-1}}
\varphi(r)dr \right\|_{\overline Y(0, \infty )} \leq C
\|\varphi\|_{\overline X_k(0, \infty )}, \quad k=1, \dots , m - h
-1,
\end{equation}
\begin{equation}\label{red5}
\left\| s^{-\frac {n-1}{\alpha}}\int _0^{s^{\frac {n-1}\alpha}}
\varphi(r)\, dr \right\|_{\overline Y(0, \infty )} \leq C
\|\varphi\|_{\overline X_0(0, \infty )},
\end{equation}
for some constant $C$, and for every non-increasing function
$\varphi : [0, \infty) \to [0, \infty)$.   Then
\begin{align}\label{red6}
\|\nabla ^hu\|_{Y(\Omega, \mu)} & \leq C' \Big(\|\nabla ^{m}
u\|_{X(\Omega )} +  \sum _{k=0}^{m -h-1} \|u\|_{\mathcal V
^{[\frac{k+h+1}{2}], \natural ( k+h )}X_k(\partial \Omega )}\Big)
\end{align}
 for every $u \in {V^{m}X(\Omega )
\cap C^{[\frac{m -1}2]}_{\rm b}(\overline \Omega )}$, for some
constant $C'=C'(n, C)$.
\end{theorem}

\begin{remark}\label{remarkred}
{\rm The statement of Theorem \ref{reduction} can be somewhat
generalized, in the sense that assumptions
\eqref{red1}--\eqref{red5} can be weakened if either $\mu (\Omega )<
\infty$, or $\mathcal L^n (\Omega ) < \infty$, or  $\hh (\partial
\Omega ) < \infty$. Specifically: if $\mu (\Omega )< \infty$, it
suffices to assume that there exists $L \in (0, \infty)$ such that
inequalities \eqref{red1}--\eqref{red5} hold with the integral
operators multiplied by $\chi_{(0, L)}$ on the left-hand sides;  if
$\mathcal L^n (\Omega ) < \infty$, it suffices to assume that
inequalities \eqref{red1}--\eqref{red2} hold with $\varphi$ replaced
by $\varphi\chi_{(0, M)}$ for some $M \in (0, \infty)$; if $\hh
(\partial \Omega ) < \infty$, it suffices to assume that
inequalities \eqref{red3}--\eqref{red5} hold with $\varphi$ replaced
by $\varphi\chi_{(0, N)}$ for some $N \in (0, \infty)$. Then
inequality \eqref{red6} holds, but with $C'$ depending also on
either on $L$ and $\mu (\Omega )$, or on $M$ and $\mathcal L^n
(\Omega )$, or on $N$ and $\hh (\partial \Omega ) < \infty$,
according to whether $\mu (\Omega )< \infty$, or $\mathcal L^n
(\Omega ) < \infty$, or $\hh (\partial \Omega ) < \infty$. }
\end{remark}



%
%

Our first application of Theorem \ref{reduction} yields  the
following Sobolev type inequality, in arbitrary domains, with usual
exponents.

\begin{theorem}\label{mainmeasure} {\bf [Sobolev inequality with measure]}
Let $\Omega$ be any  open set in $\rn$, $n \geq 2$. Assume that
$\mu$ is a measure in $\Omega$ fulfilling \eqref{measure} for some
$\alpha \in (n-1, n]$, and for some constant $C_\mu$. Let $m \in
\N$, and $h \in \N _0$ be such that
 $0<m-h < n$. If $1 < p <
\frac{n}{m -h}$, then
 there
 exists a constant $C=C(n,m, p, \alpha, C_\mu)$ such that
\begin{align}\label{mainlebdisp}
\|\nabla ^hu\|_{L^{\frac {\alpha p}{n-(m -h)p}}(\Omega, \mu)} & \leq
C \Big(\|\nabla ^{m} u\|_{L^p(\Omega )} +  \sum _{k=0}^{m -h-1}
\|u\|_{\mathcal V ^{[\frac{k+h+1}{2}], \natural ( k+h
)}L^{\frac{p(n-1)}{n-(m - h -k)p}}(\partial \Omega )}\Big)
\end{align}
for every $u \in {V^{m,p}(\Omega ) \cap C^{[\frac{m -1}2]}_{\rm
b}(\overline \Omega )}$.
\end{theorem}

The next result tells us that, as in the classical Rellich theorem,
the Sobolev embedding corresponding to inequality
\eqref{mainlebdisp} is pre-compact  if the exponent $\tfrac {\alpha
p}{n-(m -h)p}$ is replaced with any smaller one, and $\mu (\Omega) <
\infty$.

\begin{theorem}\label{maincompact} {\bf [Compact Sobolev embedding with measure]}
Let $\Omega$, $\mu$, $n$, $m$ and $h$  be as in Theorem
\ref{mainmeasure}. Assume, in addition, that $\mu (\Omega ) <
\infty$.  If $1 \leq q<\frac {\alpha p}{n-(m -h)p}$, and $\{u_i\}$
is a bounded sequence in $V^{m,p}(\Omega ) \cap
 C^{[\frac{m
-1}2]}_{\rm b}(\overline \Omega )$ endowed with the norm appearing
on the right-hand side of \eqref{mainlebdisp},  then $\{\nabla ^h
u_i\}$ is a Cauchy sequence in $L^q(\Omega, \mu)$.
\end{theorem}

The limiting case when $p= \tfrac n{m-h}$, which is excluded from
Theorem \ref{mainmeasure} , is considered in the next statement,
which provides us with a Yudovich-Pohozaev-Trudinger type inequality
in arbitrary domains.

\begin{theorem}\label{trudmeasure} {\bf [Limiting Sobolev inequality with measure]}
Let $\Omega$ and $\mu$ be as in Theorem \ref{mainmeasure} . Assume,
in addition, that $\mathcal L^n (\Omega)< \infty$, $\mu (\Omega)<
\infty$ and $ \hh (\partial \Omega) < \infty$. Let  $m \in \N$ and
$h \in \N _0$ be such that
 $0<m-h < n$.
Then there
 exists a constant $C=C(n,m, \alpha, C_\mu, \mathcal L^n (\Omega), \mu(\Omega), \hh (\partial \Omega))$ such that
\begin{align}\label{truddisp}
\|\nabla ^h u\|_{\exp L^{\frac {n}{n-(m -h)}}(\Omega, \mu)} & \leq
C\Big( \|\nabla ^{m} u\|_{L^{\frac n{m -h}}(\Omega )} + \sum
_{k=1}^{m -h-1} \|u\|_{\mathcal V ^{[\frac{k+h+1}{2}], \natural (
k+h ) }L^{\frac{n-1}{k}}({\rm log}L)^{\frac{(m
-h)(n-k-1)}{nk}}(\partial \Omega )}
\\ \nonumber & \quad \quad  +
 \|u\|_{\mathcal V^{[\frac{h+1}{2}], \natural ( h )}\exp L^{\frac
{n}{n-(m -h)}}(\partial \Omega)}\Big)
\end{align}
for every $u \in {V^{m,\frac{n}{m -h}}(\Omega ) \cap C^{[\frac{m
-1}2]}_{\rm b}(\overline \Omega )}$.
\end{theorem}

The super-limiting regime, where $p> \tfrac n{m -h}$ is the object
of the following theorem.

\begin{theorem}\label{inf} {\bf [Super-limiting Sobolev inequality]}
Let $\Omega$ be a  open set in $\rn$, $n \geq 2$, such that
 $\mathcal L^n (\Omega)< \infty$  and $ \hh (\partial \Omega) < \infty$. Assume that $m \in
\N$,  $h \in \N _0$,  and $0< m-h< n$. If $p
>\frac n{m -h}$ and $p_k > \frac {n-1}k$ for $k= 1, \dots , m - h
-1$, then
 there
 exists a  constant $C=C(n,m, p, p_1, \dots , p_{m-h-1}, \mathcal L^n (\Omega),  \hh (\partial
 \Omega))$
  such that
\begin{align}\label{infdisp}
\|\nabla ^h u\|_{L^\infty (\Omega )} & \leq C \Big(\|\nabla ^{m}
u\|_{L^{p}(\Omega )}
 +  \sum
_{k=1}^{m -h-1}\|u\|_{\mathcal V^{[\frac{k+h+1}{2}], \natural ( k+h
) } L^{p_k}(\partial \Omega )} +  \|u\|_{\mathcal
V^{[\frac{h+1}{2}], \natural ( h )} L^\infty(\partial \Omega)}\Big)
\end{align}
for every $u \in {V^{m,p}(\Omega ) \cap C^{[\frac{m -1}2]}_{\rm
b}(\overline \Omega )}$.
\end{theorem}

%
%
%
%
%

\medskip
\par\noindent
{\bf Proof of Theorem \ref{mainmeasure}}. Since, for any
 measure space $\mathcal R$, a representation space
of the Lebesgue space $L^p(\mathcal R)$ is just $L^p(0, \infty)$,
the conclusion can be easily deduced from Theorem \ref{reduction},
via standard one-dimensional Hardy type inequalities for Lebesgue
norms (see e.g. \cite[Section 1.3.2]{Mabook}). \qed

\medskip
\par\noindent
{\bf Proof of Theorem \ref{maincompact}}. Fix any $\varepsilon >0$.
Then, there exists a compact set $K\subset \Omega$ such that $\mu
(\Omega \setminus K)< \varepsilon$.  Let $\varrho \in C_0^\infty
(\Omega)$ be such that $0 \leq \varrho \leq 1$, $\varrho = 1$ in
$K$. Thus, $K \subset  {\rm supp} (\varrho)$, the support of
$\varrho$, and hence
\begin{equation}\label{comp0}
\mu ({\rm supp} (1-\varrho )) \leq \mu(\Omega \setminus K)<
\varepsilon.
\end{equation}
 Let $\Omega '$ be an open set, with a smooth boundary,
such that ${\rm supp}(\varrho) \subset \Omega ' \subset \Omega$.
 Let $\{u_i\}$ be a bounded sequence in
${V^{{m},p}(\Omega ) \cap C^{[\frac {m -1}2]}_{\rm b}(\overline
\Omega )}$.
Then, by Theorem \ref{mainmeasure} (applied with $\mu = \mathcal
L^n$), it is also bounded in the standard Sobolev space
$W^{{m},p}(\Omega ')$. By a weighted version of Rellich's
compactness theorem \cite[Theorem 1.4.6/1]{Mabook},
$\{\nabla ^h u_i\}$ is a Cauchy sequence in $L^{q}(\Omega ', \mu)$,
and hence there exists $i_0 \in \N$ such that
\begin{equation}\label{comp1}
\|\nabla ^h u_i - \nabla ^hu_j\|_{L^q(\Omega ', \mu)}<\varepsilon
\end{equation}
if $i, j > i_0$. On the other hand, by H\"older's inequality,
\begin{align}\label{comp2}
\|(1-\varrho ) (\nabla ^hu_i - \nabla ^hu_j)\|_{L^q(\Omega, \mu)} &
\leq \|\nabla ^h u_i - \nabla ^h u_j\|_{L^{\frac {\alpha
p}{n-{m}p}}(\Omega, \mu)} \mu ({\rm supp} (1-\varrho
))^{\frac{\alpha p - (n-mp)q}{\alpha pq}}
\\ \nonumber
& \leq C \big(\|u_i\|_{{V^{{m},p}(\Omega ) \cap C^{[\frac{m
-1}2]}(\overline \Omega )}} + \|u_j\|_{{V^{{m},p}(\Omega ) \cap
C^{[\frac{m -1}2]}(\overline \Omega )}}\big)
\varepsilon^{\frac{\alpha p
- (n-mp)q}{\alpha pq}}\\
\nonumber & \leq C'\varepsilon^{\frac{\alpha p - (n-mp)q}{\alpha
pq}}
\end{align}
for some constants $C$ and $C'$ independent of $i$ and $j$. From
\eqref{comp1} and \eqref{comp2} we infer that
\begin{equation}\label{comp3}
\|\nabla ^h u_i - \nabla ^h u_j\|_{L^q(\Omega , \mu)} \leq \|\nabla
^h u_i - \nabla ^h u_j\|_{L^q(\Omega ', \mu)} +  \|(1-\varrho
)(\nabla ^h u_i - \nabla ^h u_j)\|_{L^q(\Omega, \mu)} \leq
\varepsilon + C'\varepsilon^{\frac{\alpha p - (n-mp)q}{\alpha pq}}
\end{equation}
if $i, j > i_0$. Owing to the arbitrariness of $\varepsilon$,
inequality \eqref{comp3} tells us that $\{\nabla ^hu_i\}$ is a
Cauchy sequence in $L^{q}(\Omega , \mu)$. \qed

\medskip
\par\noindent
{\bf Proof of Theorem \ref{trudmeasure}}. If $\mathcal R$ is a
finite measure space , then the norm of a function $\phi$ in the
Orlicz space $\exp L^\sigma (\mathcal R)$, with $\sigma
>0$, is equivalent, up to multiplicative constants depending on
$\sigma$ and $\nu (\mathcal R)$, to the functional
$$
\big\|\big(1+ \log \tfrac {\nu(\mathcal R)}{s}\big)^{-\frac 1\sigma}
\phi _\nu ^*(s)\big\|_{L^\infty (0, \nu (\mathcal R))}.$$
 Moreover, the norm in the Orlicz space $L^p \log ^\sigma L
 (\mathcal R)$ is equivalent, up to multiplicative constants depending
 on $p$,
$\sigma$ and $\nu (\mathcal R)$ to the functional
$$
\big\|\big(1+ \log \tfrac {\nu(\mathcal R)}{s}\big)^{\frac \sigma p}
\phi_\nu ^*(s)\big\|_{L^p (0, \nu (\mathcal R))}.$$ Thus, owing to
Theorem \ref{reduction} and Remark \ref{remarkred}, inequality
 \eqref{truddisp} will follow if we
show that
\begin{equation}\label{prooftrud1}
\left\| s^{-\frac {n-(m - h)}{\alpha}}\big(1+ \log \tfrac {\mu
(\Omega)}s\big)^{-\frac{n-(m -h)}{n}}\int _0^{s^{\frac
n\alpha}}\varphi(r) dr \right\|_{L^{\infty}(0, \mu (\Omega) )} \leq
C \|\varphi\|_{L^{\frac n{m -h}}(0, \mathcal L^n (\Omega) )},
\end{equation}
\begin{equation}\label{prooftrud2}
\left\| \big(1+ \log \tfrac {\mu (\Omega)}s\big)^{-\frac{n-(m
-h)}{n}}\int _{s^{\frac n\alpha}}^\infty r^{-\frac{n-(m -h)}n}
\varphi(r) dr \right\|_{L^{\infty}(0, \mu (\Omega) )} \leq C
\|\varphi\|_{L^{\frac n{m -h}}(0, \mathcal L^n(\Omega) )},
\end{equation}
for every non-increasing function $\varphi: [0, \infty) \to [0,
\infty)$ with support in $[0, \mathcal L^n(\Omega)]$, and
\begin{multline}\label{prooftrud3}
\left\|s^{-\frac {n-k-1}{\alpha}}\big(1+ \log \tfrac {\mu
(\Omega)}s\big)^{-\frac{n-(m -h)}{n}}\int _0^{s^{\frac {n-1}\alpha}}
\varphi(r)dr \right\|_{L^{\infty}(0, \mu (\Omega) )} \\ \leq C
\|\big(1+ \log \tfrac {\hh(\partial \Omega)}s\big)^{\frac{(m
-h)(n-k-1)}{n(n-1)}}\varphi(s)\|_{L^{\frac{n-1}{k}}(0, \hh(\partial
\Omega) )}, \quad k=1, \dots , m -h -1,
\end{multline}
\begin{multline}\label{prooftrud4}
\left\|\big(1+ \log \tfrac {\mu (\Omega)}s\big)^{-\frac{n-(m
-h)}{n}} \int_{s^{\frac {n-1}\alpha}}^\infty r^{-\frac
{n-k-1}{n-1}} \varphi(r)dr \right\|_{L^{\infty}(0, \mu (\Omega) )} \\
C \|\big(1+ \log \tfrac {\hh(\partial \Omega)}s\big)^{\frac{(m
-h)(n-k-1)}{n(n-1)}}\varphi(s)\|_{L^{\frac{n-1}{k}}(0, \hh(\partial
\Omega) )}, \quad k=1, \dots , m -h -1,
\end{multline}
\begin{multline}\label{prooftrud7}
\left\| s^{-\frac {n-1}{\alpha}}\big(1+ \log \tfrac {\mu
(\Omega)}s\big)^{-\frac{n-(m -h)}{n}}\int _0^{s^{\frac {n-1}\alpha}}
\varphi(r)\, dr \right\|_{L^{\infty}(0, \mu (\Omega) )} \\ \leq C
\|\big(1+ \log \tfrac {\hh(\partial \Omega)}s\big)^{-\frac{n-(m
-h)}{n}}\varphi(s)\|_{L^{\infty}(0, \hh(\partial \Omega) )},
\end{multline}
for some constant $C$ and every non-increasing function $\varphi:
[0, \infty) \to [0, \infty)$ with support in $[0, \hh(\partial
\Omega)]$.
\\
Inequalities \eqref{prooftrud1}--\eqref{prooftrud7} are consequences
of classical weighted Hardy type inequalities  (\cite[Section
1.3.2]{Mabook}). \qed

\medskip
\par\noindent
{\bf Proof of Theorem \ref{inf}}. Inequality \eqref{infdisp} follows
from Theorem \ref{reduction} and Remark \ref{remarkred}, via
weighted Hardy type inequalities  (\cite[Section 1.3.2]{Mabook}).
 \qed



\section{Sharpness of results}\label{sharp}

In this section we work out in detail some examples,  announced in
Sections \ref{sec1} and \ref{over}, in connection with certain
sharpness features of the inequalities presented above.

\begin{example}\label{ex3}
{\rm We observed in Section \ref{over} that the term
$\|u\|_{\mathcal V^{1,0}L^{\frac{p(n-1)}{n-p}}(\partial \Omega )}$
can be dropped on the right-hand side of \eqref{mainsecond0} if
$\Omega$ is a regular domain. Here, we show that, by contrast, the
term in question is indispensable for an arbitrary domain. To this
purpose, we exhibit a domain $\Omega \subset \rn$ for which the
inequality
\begin{align}\label{mainsecondex}
\| u\|_{L^{\frac{pn}{n-2p}}(\Omega)}  \leq C \big(\|\nabla ^{2}
u\|_{L^p(\Omega )} + \|u\|_{L^{\frac{p(n-1)}{n-2p}}(\partial \Omega
)}\big)
\end{align}
fails for $1<p<\tfrac n2$, for every constant $C$ independent of
$u$. The relevant domain is the union of a sequence of axially
symmetric \lq\lq cusp-shaped " subdomains $\Omega _k$ about the
$x_n$-axis, which are connected by thin cylinders $H_k$ joining the
vertex of  $\Omega _{k}$ with the basis of  $\Omega _{k-1}$ (Figure
1, Section \ref{over}). Each subdomain $\Omega _k$ is the set of
revolution about the $x_n$-axis of the form
$$\Omega _k = \big\{x: |x'|< (x_n - x_n^k+\varepsilon _k)^\beta , x_n \in (x_n^k, x_n^k + h_k)\big\}$$
for some $x_n^k >0$ and $0< \varepsilon_k < h_k$. The cylinder $H_k$
has a basis of radius $\varepsilon_k^\beta$. Define the sequence
$\{u_k\}$ by
$$u_k(x) = 1-\tfrac{x_n - x_n^k}{h_k} \quad \hbox{for $x \in \Omega
_k$,}$$ $u_k =0$ in $\Omega _j$ for $j \neq k$ and in $H_j$ for $j
\neq k, k+1$, and is continued to $H_{k}$ and $H_{k+1}$ in such a
way that $u \in C^2(\Omega)$.
%
%
\par\noindent
One can verify that
\begin{equation}\label{ex31}
\| u_k\|_{L^{\frac{pn}{n-2p}}(\Omega)} \approx
h_k^{\frac{[(n-1)\beta + 1](n-2p)}{np}},
\end{equation}
\begin{equation}\label{ex32}
\| u_k\|_{L^{\frac{p(n-1)}{n-2p}}(\partial \Omega)} \approx
h_k^{\frac{[(n-2)\beta + 1](n-2p)}{(n-1)p}},
\end{equation}
as $k \to \infty$, and
\begin{equation}\label{ex33}
\| \nabla ^2 u_k\|_{L^{p}(\Omega)}= \| \nabla ^2
u_k\|_{L^{p}(H_{k}\cup H_{k+1})}
\end{equation}
for $k \in \N$.
%
If $\varepsilon _k$ decays to $0$ sufficiently fast as  $k \to
\infty$, the norm $ \| \nabla ^2 u_k\|_{L^{p}(H_{k}\cup H_{k+1})}$
decays arbitrarily fast to $0$. Thus, inequality
\eqref{mainsecondex} fails when tested on the sequence $u_k$,
whatever $C$ is.}
\end{example}

\begin{example}\label{ex1}
{\rm Our purpose here is to demonstrate that, whereas the seminorm
$\| u\|_{\mathcal V ^{1,0}L^{r}(\partial \Omega )} $ can be replaced
with $\| u\|_{L^{r}(\partial \Omega )}$ in \eqref{mainsecond1} when
$\Omega$ is a regular domain, this is impossible, in general, if no
regularity on $\Omega$ is retained. Precisely, we construct an
 open set $\Omega$ in $\mathbb R^2$ for which the
inequality
\begin{equation}\label{ex11}
\|\nabla u\|_{L^{q}(\Omega)} \leq C\big(\|\nabla ^2 u\|_{L^p(\Omega
)} + \|u\|_{L^\infty (\partial \Omega)}\big)
\end{equation}
for  $u \in {V^{2,p}(\Omega ) \cap C(\overline \Omega )}$ fails for
$1<p<2$ and for every $q \geq 1$. The relevant set $\Omega$ is
represented  in Figure 2, Section \ref{over}. \\ Let $u : \Omega \to
\R$ be a function
 such   that $u \in C^2(\Omega)\cap C^0(\overline \Omega)$, $u (x, y) =
1 + \tfrac 1{b_k}(y- 1- c_k)$ if $(x,y) \in R_k$, $u(x,y)=0$ if
$(x,y) \in R$, and $u(x,y)$ depends only on $y$ in $N_k$.  One has
that
\begin{equation}\label{ex12} \|u\|_{L^\infty(\partial \Omega )} = 2,
\end{equation}
and
\begin{equation}\label{ex13}\|\nabla u\|_{L^1(\Omega)} \geq \sum _{k=1}^\infty \|\nabla
u\|_{L^1(R_k)} = \sum _{k=1}^\infty \frac{a_kb_k}{b_k} = \sum
_{k=1}^\infty a_k = \infty.
\end{equation}
 On the other hand,
\begin{equation*}
\|\nabla^2 u\|_{L^p(\Omega)} = \Big(\sum _{k=1}^\infty \|\nabla^2
u\|_{L^p(N_k)}^p\Big)^{\frac 1p}.
\end{equation*}
 Thus,
\begin{equation}\label{ex14}
\|\nabla^2 u\|_{L^p(\Omega)} < \infty,
\end{equation}
provided that the sequence $c_k$ decays sufficiently fast to $0$.
Equations \eqref{ex12}--\eqref{ex14} tell us that inequality
\eqref{ex11} cannot hold in $\Omega$. }
\end{example}

\begin{example}\label{ex2}
{\rm We are concerned here with the sharpness of the exponent $q$
given by \eqref{q} in inequality \eqref{mainsecond1}. An open set
set $\Omega \subset \rn$ is produced where inequality
\eqref{mainsecond1} fails if $q$ exceeds the right-hand side of
\eqref{q}. Consider the domain $\Omega \subset \rn$, with $n \ge 3$,
depicted for $n=3$ in Figure 3, Section \ref{over}.
By the standard Sobolev inequality, one necessarily has $q \leq
\frac{np}{n-p}$. Thus, it suffices to show that
\begin{equation}\label{ex24}
q  \leq \frac{rn}{n-1}.
\end{equation}
Let $\{u_k\}$ be a sequence of functions $u_k: \Omega \to \R$
enjoying the following properties:  $u_k \in C^2(\Omega)\cap
C^0(\overline \Omega)$;
 $u_ k(x', x_n) =
1 + \tfrac 1{b_k}(x_n- 1- c_k)$ if $(x',x_n) \in R_k$;  $u_k (x',
x_n)$ depends only on $x_n$ on $N_k$;  $u_k(x',x_n)=0$ if $(x',x_n)
\notin R_k\cup N_k$. One has that, for $k \in \N$,
$$\|\nabla u_k\|_{L^q(\Omega)} \geq  \|\nabla
u_k\|_{L^q(R_k)} =  b_k^{\frac {n-q}q},$$
$$\|u_k\|_{\mathcal V^{1,0}L^r(\partial \Omega )} \leq C b_k^{\frac{n-1-r}r},$$
and
$$\|\nabla^2 u_k\|_{L^p(\Omega)} =\|\nabla^2
u_k\|_{L^p(N_k)},$$ for some constant $C$. Thus, inequality
\eqref{mainsecond1} entails that
\begin{equation}\label{ex25}
b_k^{\frac {n-q}q} \leq C\big(\|\nabla^2 u_k\|_{L^p(N_k)} +
b_k^{\frac{n-1-r}r}\big)
\end{equation}
 for some constant $C$, and for every $k \in
\N$. The norm on the right-hand side of \eqref{ex25} decays to $0$
arbitrarily fast, provided that $d_k$ tends to $0$ fast enough.
Hence, if \eqref{mainsecond1} holds, then $q$ must necessarily
satisfy \eqref{ex24}. }
\end{example}

\begin{example}\label{ex4}
{\rm We conclude by showing that the number $\big[\tfrac
{m-1}2\big]$ of derivatives to be prescribed on $\partial \Omega$,
appearing in our inequalities, is minimal, in general, for an $m$-th
order Sobolev inequality to hold in an arbitrary  domain $\Omega$.
This will be demonstrated by two examples.
 \par First, given $p>1$
and $h, i , n \in \N$ such that $p(m - h)< n$ and $0 \leq h \leq i <
\tfrac m 2$, we produce a counterexample to the inequality
\begin{equation}\label{ex41}
\|\nabla ^h u\|_{L^{\frac{pn}{n-p(m -h)}}(\Omega )} \leq C \|\nabla
^m u\|_{L^{p}(\Omega )}
\end{equation}
for all $u \in V^{m , p}(\Omega) \cap C^{i-1}(\overline \Omega)$
such that $u=\nabla u = \dots = \nabla ^{i-1}u =0$ on $\partial
\Omega$. Note that the condition $i< \frac m2$ is equivalent to
$i-1<\big[\tfrac {m-1}2\big]$.
\\
Second, in the case when $p(m - h)> n> p\max\{m - i, 2i-h\}$ and $0
\leq h <i < \tfrac m  2$ we produce a counterexample to the
inequality
\begin{equation}\label{ex42}
\|\nabla ^h u\|_{L^{\infty}(\Omega )} \leq C \|\nabla ^m
u\|_{L^{p}(\Omega )}
\end{equation}
for all $u \in V^{m , p}(\Omega) \cap C^{i-1}(\overline \Omega)$
such that $u=\nabla u = \dots = \nabla ^{i-1}u =0$ on $\partial
\Omega$.
\par To this purposes,  consider a domain $\Omega$ similar to the
one constructed in Example \ref{ex3}, save that the sequence of
cusp-shaped subdomains  $\Omega _k$ is  replaced with a sequence of
balls $B_{\delta _k} (x_k)$,  with radius $\delta _k$ to be chosen
later, again connected by thin cylinders (Figure 4).
\begin{figure}[ht]
\begin{center}
\includegraphics[height=10cm]{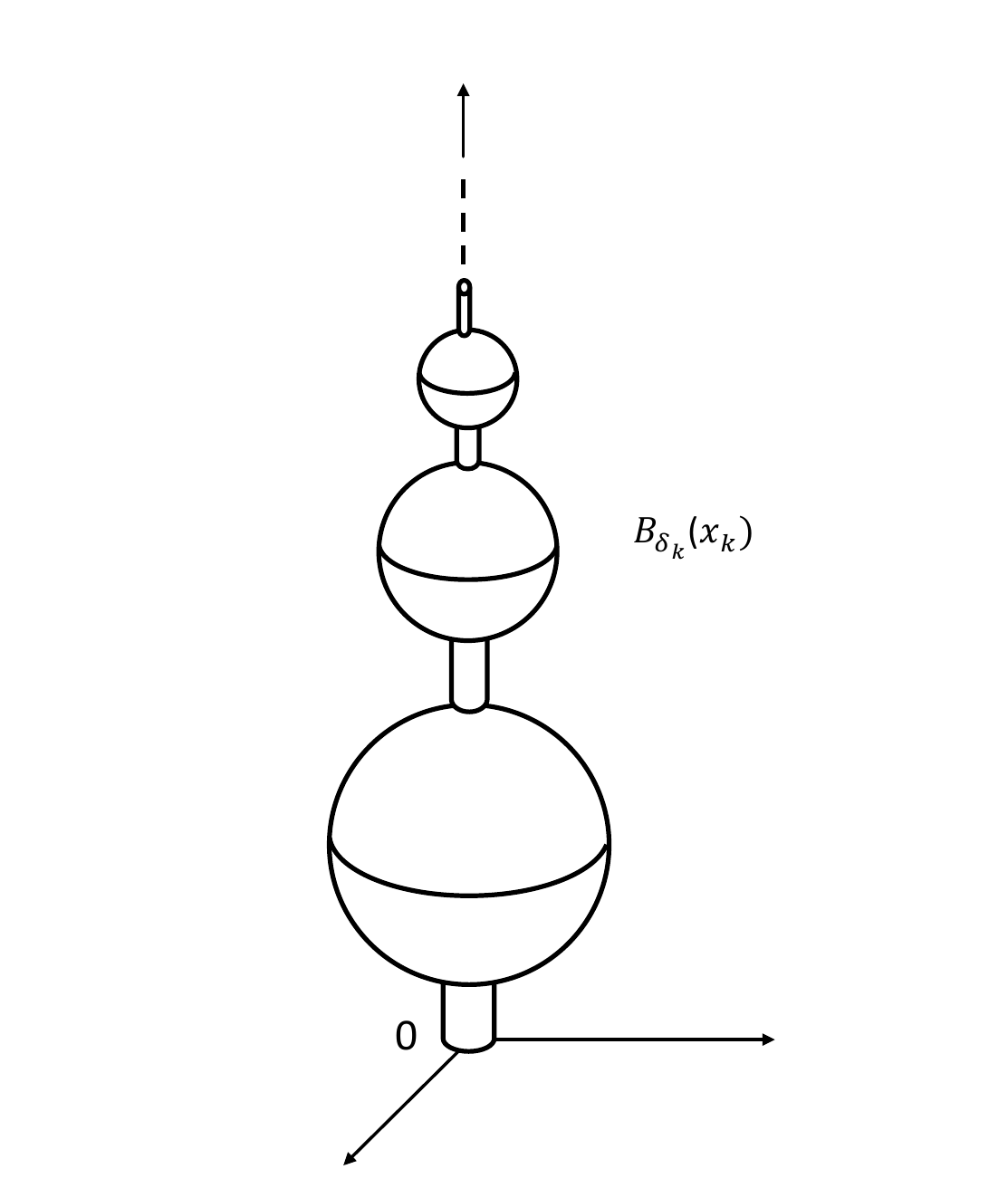}
\end{center}
        \label{Fig4}
        \caption{Example \ref{ex4}, Section \ref{sharp}}
\end{figure}
\par\noindent Let $v_k : B_{\delta _k} (x_k) \to [0, 1]$  be
the function defined as
\begin{equation}\label{ex43}
v_k (x) = \big( 1 - \tfrac {|x-x_k|^2}{\delta _k ^2}\big)^i \quad
\hbox{for $x \in B_{\delta _k} (x_k)$.}
\end{equation}
We have that $\nabla ^\ell v_k =0$ on $\partial B_{\delta _k} (x_k)$
for $0 \leq \ell \leq i-1$, and hence, given $\ell \leq i$ and
$\varepsilon _k \in (0, \tfrac {\delta _k} 2)$, there exists a
positive constant $c$ such that
\begin{equation}\label{ex44}
|\nabla ^\ell v_k| \leq c {\delta _k} ^{-i} \varepsilon _k ^{i-\ell}
\end{equation}
in an $\varepsilon$ neighborhood of $\partial B_{\delta _k} (x_k)$.
Moreover, if $\ell \leq 2i$, then there exists a positive constant
$c$ such that
\begin{equation}\label{ex45}
|\nabla ^\ell v_k| \geq c {\delta _k} ^{-\ell}
\end{equation}
in a subset of $B_{\delta _k} (x_k)$ of Lebesgue measure $\approx
{\delta _k} ^n$, whereas, if $\ell > 2i$, then
\begin{equation}\label{ex46}
\nabla ^\ell v_k =0.
\end{equation}
Next, denote by $y_k$ and $z_k$ the north and the south pole of
$B_{\delta _k} (x_k)$, respectively, and let  $\rho : [0, \infty )
\to [0, 1]$ be a smooth function, which vanishes in $[0,\tfrac 12]$
and equals $1$ in $[1, \infty)$. Let us define the function $w_k :
\Omega \to [0, \infty )$ as
\begin{equation}\label{ex47}
w_k (x) = v_{k} (x) \rho (|x-y_k|/\varepsilon_k ) \rho
(|x-z_k|/\varepsilon_k ) \quad \hbox{for $x \in B_{\delta
_k}(x_k)$,}
\end{equation}
and $w_k =0$ elsewhere, where $\varepsilon _k$ will be chosen later.
\\
If $j \leq 2i$,
\begin{equation}\label{ex48}
|\nabla ^j w_k| \geq c {\delta _k} ^{-j}
\end{equation}
in a subset of $B_{\delta _k} (x_k)$ of Lebesgue measure $\approx
{\delta _k} ^n$, and, if $j>2i$,
\begin{equation}\label{ex49}
\nabla ^j v_k =0 \quad \hbox{in $ B_{\delta _k} (x_k) \setminus
(B_{\varepsilon_k }(y_k) \cup B_{\varepsilon _k}(z_k))$.}
\end{equation}
Thus,  there exists a constant $c$ such that
\begin{equation}\label{ex410}
|\nabla ^j w_k| \leq c \sum _{\ell=0}^{\min \{j, 2i\}} \varepsilon_k
^{\ell-j}|\nabla ^\ell v_k| \quad \hbox{in $B_{\varepsilon_k} (y_k)
\cup B_{\varepsilon _k}(z_k)$.}
\end{equation}
Consequently, if $j \leq m$, then
\begin{equation}\label{ex411}
|\nabla ^j w_k| \leq c \varepsilon_k ^{i-j}{\delta _k} ^{-i} \quad
\hbox{in $B_{\varepsilon_k} (y_k) \cup B_{\varepsilon_k} (z_k)$,}
\end{equation}
for some constant $c$. Hence, if $j > 2i$, then
\begin{equation}\label{ex412}
\|\nabla ^j w_k\|_{L^p(B_{\delta _k} (x_k))}^p \leq c {\delta _k}
^{-pi} \varepsilon _k^{p(i-j)+n}
\end{equation}
for some constant $c$. On the other hand, if $j \leq 2i$, then
\begin{equation}\label{ex413}
\|\nabla ^j w_k\|_{L^p(B_{\delta _k} (x_k))}^p \leq c \big({\delta
_k} ^{n-pj} + {\delta _k} ^{-pi}\varepsilon _k^{p(i-j)+n}\big)
\end{equation}
for some constant $c$.
\par\noindent Set $\varepsilon _k = \delta _k^\alpha$, with $\alpha$ to be
chosen later. Then, by  \eqref{ex412} and \eqref{ex413},
\begin{align}\label{ex418}
\| w_k\|_{V^{m ,p}(B_{\delta _k} (x_k))}^p & \leq c \big({\delta _k}
^{-pi +\alpha [p(i-m)+n]} + {\delta _k} ^{n-2pi} + {\delta _k} ^{-pi
+\alpha (-pi+n)}\big)\\ \nonumber & \leq c \big({\delta _k} ^{-pi
+\alpha [p(i-m)+n]} + {\delta _k} ^{n-2pi}\big),
\end{align}
where the last inequality holds since $p(i-m)+n < -pi+n$, owing to
the assumption that $i < \tfrac m 2$.
\par Let us  consider \eqref{ex41}. Note that
\begin{equation}\label{ex415}
n >pi.
\end{equation}
Indeed, since we are now assuming  that $p(m - h)< n$ and $0 \leq h
\leq i < \tfrac m 2$, we have that $n > p(m - h) \geq p (m - i) > p
(2i -i ) = pi$, namely \eqref{ex415}.
\\ We may choose $\alpha$ such that
\begin{equation}\label{ex419}
-pi +\alpha [p(i-m)+n] < n-2pi.
\end{equation}
Actually, if $p(i-m)+n >0$, then inequality \eqref{ex419} holds
provided that
\begin{equation}\label{ex419bis}
1 <\alpha < \tfrac {n-pi}{p(i-m)+n}.
\end{equation}
 Note that the two inequalities in
\eqref{ex419bis} are compatible since $i < \tfrac m 2$. If, instead,
$p(i-m)+n \leq 0$, then any choice of $\alpha >1$ is admissible,
since $-pi +\alpha [p(i-m)+n] \leq -pi$, whence \eqref{ex419}
follows, owing to \eqref{ex415}.
\\ By \eqref{ex418} and \eqref{ex419},
\begin{equation}\label{ex420}
\| w_k\|_{V^{m ,p}(B_{\delta _k} (x_k))}^p \leq c {\delta _k} ^{-pi
+\alpha [p(i-m)+n]},
\end{equation}
for some constant $c$. Given a sequence $\{\lambda _k\}$, define
$w_k$ as
\begin{equation}\label{ex421}
u(x) = \sum _{k=1}^\infty \lambda _k w_k(x) \quad \hbox{for $x \in
\Omega$.}
\end{equation}
Note that
 $u=\nabla u = \dots = \nabla ^{i-1}u =0$ on $\partial \Omega$.
\\ Set $q=\frac{pn}{n-p(m -h)}$, and choose
$ \lambda _k = \delta _k ^{h - \frac nq}$ for $k \in \N$.
 By \eqref{ex48}, there exists a positive constant $c$ such that
\begin{equation}\label{ex422}
\|\nabla ^hu\|_{L^q(\Omega)}^q = \sum _{k=1}^\infty \lambda _k^q
\int _{B_{\delta _k}(x_k)}|\nabla ^h w_k|^q dx \geq c \sum
_{k=1}^\infty \lambda _k^q \delta _k^{n-hq} = c \sum _{k=1}^\infty 1
= \infty.
\end{equation}
On the other hand, by \eqref{ex420} there exists a constant $c$ such
that
\begin{equation}\label{ex424}
\| u\|_{V^{m ,p}(\Omega)}^p = \sum _{k=1}^\infty \lambda _k^p \|
w_k\|_{V^{m ,p}(B_{\delta _k} (x_k))}^p \leq c \sum _{k=1}^\infty
\lambda _k^p \delta_k ^{-pi +\alpha [p(i-m)+n]} = c \sum
_{k=1}^\infty \delta_k ^{(\alpha -1)[p(i-m)+n]}.
\end{equation}
Our assumptions ensure that $p(i-m)+n \geq p(h - m ) + n >0$. Thus,
$(\alpha -1)[p(i-m)+n] >0$, and hence the last series in
\eqref{ex424} converges, provided that $\delta _k$ decays to $0$
sufficiently fast. Clearly, equations \eqref{ex422} and
\eqref{ex424} contradict \eqref{ex41}.
  \par Let us next focus  on
\eqref{ex42}. Consider again the function $u$ given by
\eqref{ex421}. Fix $\sigma \in (0, h)$, and choose $\lambda _k =
\delta _k ^{h - \sigma}$. By \eqref{ex48},
\begin{equation}\label{ex425}
\|\nabla ^hu\|_{L^\infty(\Omega)} \geq c \lim _{k \to \infty}\lambda
_k \delta _k^{-h} = c \lim _{k \to \infty} \delta
_k^{-\sigma}=\infty.
\end{equation}
Moreover, by  \eqref{ex418},
\begin{align}\label{ex426}
\| u\|_{V^{m ,p}(\Omega)}^p   = \sum _{k=1}^\infty \lambda _k^p \|
w_k\|_{V^{m ,p}(B_{\delta _k} (x_k))}^p & \leq c \sum _{k=1}^\infty
\delta _k ^{p(h - \sigma)} \big(\delta_k ^{-pi +\alpha [p(i-m)+n]}+
\delta _k ^{n-2pi}\big)  \\ \nonumber & = c \sum _{k=1}^\infty \big(
\delta_k ^{hp - \sigma p -pi +\alpha [p(i-m)+n]}+ \delta_k ^{hp -
\sigma p + n-2pi} \big).
\end{align}
The assumption  $n> p\max\{m - i, 2i-h\}$ ensures that
$$hp  + n-2pi > 0, \quad \hbox{and} \quad hp  -pi +\alpha
[p(i-m)+n]>0,$$ provided that $\alpha $ is sufficiently large. Since
$\sigma$ can be chosen arbitrarily small, we may assume that both
exponents of $\delta _k$ in the last series of \eqref{ex426} are
positive, and hence that
\begin{equation}\label{ex427}
\| u\|_{V^{m ,p}(\Omega)} < \infty,
\end{equation}
provided that $\delta _k$ decays to $0$ fast enough. Equations
\eqref{ex425} and \eqref{ex427} contradict \eqref{ex42}.
 }
\end{example}

\section*{Appendix}

A  result in the theory of Sobolev functions tells us that, if $u$
is any weakly differentiable function in $\rn$, then
\begin{equation}\label{wg}
|u(x) -u(y)| \leq C(M(|\nabla u|)(x) + M(|\nabla u|)(y)) \quad
\hbox{for a.e. $x, y \in \rn$,}
\end{equation}
for some constant $C$. Here, $M$ denotes the maximal function
operator
 defined,  for $f \in L^1_{\rm loc} (\rn)$, as
 $$Mf (x) = \sup _{B\ni x} \frac 1{\mathcal L^n (B)} \int _B
 |f(y)|\, dy \qquad \hbox{for $x \in \rn$,}$$
where $B$ denotes a ball in $\rn$. Thus, $M(|\nabla u|)$ is an upper
gradient for $u$ in the sense of metric measure spaces, as defined
in \cite{Hajlasz}.
\\
The following proposition provides us with a higher-order
counterpart of \eqref{wg}, and gives grounds for definitions
\eqref{M} and \eqref{M1}.

\begin{proposition}\label{maximaln}
Let $n, \ell \in \N$.  Then there exists a constant $C=C(\ell,n)$
such that, if $u \in W^{2\ell-1, 1}_{\rm loc}(\rn)$, then
\begin{multline}\label{maxn1}
 \left|\sum _{|\alpha|\leq \ell-1} \frac{(2\ell-2-|\alpha|)! }{
(\ell-1-|\alpha|)!\alpha !} \frac{(y - x)^\alpha}{|y - x|^{2\ell-1}}
\big[ (-1)^{|\alpha |}D^\alpha u(y) - D ^\alpha u(x)\big]\right|
\\ \leq C \big(M(|\nabla ^{2\ell-1} u|) (x ) + M(|\nabla ^{2\ell-1}
u|) (y )\big) \quad \hbox{for a.e. $x, y \in \rn$}.
\end{multline}
\end{proposition}
%
%
%
%
%
%
%
%
%
%
%

\par\noindent {\bf Proof}. By \cite[Proposition 5.1]{Boj},   if $0\leq |\alpha| \leq
2\ell-2$, then there exists a measurable function $R_{2\ell-1 ,
\alpha }(u) : \rn \times \rn \to \R$ and a constant $C=C(\ell,n)$
such that
\begin{align}\label{maxn2}
D^\alpha u(y ) = \sum _{|\gamma |\leq 2\ell-2-|\alpha | } \frac {(y
- x )^\gamma }{\gamma !} D^{\alpha + \gamma}u(x )  + R_{2\ell-1 ,
\alpha }(u)(x , y) \quad \hbox{for a.e. $x, y \in \rn$,}
\end{align}
and
\begin{align}\label{maxn3}
|R_{2\ell-1 , \alpha }(u)(x , y)| \leq C |x - y |^{2\ell-1-
|\alpha|} \big[M(|\nabla ^{2\ell-1} u|) (x ) + M(|\nabla ^{2\ell-1}
u|) (y )\big] \quad \hbox{for a.e. $x, y \in \rn$}.
\end{align}
We claim that there exist constants $C(\alpha , \ell)$, for $|\alpha
| \leq \ell -1$, such that
\begin{multline}\label{maxn4}
\sum _{|\alpha|\leq \ell-1} \frac{(2\ell-2-|\alpha|)!
}{(\ell-1-|\alpha|)!\alpha !} \frac{(y - x)^\alpha[D ^\alpha u(x) +
(-1)^{|\alpha |+1}D^\alpha u(y)]}{|y - x|^{2\ell-1}}
\\ = \sum _{|\alpha|\leq \ell-1} \frac{C(\alpha , \ell)(y - x)^\alpha}{|y -
x|^{2\ell-1}}R_{2\ell-1 , \alpha }(u)(x , y) \quad \hbox{for a.e.
$x, y \in \rn$}.
\end{multline}
Inequality \eqref{maxn1} will then follow from \eqref{maxn4} and
\eqref{maxn3}.
\\ Let us establish \eqref{maxn4}. By \eqref{maxn2}, after
exchanging the order of summation and relabeling the indices, one
obtains that there exist constants $A(\alpha , \ell)$ and $B(\alpha
, \ell)$, for $|\alpha | \leq \ell -1$, such that
\begin{align}\label{maxn20}
\sum _{|\alpha|\leq \ell-1} & \frac{(2\ell-2-|\alpha|)! }{
(\ell-1-|\alpha|)!\alpha !} (y - x)^\alpha \frac{[D ^\alpha u(x) +
(-1)^{|\alpha |+1}D^\alpha u(y)]}{|y - x|^{2\ell-1}}
\\ \nonumber & = \sum _{|\alpha|\leq \ell-1} \frac{(2\ell-2-|\alpha|)! }{
(\ell-1-|\alpha|)!\alpha !} \frac{(y - x)^\alpha}{|y - x|^{2\ell-1}}
\\ \nonumber & \times \bigg[D ^\alpha u(x)
 + (-1)^{|\alpha
|+1} \bigg(\sum _{|\gamma |\leq 2\ell-2-|\alpha | } \frac {(y - x
)^\gamma}{\gamma !} D^{\alpha + \gamma}u(x )  + R_{2\ell-1 , \alpha
}(u)(x , y)\bigg)\bigg]
\\ \nonumber & =
\sum _{|\alpha | \leq 2\ell-2} \frac{(y - x)^\alpha}{|y -
x|^{2\ell-1}} A(\alpha , \ell) D^\alpha (u)(x) + \sum _{|\alpha |
\leq \ell-1} \frac{(y - x)^\alpha}{|y - x|^{2\ell-1}} B(\alpha ,
\ell)R_{2\ell-1 , \alpha }(u)(x , y) \end{align}
 for a.e. $x, y \in \rn$.
 \\ Now, let us choose $u=\mathcal P$ in \eqref{maxn20}, where
 $\mathcal P$ is a polynomial
 of the form
\begin{equation}\label{maxn11}
\mathcal P(y)= \sum _{|\alpha| \leq 2\ell-2} b_\alpha (y - x)
^\alpha \quad \hbox{for $y \in \rn$,}
\end{equation}
with $b_\alpha \in \R$. Clearly,
\begin{equation}\label{maxn22}
\hbox{$D^\alpha \mathcal P(x) = \alpha ! b_\alpha $ if $|\alpha |
\leq 2\ell-2$, \qquad $D^\alpha \mathcal P(x)=0$ if $|\alpha| >
2\ell-2$.}
\end{equation}
Hence, $R_{2\ell-1 , \alpha }(\mathcal P)(x , y) =0$, and from
\eqref{maxn20} we obtain that
\begin{multline}\label{maxn23}
\sum _{|\alpha|\leq \ell-1}  \frac{(2\ell-2-|\alpha|)! }{
(\ell-1-|\alpha|)!\alpha !} (y - x)^\alpha \frac{[ D ^\alpha
\mathcal P(x) + (-1)^{|\alpha |+1}D^\alpha \mathcal P(y)]}{|y -
x|^{2\ell-1}}
 =
\sum _{|\alpha | \leq 2\ell-2} \frac{(y - x)^\alpha}{|y -
x|^{2\ell-1}} A(\alpha , \ell) \alpha ! b_\alpha
\end{multline}
for a.e. $x, y \in \rn$.
%
%
%
%
We next express the leftmost side of \eqref{maxn20} in an
alternative form. Define $\varphi : \R \to \R$ as
$$\varphi (t) = u(x + t \vartheta ) \quad \hbox{for $t \in \R$,}$$
where $$\vartheta = \frac{y - x}{|y - x|}.$$  Given $j \in \{1,
\dots, 2\ell-1\}$, we have that
\begin{align}\label{maxn7}
\varphi ^{(j)}(t) & = \sum _{|\alpha |=j} \frac{j!}{\alpha !}
\vartheta ^\alpha D^\alpha u(x + t \vartheta )
 \quad \hbox{for a.e. $t
\in \R$.}
\end{align}
Thus, by the Taylor formula centered at $t=0$, if $k \in \{1, \dots,
2\ell-1\}$, then
\begin{align}\label{maxn8}
\varphi ^{(k)}(|y - x|) & = \sum _{j=k}^{2\ell-2} \frac{\varphi
^{(j)} (0)}{(j-k)!}|y - x|^{j-k} + Q_{2\ell-1, k}(\varphi )(|y - x|)
\\ \nonumber & =
 \sum _{j=k}^{2\ell-2} \frac{1} {(j-k)!} \sum _{|\alpha
|=j} \frac{j!}{\alpha !} |y - x|^{j-k}\vartheta ^\alpha D^\alpha u
(x ) + Q_{2\ell-1, k}(\varphi )(|y - x|)
\end{align}
for a.e. $x, y \in \rn$, where $Q_{2\ell-1, k}(\varphi )$ denotes
the remainder in the $(2\ell-2-k)$-th order Taylor formula for
$\varphi ^{(k)}$, centered at $t=0$.  By \eqref{maxn7} and
\eqref{maxn8}, there exist constants $A'(\alpha , \ell)$ and
$B'(\alpha , \ell)$ such that
\begin{align}\label{maxn9}
& \sum _{|\alpha|\leq \ell-1}  \frac{(2\ell-2-|\alpha|)! }{
(\ell-1-|\alpha|)!\alpha !}(y - x)^\alpha \frac{[D ^\alpha u(x)+
(-1)^{|\alpha |+1}D^\alpha u(y)]}{|y - x|^{2\ell-1}}
\\ \nonumber & =
(-1)^\ell \sum _{k=0}^{\ell-1} \frac{(2\ell-k-2)!}{k!
(\ell-k-1)!}\frac{\big[\varphi ^{(k)}(0) + (-1)^{k+1}\varphi
^{(k)}(|y - x|)\big]}{|y - x|^{2\ell-k-1}}
 \\ \nonumber & = (-1)^\ell
\sum _{k=0}^{\ell-1} \frac{(2\ell-k-2)!}{k! (\ell-k-1)!} \frac 1 {|y
- x|^{2\ell-k-1}}
\\ \nonumber & \times
\bigg[ \sum _{|\alpha |=k} \frac{k!}{\alpha !} \vartheta ^\alpha
D^\alpha u(x)
 + (-1)^{k+1} \bigg(\sum _{j=k}^{2\ell-2} \frac{1}
{(j-k)!} |y - x|^{j-k}\sum _{|\alpha |=j} \frac{j!}{\alpha !}
\vartheta ^\alpha D^\alpha u(x ) + Q_{2\ell-1, k}(\varphi )(|y -
x|)\bigg)\bigg]
\\ \nonumber & =
\sum _{|\alpha | \leq 2\ell-2} \frac{(y - x)^\alpha}{|y -
x|^{2\ell-1}} A'(\alpha , \ell) D^\alpha u(x) + \sum _{|\alpha |
\leq \ell-1} \frac{1}{|y - x|^{2\ell-|\alpha|-1}} B'(\alpha ,
\ell)Q_{2\ell-1 , |\alpha| }(\varphi)(|x - y|)
\end{align}
for a.e. $x, y \in \rn$. If $\mathcal P$ is again a polynomial in
$t$ of the form \eqref{maxn11}, then $\varphi$ is also a polynomial
of degree not
 exceeding $2\ell-2$ , and
from  \eqref{maxn22} and \eqref{maxn9}, applied with $u = \mathcal
P$, we obtain that
\begin{multline}\label{maxn25}
\sum _{|\alpha|\leq \ell-1}  \frac{(2\ell-2-|\alpha|)! }{
(\ell-1-|\alpha|)!\alpha !} (y - x)^\alpha \frac{[D ^\alpha \mathcal
P(x)+ (-1)^{|\alpha |+1}D^\alpha \mathcal P(y)]}{|y - x|^{2\ell-1}}
= \sum _{|\alpha | \leq 2\ell-2} \frac{(y - x)^\alpha}{|y -
x|^{2\ell-1}} A'(\alpha , \ell)\alpha ! b_\alpha
\end{multline} for a.e.
$x, y \in \rn$. Owing to the arbitrariness of the coefficients
$b_\alpha$, we infer from \eqref{maxn23} and \eqref{maxn25} that
\begin{equation}\label{sa1}
A(\alpha , \ell) = A'(\alpha , \ell)
\end{equation}
for every multi-index $\alpha$ such that $|\alpha | \leq 2\ell-2$.
On the other hand, by \eqref{1dim3} and \eqref{1dim4}, applied with
$ \varsigma = \psi= \varphi$, $a=0$ and $b=|y -x|$, and by
\eqref{maxn9} and \eqref{maxn25},
%
%
\begin{align}\label{maxn12}
0 & = \varphi ^{(2\ell-1)}(t)  =
\bigg[\frac{d^{\ell-1}}{dt^{\ell-1}}\bigg(\frac{\varphi (t)}{(t-|y -
x|)^\ell}\bigg)_{|t=0} +
\frac{d^{\ell-1}}{dt^{\ell-1}}\bigg(\frac{\varphi
(t)}{t^\ell}\bigg)_{|t=|y - x|}\bigg]
\\ \nonumber &
= (-1)^\ell \sum _{k=0}^{\ell-1} \frac{(2\ell-k-2)!}{k!
(\ell-k-1)!}\frac{\big[\varphi ^{(k)}(0) + (-1)^{k+1}\varphi
^{(k)}(|y - x|)\big]}{|y - x|^{2\ell-k-1}}
\\ \nonumber &
= \sum _{|\alpha|\leq \ell-1} \frac{(2\ell-2-|\alpha|)! }{
(\ell-1-|\alpha|)!\alpha !} \frac{(y - x)^\alpha[(-1)^{|\alpha
|+1}D^\alpha u(y) + D ^\alpha u(x)]}{|y - x|^{2\ell-1}}
\\ \nonumber & =
\sum _{|\alpha | \leq 2\ell-2} \frac{(y - x)^\alpha}{|y -
x|^{2\ell-1}} A'(\alpha , \ell)\alpha ! b_\alpha
\end{align} for a.e.
$x, y \in \rn$. By the arbitrariness of the coefficients $b_\alpha$
again, $A'(\alpha ,\ell)=0$ for every $\alpha$ such that $|\alpha |
\leq 2\ell-2$. Hence, owing to \eqref{sa1},
\begin{equation}\label{sa2}
A(\alpha , \ell) = 0
\end{equation}
for every  $\alpha$ such that $|\alpha | \leq 2\ell-2$. Equations
\eqref{maxn20} and \eqref{sa2}  tell us that
\begin{align}\label{maxn28}
\sum _{|\alpha|\leq \ell-1} & \frac{(2\ell-2-|\alpha|)! }{
(\ell-1-|\alpha|)!\alpha !} (y - x)^\alpha \frac{ D ^\alpha u(x) +
(-1)^{|\alpha |+1}D^\alpha u(y)]}{|y - x|^{2\ell-1}}
\\ \nonumber & = \sum _{|\alpha|\leq \ell-1}  \frac{(y - x)^\alpha}{|y -
x|^{2\ell-1}} B(\alpha , \ell) R_{2\ell-1 , \alpha }(u)(x , y) \quad
\hbox{for a.e. $x, y \in \rn$,}
\end{align}
whence \eqref{maxn4} follows with $C(\alpha ,\ell)=B(\alpha ,\ell)$.
\qed


\begin{thebibliography}{99}


%

\bibitem[Ad1]{Adams1}
D.R.Adams, Traces of potentials arising from translation invariant
operators, \emph{Ann. Sc. Norm. Super. Pisa} {\bf 25} (1971),
203--217.

\bibitem[Ad2]{Adams2}
D.R.Adams, A trace inequality  for generalized potentials,
\emph{Studia Math.} {\bf 48} (1973), 99--105.



\bibitem[AFT]{AFT}
A.Alvino, V.Ferone \& G.Trombetti, Moser-type inequalities in
Lorentz spaces, \emph{Potential Anal.} {\bf 5} (1996), 273--299.



\bibitem[AFP]{AFP} L.Ambrosio, N.Fusco \& D.Pallara,  Functions of bounded
variation and free discontinuity problems,  Oxford University Press,
Oxford, 2000.



\bibitem[AT]{AT} L.Ambrosio \& P.Tilli, \emph{Topics on Analysis in Metric Spaces},  Oxford
University Press, Oxford, 2004.



%
%



\bibitem[Au]{Aubin}
T.Aubin, Probl\`emes isop\'erimetriques et espaces de Sobolev,
\emph{J. Diff. Geom.} {\bf 11} (1976), 573--598.




%
%

\bibitem[BCR]{BCR1}
F.Barthe, P.Cattiaux \& C.Roberto,  Interpolated inequalities
between exponential and Gaussian, Orlicz hypercontractivity and
isoperimetry, \emph{Rev. Mat. Iberoam.} {\bf 22} (2006), 993--1067.

%
%


\bibitem[BWW]{BWW}
T.Bartsch, T.Weth \& M.Willem, A Sobolev inequality with remainder
term  and critical equations on domains with topology for the
polyharmonic operator, \emph{Calc. Var. Partial Differential
Equations} {\bf 18} (2003), 253--268.



%
\bibitem[BS]{BS} C.Bennett
\&  R.Sharpley, \emph{Interpolation of operators}, Academic Press,
Boston, 1988.






\bibitem[BB]{BB} A.Bj\"orn \& J.Bj\"orn, \emph{Nonlinear potential theory on metric
spaces}, European Mathematical Society (EMS), Z\"urich, 2011.


%
%
%
%


\bibitem[BH2]{BH}
S.G.Bobkov \& C.Houdr\'e,  Some connections between isoperimetric
and Sobolev-type inequalities, \emph{Mem. Am. Math. Soc.} {\bf 25}
(1997), viii+111.


\bibitem[BL]{BLbis}
S.G.Bobkov \& M.Ledoux, From Brunn-Minkowski to sharp Sobolev
inequalities, \emph{Ann. Mat. Pura Appl.} {\bf 187} (2008),
389--384.










%
%


\bibitem[Bo]{Boj} B.Bojarski, Pointwise characterization of Sobolev classes, \emph{Tr. Mat. Inst. Steklova} {\bf 255} (2006),
  71--87; English translation in \emph{Proc. Steklov Inst. Math.} {\bf 255} (2006),
  65--81.

\bibitem[Bo]{boggio}
T.Boggio, Sulle funzioni di Green d'ordine $m$, \emph{Rend. Circ.
Mat. Palermo} {\bf 20} (1905), 97--135 (Italian).




%
%
%
%




\bibitem[BL]{BL}
H.Br\'ezis \& E.Lieb, Sobolev inequalities with remainder terms,
\emph{J. Funct. Anal.} {\bf 62} (1985),  73--86.


\bibitem[BK]{BK}
S.Buckley \& P.Koskela, Sobolev-Poincar\'{e} implies John,
\emph{Math. Res. Lett.} {\bf 2} (1995), 577--593.


\bibitem[BK1]{BK1}
S.Buckley \& P.Koskela, Criteria for embeddings of
Sobolev-Poincar\'{e} type, \emph{Int. Math. Res. Not.} {\bf 18}
(1996), 881--902.





\bibitem[BZ]{BZ} Yu.D.Burago \& V.A.Zalgaller, \emph{Geometric inequalities}, Springer,
Berlin, 1988.







%
%
%


\bibitem[CDPT]{CDPT}
L.Capogna, D.Danielli, S.D.Pauls \& J.T.Tyson, \emph{An introduction
to the Heisenberg group and the sub-Riemannian isoperimetric
problem}, Birkhauser, Basel, 2007.












%
%


\bibitem[Cha]{chavel}
I.Chavel, \emph{Isoperimetric inequalities: differential geometric
aspects and analytic perspectives}, Cambridge University Press,
Cambridge, 2001.



\bibitem[Che]{Cheeger}
J.Cheeger, A lower bound for the smallest eigevalue of the
Laplacian, in \emph{Problems in analysis}, 195--199, Princeton Univ.
Press, Princeton, 1970.





%




\bibitem[Ci1]{Ci_ind}
A.Cianchi, A sharp embedding theorem for Orlicz--Sobolev spaces,
\emph{Indiana Univ. Math.  J.} {\bf 45} (1996), 39--65.

%







%


\bibitem[Ci2]{Ci1}
A.Cianchi, Symmetrization and second-order Sobolev inequalities,
\emph{Ann.\ Mat. Pura Appl.} {\bf 183} (2004), 45--77.


\bibitem[CFMP]{CFMP1}
A.Cianchi, N.Fusco,  F.Maggi \& A.Pratelli, The sharp Sobolev
inequality in quantitative form, \emph{J. Eur. Math. Soc.} {\bf 11}
(2009), 1105--1139.

%


\bibitem[CP]{CP_gauss}
A.Cianchi \& L.Pick, Optimal Gaussian Sobolev embeddings, \emph{J.
Funct. Anal.} {\bf 256} (2009), 3588--3642.





%
%



%






\bibitem[Da]{davis} P.J.Davis,  \emph{Interpolation and approximation}, Blaisdell Publishing Company, Waltham (Ma),
1963.

%
%
%
%




\bibitem[EKP]{EKP}
D.E.Edmunds, R.Kerman \& L.Pick, Optimal Sobolev imbeddings
involving rearrangement-invariant quasinorms, \emph{J. Funct. Anal.}
{\bf 170} (2000), 307--355.


\bibitem[EFKNT]{EFKNT}
L.Esposito, V.Ferone, B.Kawohl, C.Nitsch, \& C.Trombetti, The
longest shortest fence and sharp Poincar\'e-Sobolev inequalities,
\emph{Arch. Ration. Mech.  Anal.} {\bf 206} (2012), 821--851.



%
%


\bibitem[FHK]{FHK}
B.Franchi, P.Haj\l asz, P.Koskela, Definitions of Sobolev classes on
metric spaces, \emph{Ann. Inst. Fourier} {\bf 49} (1999),
1903--1924.


%


\bibitem[GGS]{GGS}
 F.Gazzola, H.-C.Grunau \& G.Sweers,
\emph{Polyharmonic boundary value problems. Positivity preserving
and nonlinear higher order elliptic equations in bounded domains},
 Springer-Verlag, Berlin, 2010.



%
%
\bibitem[Gr]{Gr} A.Grigor'yan,  Isoperimetric inequalities and capacities on Riemannian
manifolds, in The Maz'ya anniversary collection, Vol. 1 (Rostock,
1998), 139--153, \emph{Oper. Theory Adv. Appl.}, 109, Birkh�user,
Basel, 1999.
%
%


\bibitem[Ha]{Hajlasz} P.Hai\l asz, Sobolev spaces on an arbitrary metric space, \emph{Potential Anal.} {\bf 5} (1996), 403--415.



\bibitem[HaKo]{HK} P.Hai\l asz \& P.Koskela, Isoperimetric inequalites and imbedding theorems in irregular
domains, \emph{J. London Math. Soc.} {\bf 58} (1998), 425--450.






\bibitem[He]{Heb}
E.Hebey, \emph{Analysis on manifolds: Sobolev spaces and
inequalities}, Courant Lecture Notes in Mathematics \ {\bf 5}, AMS,
Providence, 1999.



\bibitem[Hei]{Hein} J.Heinonen, \emph{Lectures on analysis on metric
spaces},  Springer-Verlag, New York, 2001.


\bibitem[HeKo]{HeKo} J.Heinonen \& P.Koskela, Quasiconformal maps in metric spaces with controlled geometry, \emph{Acta Math.} {\bf 181} (1998), 1--61.




%
%


\bibitem[HS]{HS1}
D.Hoffman \& J.Spruck, Sobolev and isoperimetric inequalities for
Riemannian submanifolds,
 \emph{Comm. Pure Appl. Math.} {\bf 27} (1974),
715--727; A correction to: ``Sobolev and isoperimetric inequalities
for Riemannian submanifolds, Comm. Pure Appl. Math. 27 (1974),
715--725", \emph{Comm. Pure Appl. Math.} {\bf 28} (1975), 765--766.


\bibitem[Ho]{H} T.Holmstedt, Interpolation of quasi-normed
spaces, \emph{Math. Scand.} {\bf 26} (1970), 177--199.

\bibitem[KP]{KP}
R.Kerman \& L.Pick, Optimal Sobolev imbeddings,  \emph{Forum Math.}
{\bf 18} (2006), 535--570.


%
%


%
%
%



\bibitem[KM]{KM} T.Kilpel\"ainen \& J.Mal\'y,  Sobolev inequalities on sets with
irregular boundaries,  \emph{Z. Anal. Anwendungen} {\bf 19} (2000),
369--380.



\bibitem[Kl]{Kl}
V.S.Klimov, Imbedding theorems and geometric inequalities,
\emph{Izv. Akad. Nauk SSSR} {\bf 40} (1976), 645--671 (Russian);
English translation: Math. USSR Izv. {\bf 10} (1976), 615--638.


%


\bibitem[Kol]{Ko}
V.I.Kolyada, Estimates on rearrangements and embedding theorems,
\emph{Mt. Sb.} {\bf 136} (1988), 3--23 (Russian); English
translation: Math. USSR Sb. {\bf 64} (1989), 1--21.


\bibitem[Kos]{Koskela} P.Koskela, Metric Sobolev spaces, in \emph{Nonlinear analysis, function spaces and applications. Vol. 7},
132--147, Czech. Acad. Sci., Prague, 2003.



%
%
%

%


\bibitem[LPT]{LPT}
P.-L.Lions, F.Pacella \& M.Tricarico, Best constants in Sobolev
inequalities
 for functions vanishing on some part of the boundary and related questions,
 \emph{Indiana Univ. Math. J.} {\bf 37} (1988), 301--324.


\bibitem[LYZ]{LYZ}
E.Lutwak, D.Yang \& G.Zhang, Sharp affine $L_p$ Sobolev
inequalities,
 \emph{J. Diff. Geom.} {\bf 62} (2002),
17--38.










\bibitem[MV1]{MaggiVillani1}
F. Maggi \& C. Villani,
 Balls have the worst best Sobolev inequalities,
 \emph{J.  Geom. Anal.} {\bf 15} (2005),
83--121.



\bibitem[MV2]{MaggiVillani2}
F. Maggi \& C. Villani, Balls have the worst best Sobolev
inequalities. II. Variants and extensions,
 \emph{Calc. Var. Partial Differential Equations} {\bf 31}(2008),
47--74.



\bibitem[Mi]{M}
E.Milman, On the role of convexity in functional and isoperimetric
inequalities, \emph{Proc. London Math. Soc.} {\bf 99} (2009),
32--66.



\bibitem[Mo]{Moser}
J.Moser, A sharp form of an inequality by Trudinger, \emph{Indiana
Univ. Math. J.} {\bf 20} (1971), 1077--1092.


\bibitem[Mat]{Mattila} P.Mattila, \emph{Geometry of sets and measures in Euclidean spaces},
Cambridge University
Press, Cambridge, 1995.




\bibitem[Ma1]{Ma1960} V.G.Maz'ya, Classes of regions and imbedding theorems for function spaces,
{\em Dokl. Akad. Nauk. SSSR} {\bf 133} (1960), 527--530 (Russian);
English translation: {\em Soviet Math. Dokl.} {\bf 1} (1960),
882--885.
\bibitem[Ma3]{Ma1961} V.G.Maz'ya, On p-conductivity and theorems on embedding certain functional
spaces into a C-space,  \emph{Dokl. Akad. Nauk SSSR} {\bf 140}
(1961), 299--302 (Russian).



\bibitem[Ma8]{Mabook} V.G.Maz'ya,  \emph{Sobolev spaces with applications to elliptic partial differential equations}, Springer, Heidelberg, 2011.

\bibitem[MP1]{MP1} V.G.Maz'ya \& S.V.Poborchi,  \emph{Differentiable functions on bad domains}, World Scientific, Singapore, 1997.





%
%




%
%
%


%
%
%

\bibitem[Sa]{Saloff}
L.Saloff-Coste, \emph{Aspects of Sobolev-type inequalities},
Cambridge University Press, Cambridge, 2002.




%

%
%
%
%





\bibitem[Ta]{Ta}
G.Talenti, Best constant in Sobolev inequality, \emph{Ann. Mat. Pura
Appl.} {\bf 110} (1976), 353--372.

%
%
%
%
%

\bibitem[Zh]{Zh}
G.Zhang, The affine Sobolev inequality, \emph{J. Diff. Geom.} {\bf
53} (1999), 183--202.



%
%
%


%

%


\end{thebibliography}
\end{document}